\documentclass[a4paper,leqno,12pt]{article}

\usepackage{amssymb,latexsym,amsxtra,amscd,amsfonts,amsthm,amsmath,verbatim}

\def\proclaim #1. #2\par{\medbreak\indent{\bf#1.\enspace}
{\sl#2}\par\medbreak}

\begin{document}

\title{Cohomology and $L$-values}
\author {By Hiroyuki Yoshida\thanks{The author was supported
by Grant-in-Aid for Scientific Research (C) (No. 21540014),
Japan Society for Promotion of Science.}}
\date{}
\maketitle

\begin{abstract}
In a paper published in 1959, Shimura presented an elegant calculation
of the critical values of $L$-functions attached to elliptic modular forms
using the first cohomology group. We will show that a similar calculation
is possible for Hilbert modular forms over real quadratic fields using
the second cohomology group. We present explicit numerical examples
calculated by this method.
\end{abstract}

\bigskip

In a celebrated paper [Sh1] published in 1959, Shimura showed that
ratios of critical values of the $L$-function attached to an elliptic modular
form can be calculated explicitly using the cohomology group.
This method was developed into the theory of modular symbols by Manin [Man].
Though there have been great advances during the next  half century in
understanding the relationship of automorphic forms and group cohomologies,
it seems that no explicit calculations of $L$-values using cohomology groups
were performed beyond the one dimensional case. The purpose of this paper is to show
that we can use cohomology groups effectively for calculations of
$L$-values even in higher dimensional cases.

To explain our ideas and results, it is best to review first the calculation in [Sh1].
Let $\mathfrak H$ be the complex upper half plane.
Let $\Gamma$ be a Fuchsian group and let $\Omega$ be a cusp form of weight $k \geq 2$
with respect to $\Gamma$. Put $l=k-2$ and let $\rho _l$ be
the symmetric tensor representation of ${\rm GL}(2, \mathbf C)$ of degree $l$
on a vector space $V$. We regard $V$ as a $\Gamma$-module. Put $\rho =\rho _l$.
We consider a $V$-valued differential form on $\mathfrak H$:
$$
\mathfrak d(\Omega )=\Omega (z)\begin{bmatrix} z \\ 1 \end{bmatrix} ^l dz.
$$
Here $\begin{bmatrix} z \\ 1 \end{bmatrix} ^l$ denotes the column vector of dimension
$l+1$ whose components are $z^l$, $z^{l-1}, \ldots , 1$.
We have $\mathfrak d(\Omega )\circ\gamma =\rho (\gamma )\mathfrak d(\Omega )$ for every
$\gamma \in \Gamma$. Here $\mathfrak d(\Omega )\circ\gamma$ denotes the transform of
$\mathfrak d(\Omega )$ by $\gamma$. Take a point of the complex upper half plane or
a cusp of $\Gamma$ and denote it by $z_0$. For $\gamma \in \Gamma$,
we consider the integral
\begin{equation}
f(\gamma )=\int _{z_0}^{\gamma z_0} \mathfrak d(\Omega ).
\tag{1}
\end{equation}
Then $f$ satisfies the $1$-cocycle condition:
$$
f(\gamma _1\gamma _2)
=f(\gamma _1)+\rho (\gamma _1)f(\gamma _2).
$$
The cohomology class of $f$ in $H^1(\Gamma , V)$ does not depend on
the choice of $z_0$. Let $p \in \Gamma$ be a parabolic element and
$z_0^{\prime}$ be the cusp fixed by $p$. Then we have
$$
f(p)=(\rho (p)-1)\int _{z_0^{\prime}}^{z_0} \mathfrak d(\Omega ).
$$
Thus $f(p)$ looks like a coboundary,
which is the parabolic condition on $f$.

Now let $\Gamma ={\rm SL}(2, \mathbf Z)$ and $z_0=i\infty$. Put
$$
\sigma =\begin{pmatrix} 0 & 1 \\ -1 & 0 \end{pmatrix}, \qquad
\tau =\begin{pmatrix} 1 & 1 \\ 0 & 1 \end{pmatrix}.
$$
Then we find
\begin{equation}
f(\sigma\tau )=-\bigg(\int _0^{i\infty} \Omega (z)z^tdz\bigg)_{0 \leq t \leq l}
=-\bigg(i^{t+1}R(t+1, \Omega )\bigg)_{0 \leq t \leq l},
\tag{2}
\end{equation}
where $R(s,\Omega )=(2\pi )^{-s}\Gamma (s)L(s, \Omega )$ with the $L$-function
$L(s, \Omega )$ of $\Omega$. Since $(\sigma\tau )^3=1$, the $1$-cocycle condition gives
\begin{equation}
[1+\rho (\sigma\tau )+\rho ((\sigma\tau )^2)]f(\sigma\tau )=0.
\tag{3}
\end{equation}
In other words, $f(\sigma\tau )$ is annihilated by the element
$1+\sigma\tau +(\sigma\tau )^2$ of the group ring
$\mathbf Z[{\rm SL}(2, \mathbf Z)]$. This gives a constraint on
the critical values of $L(s, \Omega )$. For $k=12$ and $\Omega =\Delta$,
Shimura obtained that
$$
R(8, \Delta) =\frac {5}{4}R(6, \Delta ), \qquad R(10, \Delta )=\frac {12}{5}R(6, \Delta ),
\quad {\rm etc}.
$$

In this paper, we will treat the case of Hilbert modular forms over
a real quadratic field $F$. Let $\mathcal O_F$ be the ring of integers of $F$ and
$\Gamma$ be a congruence subgroup of ${\rm SL}(2, \mathcal O_F)$. 
Let $\Omega$ be a Hilbert modular cusp form of weight $(k_1, k_2)$ with respect to
$\Gamma$. We assume $2 \leq k_2\leq k_1$ and put $l_i=k_i-2$, $i=1$, $2$.
The first step is to attach an explicitly given $2$-cocycle of $\Gamma$ to $\Omega$.
This is given in the author's book [Y3] as follows.
Let $\rho =\rho _{l_1}\otimes\rho _{l_2}$ and $V$ be the representation space of $\rho$.
We consider a $V$-valued differential form on $\mathfrak H^2$:
$$
\mathfrak d(\Omega )=\Omega (z) \begin{bmatrix} z_1 \\ 1 \end{bmatrix} ^{l_1} \otimes
\begin{bmatrix} z_2 \\ 1 \end{bmatrix} ^{l_2} dz_1dz_2,
\qquad z=(z_1, z_2) \in \mathfrak H^2.
$$
We have
\begin{equation}
\mathfrak d(\Omega )\circ\gamma =\rho (\gamma )\mathfrak d(\Omega ),
\qquad \gamma \in \Gamma .
\tag{4}
\end{equation}
Take a point $w=(w_1, w_2)$ on $\mathfrak H^2$. For
$\gamma _1$, $\gamma _2 \in \Gamma$, we consider the integral
\begin{equation}
f(\gamma _1, \gamma _2)
=\int _{\gamma _1\gamma _2w_1}^{\gamma _1w_1}
\int _{w_2}^{\gamma _1^{\prime}w_2} \mathfrak d(\Omega ).
\tag{5}
\end{equation}
Here $\gamma _1^{\prime}$ denotes the conjugate of $\gamma_1$ by
${\rm Gal}(F/\mathbf Q)$. Then $f$ is a $2$-cocycle of $\Gamma$
taking values in $V$. The cohomology class of $f \in H^2(\Gamma , V)$
does not depend on the choice of $w$. Let $p \in \Gamma$ be
a parabolic element and let $(w_1^{\ast}, w_2^{\ast})$ be the cusp fixed by $p$.
Since $\Omega$ is a cusp form, we may replace $w_2$ by $w_2^{\ast}$.
By this procedure, we find the parabolic condition satisfied by $f$.

Next let $\Gamma ={\rm SL}(2, \mathcal O_F)$ and let $\epsilon$ be
the fundamental unit of $F$. We assume that $l_1 \equiv l_2 \mod 2$
and replace $\Gamma$ by ${\rm PSL}(2, \mathcal O_F)$. Put
$$
\sigma =\begin{pmatrix} 0 & 1 \\ -1 & 0 \end{pmatrix}, \qquad
\mu =\begin{pmatrix} \epsilon & 0 \\ 0 & \epsilon ^{-1} \end{pmatrix}.
$$
We choose $w_1=i\epsilon ^{-1}$, $w_2=i\infty$. Then we have
\begin{equation*}
f(\sigma , \mu )=f(\sigma , \sigma )
=-\int _{i\epsilon ^{-1}}^{i\epsilon} \int _0^{i\infty} \mathfrak d(\Omega ).
\end{equation*}
For $0 \leq s \leq l_1$, $0 \leq t \leq l_2$, we put
\begin{equation*}
P_{s,t}=\int _{i\epsilon ^{-1}}^{i\epsilon}\int _0^{i\infty}
\Omega (z)z_1^sz_2^t dz_1dz_2.
\end{equation*}
The $(l_1+1)(l_2+1)$ components of $f(\sigma , \mu )$ are given by
$-P_{s, t}$. We have

\begin{equation}
P_{m, m-(k_1-k_2)/2}=(-1)^{m+1}i^{-(k_1-k_2)/2}
(2\pi )^{(k_1-k_2)/2}R(m+1, \Omega )
\tag{6}
\end{equation}
where $R(s,\Omega )=(2\pi )^{-2s}\Gamma (s)\Gamma (s-(k_1-k_2)/2)L(s, \Omega )$
with the $L$-function $L(s, \Omega )$ of $\Omega$.
The formula (6) gives a generalization of (2); (5) and (6) were known
to the author eight years ago.

The $L$-value $L(m, \Omega )$, $m \in \mathbf Z$ is a critical value if and only if
\begin{equation*}
\frac {l_1-l_2}{2}+1 \leq m \leq \frac {l_1+l_2}{2}+1.
\end{equation*}
Since all of them appear as components of $f(\sigma , \mu )$,
we expect that we can deduce information on critical values
once we know the second cohomology group $H^2(\Gamma , V)$ well.
Before to materialize this hope, we need to answer the following conceptual question:
Can we annihilate the effect of adding a coboundary to $f$?
We can give an affirmative answer to this question by using the parabolic condition.
Put
$$
P=\left\{ \begin{pmatrix} u & v \\ 0 & u^{-1} \end{pmatrix} \bigg|
\ u \in \mathcal O_F^{\times}, \ v \in \mathcal O_F \right\}/ \{ \pm 1_2 \} \subset \Gamma .
$$
Then we have
\begin{equation}
f(p\gamma _1, \gamma _2)=p f(\gamma _1, \gamma _2)
\qquad \text {for every} \quad p \in P, \ \gamma _1, \gamma _2 \in \Gamma .
\tag{7}
\end{equation}
This is the parabolic condition on $f$ when $\Gamma ={\rm PSL}(2, \mathcal O_F)$.
A $2$-cocycle which satisfies (7) will be called a {\it parabolic $2$-cocycle}.
In section 3, we will prove:

\proclaim Theorem. Let $i=1$ or $2$. Then
$$
\dim H^i(P,V)=
\begin{cases}
0 \qquad {\rm if} \quad l_1\neq l_2 \ {\rm or} \ N(\epsilon )^{l_1}=-1,\\
1 \qquad {\rm if} \quad l_1=l_2 \ {\rm and} \ N(\epsilon )^{l_1}=1.
\end{cases}
$$
\par

\noindent Now suppose that we add a coboundary to $f$ keeping
the parabolic condition (7). In section 4, using this theorem for the case $i=1$,
we will show: If $l_1 \neq l_2$, the components of
$f(\sigma , \mu )$ related to the critical values do not change.
If $l_1=l_2$, the same assertion holds except for the critical values on the edges:
$L(1, \Omega )$ and $L(l_1+1, \Omega )$. Therefore
we can deduce information on critical values $L(m, \Omega )$
once we know a parabolic $2$-cocycle corresponding to $\Omega$.

The final step is to find constraints on $f(\sigma , \mu )$ which generalizes (3).
This is technically the most difficult step.
Let $\mathcal O_F=\mathbf Z+\mathbf Z\omega$ and put
$$
\tau =\begin{pmatrix} 1 & 1 \\ 0 & 1 \end{pmatrix}, \qquad
\eta =\begin{pmatrix} 1 & \omega \\ 0 & 1 \end{pmatrix}.
$$
It is known ([V]) that $\Gamma$ is generated by $\sigma$, $\mu$, $\tau$ and $\eta$.
Let $\mathcal F$ be the free group on four letters $\widetilde\sigma$,
$\widetilde\mu$, $\widetilde\tau$, $\widetilde\eta$. Define a surjective homomorphism
$\pi : \mathcal F \longrightarrow \Gamma$ by $\pi (\widetilde\sigma )=\sigma$,
$\pi (\widetilde\mu )=\mu$, $\pi (\widetilde\tau )=\tau$,
$\pi (\widetilde\eta )=\eta$ and let $R$ be the kernel of $\pi$.
Then we have $\Gamma =\mathcal F/R$ and (cf. \S1.4)
\begin{equation}
H^2(\Gamma , V)\cong H^1(R, V)^{\Gamma}/{\rm Im}(H^1(\mathcal F, V)).
\tag{8}
\end{equation}
Here we have
\begin{equation*}
H^1(R, V)^{\Gamma}=\bigg\{ \varphi \in {\rm Hom}(R, V) \mid
\varphi (grg^{-1})=g\varphi (r), \quad g \in \mathcal F, r \in R \bigg\} .
\end{equation*}
We write
\begin{equation*}
\epsilon ^2=A+B\omega, \qquad
\epsilon ^2\omega =C+D\omega .
\end{equation*}
We have relations:
(i) $\sigma ^2=1$.
(ii) $(\sigma\tau ) ^3=1$.
(iii) $(\sigma\mu ) ^2=1$.
(iv) $\tau\eta =\eta\tau$.
(v) $\mu\tau\mu ^{-1}=\tau ^A\eta ^B$.
(vi) $\mu\eta\mu ^{-1}=\tau ^C\eta ^D$. For
$t \in \mathcal O_F^{\times}$, we have
\begin{equation}
\sigma\begin{pmatrix} 1 & t \\ 0 & 1 \end{pmatrix}\sigma
=\begin{pmatrix} 1 & -t^{-1} \\ 0 & 1 \end{pmatrix}\sigma
\begin{pmatrix} -t & 1\\ 0 & -t^{-1} \end{pmatrix}.
\tag{vii}
\end{equation}
The relation (ii) follows from (vii) by taking $t=1$.
We call the relation group $R$ {\it minimal} if it is generated by
the elements corresponding to (i) $\sim$ (vii) and their conjugates.
We see that $\mu$, $\tau$ and $\eta$ genarate $P$ and (iv) $\sim$ (vi) are
their fundamental relations.

Now let $\varphi \in H^1(R, V)^{\Gamma}$ be a corresponding element to $f$.
Adding an element of ${\rm Im}(H^1(\mathcal F, V))$, we may assume that
$\varphi (\widetilde\sigma ^2)=0$. Then we find (cf. (5.3))
$$
f(\sigma , \mu )=-\varphi ((\widetilde\sigma\widetilde\mu )^2).
$$
Our problem is reduced to find constraints on
$\varphi ((\widetilde\sigma\widetilde\mu )^2)$. We have an obvious constraint
$\sigma\mu \varphi ((\widetilde\sigma\widetilde\mu )^2)
=\varphi ((\widetilde\sigma\widetilde\mu )^2)$ but of course it is not enough.

To proceed further, we assume that $l_1$ and $l_2$ are even and change $\rho$ to
$\rho ^{\prime}=\rho _{l_1}^{\prime}\otimes \rho _{l_2}^{\prime}$
where $\rho _l^{\prime}(g)=\rho _l(g)\det (g)^{-l/2}$ and regard $V$ as
a ${\rm PGL}(2, \mathcal O_F)$-module. The $\Gamma$-module structure
does not change. We put
$$
\nu =\begin{pmatrix} \epsilon & 0 \\ 0 & 1 \end{pmatrix}, \qquad
\delta =\begin{pmatrix} -1 & 0 \\ 0 & 1 \end{pmatrix}.
$$
These two elements act on $\Gamma$ as outer automorphisms and induce
automorphisms of $H^2(\Gamma , V)$ of order $2$.
Hence $H^2(\Gamma , V)$ decomposes into four pieces under their actions.
Let $\Gamma ^{\ast}$ be the subgroup of ${\rm PGL}(2, \mathcal O_F)$
generated by $\Gamma$ and $\nu$. The transfer map gives an isomorphism
of the plus part of $H^2(\Gamma , V)$ under the action of $\nu$ onto
$H^2(\Gamma ^{\ast}, V)$. For simplicity suppose that we can take $\omega =\epsilon$.
Then $\sigma$, $\nu$ and $\tau$ generate $\Gamma ^{\ast}$.
Let $\mathcal F^{\ast}$ be the free group on three letters $\widetilde\sigma$,
$\widetilde\nu$, $\widetilde\tau$. Define a surjective homorphism
$\pi ^{\ast}: \mathcal F^{\ast} \longrightarrow \Gamma ^{\ast}$
by $\pi ^{\ast}(\widetilde\sigma )=\sigma$,
$\pi ^{\ast}(\widetilde\nu )=\nu$, $\pi ^{\ast}(\widetilde\tau )=\tau$
and let $R^{\ast}$ be the kernel of $\pi ^{\ast}$.
Then we have $\Gamma ^{\ast}=\mathcal F^{\ast}/R^{\ast}$ and
\begin{equation}
H^2(\Gamma ^{\ast}, V)\cong H^1(R^{\ast}, V)^{\Gamma^{\ast}}
/{\rm Im}(H^1(\mathcal F^{\ast}, V)).
\tag{8$^{\ast}$}
\end{equation}
Let $f^{\ast}$ be the transfer of $f$ to $\Gamma ^{\ast}$
and let $f^+$ be the restriction of $f^{\ast}$ to $\Gamma$.
Then $f^+$ is the projection of $f$ to the plus part.
(We perform this procedure on the cocycle level.)
We have
$$
f^{\ast}(\sigma , \mu )=f^+(\sigma , \mu )=(1+\nu )f(\sigma , \mu ).
$$
In $\Gamma ^{\ast}$, $\sigma$, $\nu$ and $\tau$ satisfy the relations (i), (ii) and
(iii$^{\ast}$): $(\sigma\nu )^2=1$, (iv$^{\ast}$):
$\tau\nu\tau\nu ^{-1}=\nu\tau\nu ^{-1}\tau$,
(v$^{\ast}$): $\nu ^2\tau\nu ^{-2}=\tau ^A(\nu\tau\nu ^{-1})^B$.
Let $P^{\ast}$ be the subgroup of $\Gamma ^{\ast}$ generated by
$P$ and $\nu$. We see that $P^{\ast}$ is generated by $\nu$ and $\tau$
and (iv$^{\ast}$) and (v$^{\ast}$) are the fundamental relations between
generators $\nu$ and $\tau$. Let
$\varphi ^{\ast} \in H^1(R^{\ast}, V)^{\Gamma^{\ast}}$
be a corresponding element to $f^{\ast}$. By the parabolic condition on $f$,
we may assume that $\varphi ^{\ast}$ vanishes on the elements of $R^{\ast}$
corresponding to (iv$^{\ast}$) and (v$^{\ast}$).
Adding an element of
${\rm Im}(H^1(\mathcal F^{\ast}, V))$, we may also assume that
$\varphi ^{\ast}(\widetilde\sigma ^2)=0$. Then we have (cf. (6.6))
$$
f^{\ast}(\sigma , \mu )
=-(1+\nu ^{-1})\varphi ^{\ast}((\widetilde\sigma\widetilde\nu )^2)
$$
and two quantities
$$
A=\varphi ^{\ast}((\widetilde\sigma\widetilde\nu )^2), \qquad
B=\varphi ^{\ast}((\widetilde\sigma\widetilde\tau )^3)
$$
remain to be determined. The Hecke operators act on
$H^2(\Gamma ^{\ast}, V)$. We can analyze its action on the right-hand side
of (8$^{\ast}$) and will give a simple formula for it.
The quantity $A$ is related to the critical values of $L(s, \Omega )$.
We may assume that the class of $f^{\ast}$ is in the plus space
of $H^2(\Gamma ^{\ast}, V)$ under the action of $\delta$.
Then $A$ must satisfy the constraints
\begin{equation}
(\sigma\nu -1)A=0, \qquad (\delta -1)A=0.
\tag{9}
\end{equation}

We will execute the determination of $A$ for $F={\mathbf Q}(\sqrt {5})$
and $F={\mathbf Q}(\sqrt {13})$. First assume $F={\mathbf Q}(\sqrt {5})$.
In this case, we can show that $R$ is minimal and that $R^{\ast}$ is generated
by the elements corresponding to the relations (i), (ii), (iii$^{\ast}$),
(iv$^{\ast}$), (v$^{\ast}$) and their conjugates.
Calculating the action of the Hecke operator $T(2)$
on the right hand side of (8$^{\ast}$),
we find a certain element $x \in \mathcal F^{\ast}$ such that $\pi ^{\ast}(x)^3=1$.
We can give an explicit formula expressing $\varphi ^{\ast}(x^3)$ in terms of
$A$ and $B$. In every case examined, we find by numerical computations
that we may assume that $B=0$ by adding an element of
${\rm Im}(H^1(\mathcal F^{\ast}, V))$. Therefore 
\begin{equation}
(x-1)\varphi ^{\ast}(x^3)=0
\tag{10}
\end{equation}
gives a new constraint on $A$. Let $Z_A^+$ be the subspace of $V$
consisting of all $A$ which satisfy (9) and (10) and let $B_A^+$
be the subspace of $Z_A^+$ which represents the contribution from
${\rm Im}(H^1(\mathcal F^{\ast}, V))$. Again in every case examined,
we find by numerical computations that $\dim S_{l_1+2, l_2+2}(\Gamma )
=\dim Z_A^+/B_A^+$. If this is one dimensional, we can deduce information
on $L$-values by calculating $Z_A^+$. In general,
calculating the action of $T(2)$ on $Z_A^+/B_A^+$ and taking eigenvectors,
we can obtain many examples on $L$-values. 
Actually by considering $f^+$, we are losing half of the information
on critical values (cf. \S5.6).
To treat all critical values, we need to consider $f^-$,
the projection of $f$ to the minus part of $H^2(\Gamma , V)$ under
the action of $\nu$. To handle $f^-$ is a somewhat more complicated task and we leave
the explanation of it to the text.
Next let $F={\mathbf Q}(\sqrt {13})$. The procedure is almost the same.
Let $\mathfrak p$ be the prime ideal
generated by $4-\sqrt {13}$. Calculating the action of
the Hecke operator $T(\mathfrak p)$,
we obtain a certain element $x \in \mathcal F^{\ast}$ such that
$\pi ^{\ast}(x)$ is of order $3$. Then the constraint $(x-1)\varphi ^{\ast}(x^3)=0$
obtained from $x$ is sufficient. Here remarkably we can perform rigorous calculations
without proving that $R$ is minimal.
(This is actually true also for the case $F=\mathbf Q(\sqrt {5})$.)
We have used Pari [PARI] for the numerical calculations in sections 6 and 7.

To calculate the ratios of critical values of L-functions, there is another method
initiated by Shimura [Sh3] which employs the Rankin-Selberg convolution
and differential operators. 
A comparison of this method and the cohomological method will be discussed
in section 8.

\medskip

Now let us explain the organization of this paper briefly.
In section 1, we will review several facts on cohomology of a group
which will be repeatedly used in later sections.
In section 2, we will review Hilbert modular forms. We will prove
(5) and (6). In section 3, we will study cohomology groups of $P$
and will prove the theorem stated above.
In section 4, we will examine the parabolic condition on a cocycle
applying results in section 3. We will prove the nonvanishing of the cohomology class
of $f$ under mild conditions. In section 5, we will study the decomposition of
$H^2(\Gamma , V)$ under the action of outer automorphisms of $\Gamma$.
It decomposes into four pieces under this action. In section 6, we will study
the case $F=\mathbf Q(\sqrt {5})$ in detail and will give many examples.
In section 7, we will study the case $F=\mathbf Q(\sqrt {13})$.
We devote section 8 to the comparison of two methods mentioned above.
In Appendix, we will prove that the relation group $R$ is minimal when
$F=\mathbf Q(\sqrt {5})$.

\bigskip

{\bf Notation.} For an associative ring $A$ with identity element,
$A^{\times}$ denotes the group of all invertible elements of $A$.
Let $R$ be a commutative ring with identity element.
We denote by $M(n, R)$ the ring of all $n \times n$ matrices
with entries in $R$. We define ${\rm GL}(n, R)=M(n,R)^{\times}$,
${\rm SL}(n, R)=\{ g \in {\rm GL}(n, R) \mid \det g=1 \}$.
The quotient group of ${\rm GL}(n, R)$ (resp. ${\rm SL}(n, R)$) by its center
is denoted by ${\rm PGL}(n, R)$ (resp. ${\rm PSL}(n, R)$).
Let $G$ be a group. The subgroup of $G$ generated by $g_1, \ldots , g_n \in G$
is denoted by $\langle g_1, \ldots , g_n\rangle$. When $G$ acts on a module $M$,
$M^G$ denotes the submodule of $M$ consisting of all elements fixed by
$G$. For an algebraic number field $F$, $\mathcal O_F$ denotes
the ring of integers of $F$. For $a \in \mathcal O_F$, the ideal
$a\mathcal O_F$ generated by $a$ is denoted by $(a)$.
We denote by $E_F$ the group of units of $F$,
i.e., $E_F=\mathcal O_F^{\times}$. When $F$ is totally real and $\alpha \in F$,
$\alpha \gg 0$ means that $\alpha$ is totally positive. We denote by
$\mathfrak H$ the complex upper half plane. The set of all positive real numbers
is denoted by $\mathbf R_+$.

\bigskip

\centerline{\S1. Preparations on cohomology groups}

\medskip

In this section, we will review group cohomology. Most of the results,
except for the results presented in subsection 1.5,
can be found in standard text books such as Cartan--Eilenberg [CE],
Serre [Se1], Suzuki[Su].

\medskip

{\bf 1.1.} Let $G$ be a group, $M$ be a left $G$-module.
We set $C^0(G,M)=M$, and for $0<n \in \mathbf Z$, let $C^n(G, M)$
be the abelian group consisting of all mappings of $G^n$ into $M$.
We define the coboundary operator
$d_n: C^n(G, M) \longrightarrow C^{n+1}(G, M)$ by the usual formula
\begin{equation*}
\begin{aligned}
(d_nf)(g_1, \ldots , & g_{n+1})= g_1f(g_2, \ldots , g_{n+1})
+(-1)^{n+1}f(g_1, \ldots , g_n) \\
+ & \sum _{i=1}^n (-1)^i f(g_1, \ldots , g_ig_{i+1}, \ldots , g_{n+1}).
\end{aligned}
\tag{1.1}
\end{equation*}
We set
$$
Z^n(G, M)={\rm Ker}(d_n), \qquad B^n(G, M)={\rm Im}(d_{n-1}).
$$
Here we understand $B^0(G, M)=\{ 0 \}$.
An element of $C^n(G,M)$ (resp. $Z^n(G, M)$, resp. $B^n(G, M)$)
 is called an $n$-cochain (resp. $n$-cocycle, resp. $n$-coboundary).
The cohomology group $H^n(G,M)$ is that of the complex
$\{ C^n(G,M), d_n \}$, i.e., $H^n(G,M)=Z^n(G,M)/B^n(G,M)$.

Let $G^{\prime}$ be a group and $M^{\prime}$ be a left $G^{\prime}$-module.
Let $\varphi : G \longrightarrow G^{\prime}$ be a group homomorphism
and $\psi : M^{\prime} \longrightarrow M$ be a homomorphism of abelian groups.
We assume that $\varphi$ and $\psi$ are compatible, that is
$$
\psi (\varphi (g)m^{\prime})=g(\psi (m^{\prime})),
\qquad m^{\prime} \in M^{\prime}, \quad g \in G.
$$
For $f \in C^n(G^{\prime}, M^{\prime})$, define
$\omega _nf \in C^n(G, M)$ by the formula
\begin{equation}
(\omega _nf)(g_1, g_2, \ldots , g_n)
=\psi (f(\varphi (g_1), \varphi (g_2), \ldots, \varphi (g_n)).
\tag{1.2}
\end{equation}
Then we can check easily that the following diagram is commutative.
\begin{equation*}
\begin{CD}
C^n(G^{\prime} , M^{\prime}) @>\omega _n>> C^n(G, M) \\
@VV d_n V @VV d_n V \\
C^{n+1}(G^{\prime} , M^{\prime}) @>\omega _{n+1}>> C^{n+1}(G, M).
\end{CD}
\end{equation*}
Therefore $\omega _n$ sends $Z^n(G^{\prime}, M^{\prime})$
(resp. $B^n(G^{\prime}, M^{\prime})$) into
$Z^n(G, M)$ (resp. $B^n(G, M)$) and induces a homomorphism
$H^n(G^{\prime}, M^{\prime}) \longrightarrow H^n(G, M)$.

Now let $N$ be a subgroup of $G$. Let $g \in G$.
We define
$$
\varphi (n)=gng^{-1}, \quad n \in g^{-1}Ng, \qquad
\psi (m)=g^{-1}m, \quad m \in M.
$$
Then $\varphi$ is an isomorphism of $g^{-1}Ng$ onto $N$;
$\varphi$ and $\psi$ are compatible.
Hence we obtain an isomorphism of $H^p(N, M)$ onto $H^p(g^{-1}Ng,M)$,
which is induced by sending $f \in Z^p(N, M)$ to
$f^{\prime} \in Z^p(g^{-1}Ng, M)$:
\begin{equation}
f^{\prime}(n_1, n_2, \ldots , n_p)
=g^{-1}f(gn_1g^{-1}, gn_2g^{-1}, \ldots , gn_pg^{-1}).
\tag{1.3}
\end{equation}

{\bf 1.2.} Let $H$ be a subgroup of $G$ of finite index.
We are going to consider an explicit form of the transfer map
$H^n(H, M) \longrightarrow H^n(G, M)$ (cf. Eckmann [E]). To this end,
it is convenient to go back to a more conceptual definition of group cohomology:
$$
H^n(G, M)={\rm Ext}_G^n(\mathbf Z, M).
$$
The right-hand side can be computed as follows. We take a resolution of
$\mathbf Z$ by projective $G$-modules $P_n$.
\footnote{The $d_n$ in the diagram below should not be confused with
$d_n$ in (1.1).}
\begin{equation*}
\begin{CD}
\cdots  @>>> P_2 @>d_1>>P_1 @>d_0>>P_0 @>\epsilon >> \mathbf Z @>>>0.
\end{CD}
\end{equation*}
Then we obtain a complex
\begin{equation*}
\begin{CD}
0  @>>> {\rm Hom}_G(P_0, M)  @>d_0^{\ast}>>
{\rm Hom}_G(P_1, M)  @>d_1^{\ast}>>
{\rm Hom}_G(P_2, M)  @>d_2^{\ast}>> \cdots
\end{CD}
\end{equation*}
which gives the cohomology group
$H^n(G, M)={\rm Ker} (d_n^{\ast})/{\rm Im}(d_{n-1}^{\ast})$.
(If $n=0$, we understand ${\rm Im}(d_{n-1}^{\ast})=0$.)
As is well known, an explicit form of a resolution can be given as follows.
Let $P_n$ be the free abelian group on the base $G^{n+1}$. We give
$P_n$ a $G$-module structure by
$$
g(g_0, g_1, \ldots , g_n)=(gg_0, gg_1, \ldots , gg_n), \qquad g \in G
$$
and define $d_{n-1}: P_n \longrightarrow P_{n-1}$ by
$$
d_{n-1}(g_0, g_1, \ldots , g_n)=\sum _{i=0}^n
(-1)^i(g_0, \ldots , g_{i-1}, g_{i+1}, \ldots , g_n), \qquad n \geq 1.
$$
We set $\epsilon (g_0)=1$. We have
$$
{\rm Hom}_G(P_n, M)
=\left\{ \varphi : G^{n+1} \longrightarrow M
\mid \varphi (gg_0, \ldots , gg_n)=g\varphi (g_0, \ldots , g_n) \right\}.
$$
An element of ${\rm Ker}(d_n^{\ast})$ is called a {\it homogeneous $n$-cocycle}.
To $\varphi \in {\rm Hom}_G(P_n, M)$, we let $f \in C^n(G, M)$
correspond by the formula
\begin{equation}
f(g_1, g_2, \ldots , g_n)=\varphi (1, g_1, g_1g_2, \ldots , g_1g_2\cdots g_n).
\tag{1.4}
\end{equation}
Then $\varphi \mapsto f$ gives an isomophism of abelian groups.
The coboundary operator $d_n^{\ast}$ induces $d_n$ on $C^n(G,M)$
which is given by (1.1).

Now let $H$ be a subgroup of $G$ of finite index $r$ and let
$$
G=\sqcup _{i=1}^r x_i H
$$
be a coset decomposition. For $G$-modules $A$ and $B$, we define
a homomorphism
$t: {\rm Hom}_H(B, A) \longrightarrow {\rm Hom}_G(B, A)$ by
$$
(t\varphi )(b)=\sum _{i=1}^r x_i\varphi (x_i^{-1}b).
$$
We denote $P_n$ for $G$ (resp. $H$) by $P_n^G$ (resp. $P_n^H$).
Since $P_n^G$ is a free $H$-module, the complex
$\{ {\rm Hom}_H(P_n^G, M), d_n^{\ast} \}$ gives the cohomology group
$H^n(H, M)$. Define $P_n^G \longrightarrow P_n^H$ by
$(g_0, \ldots , g_n) \mapsto (h_0, \ldots , h_n)$ where $h_i$ is determined by
\begin{equation}
g_i=h_ix_{j(i)}^{-1}, \qquad h_i \in H, \quad 0 \leq i \leq n.
\tag{1.5}
\end{equation}
This is an $H$-homomorphism and commutes with $d_n$;
it induces an isomophism between the cohomology groups
of complexes $\{ {\rm Hom}_H(P_n^H, M), d_n^{\ast} \}$ and
$\{ {\rm Hom}_H(P_n^G, M), d_n^{\ast} \}$.
Now the following diagram is commutative.
\begin{equation*}
\begin{CD}
0  @>>> {\rm Hom}_G(P_0^G, M)  @>d_0^{\ast}>>
{\rm Hom}_G(P_1^G, M)  @>d_1^{\ast}>>
{\rm Hom}_G(P_2^G, M)  @>d_2^{\ast}>> \cdots \\
@. @AA t A @AA t A @AA t A \\
0  @>>> {\rm Hom}_H(P_0^G, M)  @>d_0^{\ast}>>
{\rm Hom}_H(P_1^G, M)  @>d_1^{\ast}>>
{\rm Hom}_H(P_2^G, M)  @>d_2^{\ast}>> \cdots
\end{CD}
\end{equation*}
Hence it gives the transfer map $T:H^n(H, M) \longrightarrow H^n(G, M)$.

Take $\varphi \in {\rm Hom}_H(P_n^H, M)$. The corresponding element
$\varphi ^{\ast} \in {\rm Hom}_H(P_n^G, M)$ to $\varphi$ is given by
$$
\varphi ^{\ast}(g_0, g_1, \ldots , g_n)
=\varphi (g_0x_{j(0)}, g_1x_{j(1)}, \ldots , g_nx_{j(n)}),
$$
where $x_{j(i)}$ is determined by (1.5). Then
$\tilde\varphi =t\varphi ^{\ast}\in {\rm Hom}_G(P_n^G, M)$ is given by
$$
\tilde\varphi (g_0, g_1, \ldots , g_n)
=\sum _{j=1}^r x_j\varphi ^{\ast}(x_j^{-1}g_0, x_j^{-1}g_1, \ldots , x_j^{-1}g_n).
$$
Therefore we obtain the following result.

\proclaim Proposition 1.1. Let $G$ be a group, $H$ be a subgroup of finite index
and $M$ be a left $G$-module. Let $G=\sqcup _{i=1}^r x_iH$
be a coset decomposition. Let $c \in H^n(H,M)$ and
let $\varphi \in {\rm Hom}_H(P_n^H, M)$ be a homogeneous $n$-cocycle
representing $c$. Then a homogeneous $n$-cocycle
$\tilde\varphi \in {\rm Hom}_G(P_n^G, M)$ which represents $T(c)$
is given by
$$
\tilde\varphi (g_0, g_1, \ldots g_n)
=\sum _{i=1}^r x_i\varphi (x_i^{-1}g_0x_{j_i(0)},
x_i^{-1}g_1x_{j_i(1)}, \ldots , x_i^{-1}g_nx_{j_i(n)}).
$$
Here $x_{j_i(k)}$ is chosen so that $x_i^{-1}g_kx_{j_i(k)} \in H$.
\par

Then using (1.4), we immediately deduce:

\proclaim Proposition 1.2. Let the notation be the same as in Proposition 1.1.
Let $f \in Z^n(H, M)$ be an $n$-cocycle representing $c \in H^n(H,M)$.
Then an $n$-cocycle $\tilde f \in Z^n(G, M)$ which represents
$T(c) \in H^n(G, M)$ is given by
$$
\tilde f(g_1, g_2, \ldots , g_n)
=\sum _{i=1}^r x_i f(x_i^{-1}g_1x_{p_i(1)}, x_{p_i(1)}^{-1}g_2x_{p_i(2)},
\ldots , x_{p_i(n-1)}^{-1}g_nx_{p_i(n)}).
$$
Here $x_{p_i(l)}$ is chosen so that
$$
x_i^{-1}g_1x_{p_i(1)} \in H, \qquad
x_{p_i(l-1)}^{-1}g_lx_{p_i(l)} \in H, \quad 2 \leq l \leq n.
$$
\par

Let ${\rm Res}: H^n(G, M) \longrightarrow H^n(H,M)$ be the restriction
homomorphism. Then we have the well-known result:
\begin{equation}
T\circ {\rm Res}(c)=[G:H]c, \qquad c \in H^n(G,M).
\tag{1.6}
\end{equation}

\medskip

{\bf 1.3.} We are going to consider the action of Hecke operators
on cohomology groups. Let $\widetilde G$ be a group and $G$ be a subgroup.
Let $M$ be a $\widetilde G$-module. We assume that
$G$ and $tGt^{-1}$ are commensurable for every $t \in \tilde G$.
For $t \in \widetilde G$, we put
$$
G_t=G \cap t^{-1}Gt.
$$
Let
$$
{\rm conj}: H^n(G, M) \longrightarrow H^n(t^{-1}Gt, M)
$$
be the isomorphism induced by (1.3). Let {\rm Res} be the restriction map
from $H^n(t^{-1}Gt, M)$ to $H^n(G_t, M)$ and let
$T: H^n(G_t, M) \longrightarrow H^n(G,M)$ be the transfer map.
Then we define
\begin{equation}
[GtG]=T\circ{\rm Res}\circ{\rm conj}.
\tag{1.7}
\end{equation}
(It is not difficult to check that the right-hand side of (1.7) depends only on
the double coset $GtG$ and that (1.7) defines a homomorphism
of the Hecke ring $H(G, \widetilde G)$ into ${\rm End}(H^n(G, M))$.)
Let us write an explicit form of this operator when $n=2$, which
will be necessary for our later computation. Let
$$
G=\sqcup _{i=1}^d G_t\alpha _i
$$
be a coset decomposition. Then we have
$$
GtG=\sqcup _{i=1}^d Gt\alpha _i.
$$
Put $\beta _i=t\alpha _i$. Let $c \in H^2(G, M)$ and let $f \in Z^2(G, M)$
be a $2$-cocycle which represents $c$. By (1.3), ${\rm conj}(c)$
is represented by the $2$-cocycle
$f^{\prime}(g_1, g_2)=t^{-1}f(tg_1t^{-1}, tg_2t^{-1})$. By Proposition 1.2,
$[GtG](c)$ is represented by the $2$-cocycle
$$
f^{\prime\prime}(g_1, g_2)=\sum _{i=1}^d \alpha _i^{-1}
f^{\prime}(\alpha _ig_1\alpha _{j(i)}^{-1}, \alpha _{j(i)}g_2\alpha _{k(j(i))}^{-1}),
$$
since $G=\sqcup _{i=1}^d \alpha _i^{-1}G_t$.
Here, for $1 \leq i \leq d$, we choose $j(i)$ and $k(i)$ so that
$$
\alpha _ig_1\alpha _{j(i)}^{-1} \in G_t, \qquad
\alpha _ig_2\alpha _{k(i)}^{-1} \in G_t.
$$
Writing the result in terms of $\beta _i$, we obtain the following proposition.

\proclaim Proposition 1.3. Let $c \in H^2(G, M)$ and let $f \in Z^2(G, M)$ be
a $2$-cocycle representiong $c$. Let
$GtG=\sqcup _{i=1}^d G\beta _i$ be a coset decomposition.
Then a $2$-cocycle $h \in Z^2(G, M)$ representing $[GtG](c)$ is given by
$$
h(g_1, g_2)=\sum _{i=1}^d \beta _i^{-1}f(\beta _ig_1\beta _{j(i)}^{-1},
\beta _{j(i)}g_2\beta _{k(j(i))}^{-1}).
$$
Here, for $1 \leq i \leq d$, we choose $j(i)$ and $k(i)$ so that
$$
\beta _ig_1\beta _{j(i)}^{-1} \in G, \qquad
\beta _ig_2\beta _{k(i)}^{-1} \in G.
$$
\par

{\bf 1.4.} Let $G$ be a group and $M$ be a left $G$-module.
Let $N$ be a normal subgroup of $G$. Then we have
the Hochschild-Serre spectral sequence
\begin{equation}
E_2^{p,q}=H^p(G/N, H^q(N,M)) \Longrightarrow H^n(G,M).
\tag{1.8}
\end{equation}
In low dimensions, this gives an exact sequence
\begin{equation}
\begin{aligned}
\begin{CD}
0 @> >> H^1(G/N, M^N) @> >> H^1(G,M)
@> >> H^1(N,M)^{G/N}
\end{CD}
\\
\begin{CD}
@> >> H^2(G/N, M^N) @> >> H^2(G,M).
\end{CD}
\end{aligned}
\tag{1.9}
\end{equation}

Now we are going to describe a method to calculate $H^2(G, M)$,
which is originally due to MacLane (cf. [K], \S50).
Taking a free group $\mathcal F$, we write $G=\mathcal F/R$.
Let $\pi : \mathcal F \longrightarrow G$ be the canonical homomorphism
such that ${\rm Ker}(\pi )=R$. We regard $M$ as an $\mathcal F$-module by
$gm=\pi (g)m$, $g \in \mathcal F$, $m \in M$. Since
\begin{equation}
H^i(\mathcal F,M)=0, \qquad i \geq 2,
\tag{1.10}
\end{equation}
(1.9) yields an exact sequence
\begin{equation*}
\begin{CD}
0 @> >> & H^1(G, M) @> >> H^1(\mathcal F,M)
@> >> H^1(R,M)^G @> >> H^2(G ,M) @> >> 0.
\end{CD}
\end{equation*}
Therefore we have
\begin{equation}
H^2(G ,M) \cong H^1(R,M)^G/{\rm Im}(H^1(\mathcal F,M)).
\tag{1.11}
\end{equation}
Since $R$ acts on $M$ trivially, we have $B^1(R,M)=0$ and
$H^1(R,M)={\rm Hom}(R,M)$. Therefore we have
\begin{equation*}
H^1(R,M)^G=\{ \varphi \in {\rm Hom}(R, M) \mid
\varphi (grg^{-1})=g\varphi (r), \quad g \in \mathcal F, r \in R \} .
\end{equation*}

The isomorphism (1.11) is explicitly given as follows.
For $g \in \mathcal F$, we put $\pi (g)=\bar g$.
Take a $2$-cocycle $f \in Z^2(G , M)$. The mapping
$(g_1, g_2) \longrightarrow f(\bar g_1, \bar g_2)$ is an $M$-valued $2$-cocycle of
$\mathcal F$. By (1.10), there exists a $1$-cochain $a \in C^1(\mathcal F, M)$
such that
\begin{equation}
f(\bar g_1, \bar g_2)=g_1a(g_2)+a(g_1)-a(g_1g_2), \qquad
g_1, g_2 \in \mathcal F.
\tag{1.12}
\end{equation}
Let $\varphi =a\vert R$, the restriction of $a$ to $R$.
We may assume that $f$ is normalized, i.e.,
$$
f(1, g)=f(g, 1)=0 \qquad \text{for all} \ g \in G.
$$
If $r_1$, $r_2 \in R$, then, by (1.12), we have
$$
a(r_2)+a(r_1)-a(r_1r_2)=0.
$$
Therefore we get $\varphi \in Z^1(R,M)={\rm Hom}(R, M)$. By (1.12), we have
\begin{equation}
a(gr)=ga(r)+a(g), \qquad g \in \mathcal F, \ r \in R.
\tag{1.13}
\end{equation}
Again by (1.12), we have
\begin{equation*}
\begin{aligned}
a(grg^{-1}) & =gra(g^{-1})+a(gr)-f(\bar g, \bar g^{-1}) \\
& =ga(g^{-1})+ga(r)+a(g)-f(\bar g, \bar g^{-1})
\end{aligned}
\end{equation*}
for $g \in \mathcal F$, $r \in R$.
Using (1.12) with $g_1=g$, $g_2=g^{-1}$ and noting $a(1)=0$, we obtain
\begin{equation}
\varphi (grg^{-1})=g\varphi (r), \qquad g \in \mathcal F, r \in R.
\tag{1.14}
\end{equation}
This formula shows that $\varphi$ belongs to $H^1(R,M)^G$.
Suppose that $a^{\prime}$ is another $1$-cochain satisfying (1.12).
Put $\varphi ^{\prime}=a^{\prime}\vert R$, $a^{\prime}=a+b$.
Then $b \in Z^1(\mathcal F, M)$.
Hence the classes of $\varphi$ and $\varphi ^{\prime}$ in
$H^1(R, M)^G/{\rm Im}(H^1(\mathcal F, M))$ are the same.
Suppose that we add the coboundary of a $1$-cochain $c$ to $f$.
Then (1.12) holds when we replace $a(g)$ by $a(g)+c(\bar g)$.
Then $a\vert R$ does not change. Thus we have defined a homomorphism
$$
\omega : H^2(G ,M) \longrightarrow H^1(R,M)^G/{\rm Im}(H^1(\mathcal F,M)).
$$

Next suppose that $\varphi \in H^1(R,M)^G$. Take a coset decomposition
$\mathcal F=\sqcup _if_iR$. We assume that if $f_iR=R$, then $f_i=1$. 
We extend $\varphi$ to a mapping from $\mathcal F$ to $M$ as follows.
Choose $a(f_i) \in M$ in arbitrary way. Then put
\begin{equation}
a(f_ir)=f_i\varphi (r)+a(f_i), \qquad r \in R.
\tag{1.15}
\end{equation}
For $g_1=f_ir_1$, $g_2=f_jr_2$, $r_1$, $r_2 \in R$,
a direct calculation shows that
$$
g_1a(g_2)+a(g_1)-a(g_1g_2)=f_ia(f_j)+a(f_i)-a(f_k)-\varphi (f_if_jf_k^{-1}).
$$
Here $f_if_j=f_kr_3$, $r_3 \in R$. Note that $f_k$ does not depend on
$r_1$ and $r_2$.
Therefore we can define a $2$-cochain $f \in C^2(G , M)$ by (1.12)
Then it is immediate to see that $f \in Z^2(G , M)$ and that
$f$ is normalized (see Lemma 1.4 below). When we add an element of
${\rm Im}(H^1(\mathcal F, M))$ to $\varphi$, the cohomology class of $f$
does not change. Thus we have defined a homomorphism
$$
\eta : H^1(R,M)^G/{\rm Im}(H^1(\mathcal F,M))
\longrightarrow H^2(G ,M).
$$
We can check easily that $\omega$ and $\eta$ are inverse mappings to each other.
This finishes an explicit description of the isomorphism (1.11).

\medskip

{\bf 1.5.} Let $f \in Z^2(G, M)$ be a normalized cocycle. Take $a \in C^1(\mathcal F, M)$
which satisfies (1.12) and put $\varphi =a\vert R \in H^1(R, M)^G$.
For every $g \in G$, we choose $\widetilde g \in \mathcal F$ such that
$\pi (\widetilde g)=g$. The formula (1.12) can be written as
$$
f(g_1, g_2)=g_1a(\widetilde g_2)+a(\widetilde g_1)
-a(\widetilde g_1\widetilde g_2),
\qquad g_1, g_2 \in G.
$$
By (1.13), we have
$$
a(\widetilde{g_1g_2}(\widetilde{g_1g_2})^{-1}\widetilde g_1\widetilde g_2)
=g_1g_2\varphi((\widetilde{g_1g_2})^{-1}\widetilde g_1\widetilde g_2)
+a(\widetilde{g_1g_2}).
$$
Then, using (1.14), we have
$$
a(\widetilde g_1\widetilde g_2)=a(\widetilde{g_1g_2})
+\varphi (\widetilde g_1\widetilde g_2(\widetilde{g_1g_2})^{-1}).
$$
Therefore we obtain
\begin{equation}
f(g_1, g_2)=g_1a(\widetilde g_2)+a(\widetilde g_1)-a(\widetilde {g_1g_2})
-\varphi (\widetilde g_1\widetilde g_2(\widetilde {g_1g_2})^{-1}),
\qquad g_1, g_2 \in G.
\tag{1.16}
\end{equation}
This formula shows that, adding a coboundary to $f$, we may assume that
\begin{equation}
f(g_1, g_2)
=-\varphi (\widetilde g_1\widetilde g_2(\widetilde {g_1g_2})^{-1}).
\tag{1.17}
\end{equation}
Conversely we note the following Lemma.

\proclaim Lemma 1.4. Let $\varphi \in H^1(R. M)^G$.
For $g_1$, $g_2 \in G$, define $f(g_1, g_2)$ by (1.17).
Then $f \in Z^2(G, M)$. If $\widetilde 1=1$, $f$ is normalized.
\par
{\bf Proof.} The cocycle condition is
$$
g_1f(g_2, g_3)-f(g_1g_2, g_3)+f(g_1, g_2g_3)-f(g_1, g_2)=0.
$$
We have
\begin{equation*}
\begin{aligned}
& g_1\varphi (\widetilde g_2\widetilde g_3(\widetilde {g_2g_3})^{-1})
-\varphi (\widetilde {g_1g_2}\widetilde g_3(\widetilde {g_1g_2g_3})^{-1})
+\varphi (\widetilde g_1\widetilde {g_2g_3}(\widetilde {g_1g_2g_3})^{-1}) \\
-\, & \varphi (\widetilde g_1\widetilde g_2(\widetilde {g_1g_2})^{-1}) \\
=\, & \varphi (\widetilde g_1\widetilde g_2\widetilde g_3(\widetilde {g_2g_3})^{-1}
\widetilde g_1^{-1})
+\varphi (\widetilde g_1\widetilde {g_2g_3}(\widetilde {g_1g_2g_3})^{-1})
+\varphi (\widetilde {g_1g_2g_3}\widetilde g_3^{-1}(\widetilde {g_1g_2})^{-1}) \\
+\, & \varphi (\widetilde {g_1g_2}\widetilde g_2^{-1}\widetilde g_1^{-1}) \\
=\, & \varphi (\widetilde g_1\widetilde g_2\widetilde g_3(\widetilde {g_2g_3})^{-1}
\widetilde g_1^{-1})
+\varphi (\widetilde g_1\widetilde {g_2g_3}\widetilde g_3^{-1}
(\widetilde {g_1g_2})^{-1})
 +\varphi (\widetilde {g_1g_2}\widetilde g_2^{-1}\widetilde g_1^{-1}) \\
=\, & \varphi (\widetilde g_1\widetilde g_2(\widetilde {g_1g_2})^{-1})
+\varphi (\widetilde {g_1g_2}\widetilde g_2^{-1}\widetilde g_1^{-1}) =0.
\end{aligned}
\end{equation*}
Hence the cocycle condition holds. The latter assertion is obvious.
This completes the proof.

\medskip

We are going to write the action of Hecke operators
on the right-hand side of (1.11) explicitly. Let the notation be the same
as in subsections 1.3 and 1.4. Let $f \in Z^2(G, M)$ be a normalized
$2$-cocycle of the cohomology class $c$. Let $h$ be the $2$-cocycle
given by Proposition 1.3 which represents the class $[GtG](c)$.
Clearly $h$ is normalized. There exists a $1$-cochain $b \in C^1(\mathcal F, M)$
such that
\begin{equation*}
h(\bar g_1, \bar g_2)=g_1b(g_2)+b(g_1)-b(g_1g_2), \qquad
g_1, g_2 \in \mathcal F.
\end{equation*}

\proclaim Proposition 1.5. Let $\varphi \in H^1(R, M)^G$
and let a normalized $2$-cocycle $f \in Z^2(G,M)$ be given by (1.17).
Suppose $g_j\in G$ are given for $1 \leq j \leq m$.
For every $j$, we define a permutation on $d$ letters $p_j \in S_d$ by
$$
\beta _ig_j\beta _{p_j(i)}^{-1} \in G , \qquad 1 \leq i \leq d.
$$
We define $q_j \in S_d$ inductively by
$$
q_1=p_1, \qquad q_k=p_kq_{k-1}, \quad 2 \leq k \leq m.
$$
We assume that $b(\widetilde g_j)=0$ for $1 \leq j\leq m$.
Then we have
\begin{equation*}
\begin{aligned}
& b(\widetilde g_1\widetilde g_2\cdots\widetilde g_m) \\
= \, & \sum _{i=1}^d \beta _i^{-1}
\varphi (\widetilde{\beta _ig_1\beta _{q_1(i)}^{-1}}
\widetilde{\beta _{q_1(i)}g _2\beta _{q_2(i)}^{-1}}\cdots
\widetilde{\beta _{q_{m-1}(i)}g_m\beta _{q_m(i)}^{-1}}
(\widetilde{\beta _ig_1g_2\cdots g_m\beta _{q_m(i)}^{-1}})^{-1}).
\end{aligned}
\tag{1.18}
\end{equation*}
\par
{\bf Proof.} If $m=1$, the left-hand side of (1.18) is $0$
and the right-hand side is $0$ since $\varphi (1)=0$.
We assume that $m \geq 2$ and the formula is valid for $m-1$.
Then, by Proposition 1.3 and (1.17), we have
\begin{equation*}
\begin{aligned}
& b(\widetilde g_1\widetilde g_2\cdots\widetilde g_{m-1}\widetilde g_m) \\
= \, & g_1g_2\cdots g_{m-1}b(\widetilde g_m)
+b(\widetilde g_1\widetilde g_2\cdots\widetilde g_{m-1})
-h(g_1\cdots g_{m-1}, g_m) \\
= \, & \sum _{i=1}^d \beta _i^{-1}
\varphi (\widetilde{\beta _ig_1\beta _{q_1(i)}^{-1}}\cdots
\widetilde{\beta _{q_{m-2}(i)}g_{m-1}\beta _{q_{m-1}(i)}^{-1}}
(\widetilde{\beta _ig_1g_2\cdots g_{m-1}\beta _{q_{m-1}(i)}^{-1}})^{-1}) \\
+ \, & \sum _{i=1}^d \beta _i^{-1}
\varphi (\widetilde{\beta _ig_1g_2\cdots g_{m-1}
\beta _{q_{m-1}(i)}^{-1}}
\widetilde{\beta _{q_{m-1}(i)}g_m\beta _{q_m(i)}^{-1}}
(\widetilde{\beta _ig_1g_2\cdots g_m\beta _{q_m(i)}^{-1}})^{-1}) \\
= \, & \sum _{i=1}^d \beta _i^{-1}
\varphi (\widetilde{\beta _ig_1\beta _{q_1(i)}^{-1}}\cdots
\widetilde{\beta _{q_{m-1}(i)}g_m\beta _{q_m(i)}^{-1}}
(\widetilde{\beta _ig_1g_2\cdots g_m\beta _{q_m(i)}^{-1}})^{-1})
\end{aligned}
\end{equation*}
since $b(\widetilde g_m)=0$. This completes the proof.

\medskip

We have
\begin{equation*}
b(g_1g_2)=g_1b(g_2)+b(g_1)-h(\bar g_1, \bar g_2), \qquad
g_1, g_2 \in \mathcal F.
\end{equation*}
We may take $b(g)=0$ for a fixed set of generators of $\mathcal F$
and we can apply the above formula to determine the value of
$b(g)$ according to the length of $g \in \mathcal F$.
But Proposition 1.5 is useful beyond this case as will be seen after section 5.

\bigskip

\centerline{\S2. Hilbert modular forms}

\medskip

{\bf 2.1.} In this subsection, we follow the exposition given in Shimura [Sh4].
Let $F$ be a totally real algebraic number field of degree $n$.
Let $\mathfrak d_F$ denote the different of $F$ over $\mathbf Q$ and
let $\{ \sigma _1, \sigma _2, \ldots , \sigma _n \}$ be the set of
all isomorphisms of $F$ into $\mathbf R$. For $\xi \in F$, we put
$\xi ^{(\nu )}=\xi ^{\sigma _{\nu}}$. For
$z=(z_1, z_2, \ldots , z_n) \in \mathfrak H^n$, we put
$$
\mathbf e_F(\xi z)=\exp (2\pi i\sum _{\nu =1}^n \xi ^{(\nu )}z_{\nu}).
$$
Let $k=(k_1, k_2, \ldots , k_n) \in \mathbf Z^n$.
For $g=\begin{pmatrix} a & b \\ c & d \end{pmatrix} \in {\rm GL}(2, \mathbf R)_+$
and $z \in \mathfrak H$, we put $gz=(az+b)/(cz+d)$, $j(g, z)=cz+d$, where
${\rm GL}(2, \mathbf R)_+=\{ g \in {\rm GL}(2, \mathbf R) \mid \det g>0 \}$;
${\rm GL}(2, \mathbf R)_+^n$ acts on $\mathfrak H^n$.
For a function $\Omega$ on $\mathfrak H^n$,
$g=(g_1, \ldots , g_n) \in {\rm GL}(2, \mathbf R)_+^n$
and $z=(z_1, \ldots , z_n) \in \mathfrak H^n$, we define
a function $\Omega |_k \, g$ on $\mathfrak H^n$ by the formula
$$
(\Omega |_k \, g)(z)=\prod _{\nu =1}^n  \det (g_{\nu})^{k_{\nu}/2}
j(g_{\nu}, z_{\nu})^{-k_{\nu}}\Omega(gz).
$$
We embed ${\rm GL}(2, F)$ into ${\rm GL}(2, \mathbf R)^n$ by
$$
{\rm GL}(2, F) \ni \begin{pmatrix} a & b \\ c & d \end{pmatrix} \mapsto
\big(
\begin{pmatrix} a^{(1)} & b^{(1)} \\ c^{(1)} & d^{(1)} \end{pmatrix}, \ldots ,
\begin{pmatrix} a^{(n)} & b^{(n)} \\ c^{(n)} & d^{(n)} \end{pmatrix} \big )
\in {\rm GL}(2, \mathbf R)^n.
$$

Let $\Gamma$ be a congruence subgroup of $SL(2, \mathcal O_F)$.
A holomorphic function $\Omega$ on $\mathfrak H^n$ is called
a Hilbert modular form of weight $k$ with respect to $\Gamma$ if
$$
\Omega |_k \, \gamma =\Omega
$$
holds for every $\gamma\in \Gamma$, and usual conditions
at cusps when $F=\mathbf Q$. For every $g \in {\rm SL}(2, F)$,
$\Omega |_k \, g$ has a Fourier expansion of the form
$(\Omega |_k \, g)(z)=\sum _{\xi \in L} a_g(\xi )\mathbf e_F(\xi z)$, where $L$ is a lattice in $F$.
We have $a_g(\xi )=0$ if $\xi \neq 0$ is not totally positive.     
We call $\Omega$ {\it a cusp form} if the constant term $a_g(0)$ vanishes for every
$g \in {\rm SL}(2, F)$.
We denote the space of Hilbert modular forms (resp. cusp forms) of weight $k$
with respect to $\Gamma$ by
$M_k(\Gamma )=M_{k_1, k_2, \ldots , k_n}(\Gamma )$
(resp. $S_k(\Gamma )=S_{k_1, k_2, \ldots , k_n}(\Gamma )$).

Hereafter until the end of this subsection, we assume that
$\Gamma ={\rm SL}(2, \mathcal O_F)$ and $0 \neq \Omega \in S_k(\Gamma )$. 
The Fourier expansion of $\Omega$ takes the form
\begin{equation}
\Omega(z)=\sum _{0 \ll \xi \in \mathfrak d_F^{-1}} a(\xi )\mathbf e_F(\xi z).
\tag{2.1}
\end{equation}
Since
$\begin{pmatrix} u & 0 \\ 0 & u ^{-1} \end{pmatrix} \in \Gamma$
for $u \in E_F$, we have
$$
u^k \sum _{0 \ll \xi \in \mathfrak d_F^{-1}} a(\xi )\mathbf e_F(\xi u^2z)
=\sum _{0 \ll \xi \in \mathfrak d_F^{-1}} a(\xi )\mathbf e_F(\xi z),
$$
where we put $u^k=\prod _{\nu =1}^n (u^{(\nu )})^{k_{\nu}}$.
Therefore we have
\begin{equation}
a(u^2\xi )=u^ka(\xi ), \qquad u \in E_F.
\tag{2.2}
\end{equation}
In particular, taking $u=-1$, we have
\begin{equation}
\sum _{\nu =1}^n k_{\nu} \equiv 0 \mod 2.
\tag{2.3}
\end{equation}
For the sake of simplicity, we assume that
\begin{equation}
u^k>0 \qquad \text {for every $u \in E_F$}.
\tag{A}
\end{equation}
Put
$$
k_0=\max (k_1, k_2, \ldots , k_n), \qquad k_{\nu}^{\prime}=k_0-k_{\nu},
\qquad k^{\prime}=(k_1^{\prime}, k_2^{\prime}, \ldots , k_n^{\prime}).
$$
We define the $L$-function of $\Omega$ by
\begin{equation}
L(s, \Omega)=\sum _{\xi E_F^2} a(\xi )\xi ^{k^{\prime}/2}N(\xi )^{-s}, \qquad
\xi ^{k^{\prime}/2}=\prod _{\nu =1}^n (\xi ^{(\nu )})^{k_{\nu}^{\prime}/2}.
\tag{2.4}
\end{equation}
Here the summation extends over all cosets $\xi E_F^2$ with $\xi$ satisfying
$0 \ll \xi \in \mathfrak d_F^{-1}$. By (2.2) and (A), we see that the sum is well defined.
The series (2.4) converges when $\Re (s)$ is sufficiently large.
We put
\begin{equation}
R(s, \Omega)=(2\pi )^{-ns}\prod _{\nu =1}^n
\Gamma (s-\frac {k_{\nu}^{\prime}}{2}) L(s, \Omega).
\tag{2.5}
\end{equation}
By the standard calculation, we obtain the integral representation
\begin{equation}
 \int _{\mathbf R_+^n/E_F^2} \Omega(iy_1, iy_2, \ldots , iy_n)
\prod _{\nu =1}^n y_{\nu}^{s-k_{\nu}^{\prime}/2 -1}dy_{\nu}
=(2\pi )^{\sum _{\nu =1}^n k_{\nu}^{\prime}/2}R(s, \Omega)
\tag{2.6}
\end{equation}
when $\Re (s)$ is sufficiently large.
By a suitable transformation of this integral, we can show that
$R(s, \Omega)$ is an entire function of $s$ and satisfies the functional equation
\begin{equation}
R(s, \Omega)=(-1)^{\sum _{\nu =1}^n k_{\nu}/2}R(k_0-s, \Omega).
\tag{2.7}
\end{equation}

\medskip

{\bf 2.2.} In [Y3], Chapter V, \S5, we gave an explicit  method to attach a cohomology class
to a Hilbert modular form.  We will review it in this subsection. For
$0 \leq l \in \mathbf Z$ and 
$\begin{bmatrix} u \\ v \end{bmatrix} \in \mathbf C^2$, put
$$
\begin{bmatrix} u \\ v \end{bmatrix} ^l
={}^t (u^l \ u^{l-1}v \ldots uv^{l-1} \ v^l).
$$
Define a representation $\rho _l: GL(2, \mathbf C) \longrightarrow GL(l+1, \mathbf C)$ by
$$
\rho _l(g) \begin{bmatrix} u \\ v \end{bmatrix} ^l
=(g \begin{bmatrix} u \\ v \end{bmatrix} )^l.
$$
Let $\Gamma$ be a congruence subgroup of $SL(2, \mathcal O_F)$.
Let $l_1$, $l_2, \ldots , l_n$ be nonnegative integers.
Let $V$ be the representation space of
$\rho _{l_1} \otimes \rho _{l_2} \otimes \cdots \otimes \rho _{l_n}$.
Let $\Omega \in M_{l_1+2, l_2+2, \ldots , l_n+2}(\Gamma )$ be
a Hilbert modular form of weight $(l_1+2, l_2+2, \ldots , l_n+2)$.
Define a holomorphic $V$-valued $n$-form $\mathfrak d(\Omega)$ on $\mathfrak H^n$ by
\begin{equation}
\mathfrak d(\Omega)=\Omega(z) \begin{bmatrix} z_1 \\ 1 \end{bmatrix} ^{l_1} \otimes
\begin{bmatrix} z_2 \\ 1 \end{bmatrix} ^{l_2} \otimes \cdots \otimes
\begin{bmatrix} z_n \\ 1 \end{bmatrix} ^{l_n} dz_1dz_2 \cdots dz_n.
\tag{2.8}
\end{equation}
We put $\rho =\rho _{l_1} \otimes \rho _{l_2} \otimes \cdots \otimes \rho _{l_n}$.

Let $g=(g_1, \ldots , g_n) \in {\rm GL}(2, \mathbf R)_+^n$. Under the action of $g$
on $\mathfrak H^n$, $\mathfrak d(\Omega)$ transforms to $\mathfrak d(\Omega)\circ g$,
where
$$
\mathfrak d(\Omega)\circ g
=\Omega (g(z)) \begin{bmatrix} g_1z_1 \\ 1 \end{bmatrix} ^{l_1} \otimes
\cdots \otimes
\begin{bmatrix} g_nz_n \\ 1 \end{bmatrix} ^{l_n} (dz_1\circ g_1)
\cdots (dz_n\circ g_n).
$$
Since
$$
\begin{bmatrix} g_{\nu}z_{\nu} \\ 1 \end{bmatrix} ^{l_{\nu}}
=j(g_{\nu}, z_{\nu})^{-l_{\nu}}\rho _{l_{\nu}}(g_{\nu})
\begin{bmatrix} z_{\nu} \\ 1 \end{bmatrix} ^{l_{\nu}},
$$
$$
dz_{\nu}\circ g_{\nu}=(\det g_{\nu})j(g_{\nu}, z_{\nu})^{-2}dz_{\nu},
$$
we obtain
\begin{equation}
\mathfrak d(\Omega) \circ g
=\prod _{\nu =1}^n (\det g_{\nu})^{-l_{\nu}/2}
\rho (g)\mathfrak d(\Omega |_k \, g),
\quad g \in {\rm GL}(2, \mathbf R)_+^n \cap {\rm GL}(2, F).
\tag{2.9$_a$}
\end{equation}
In particular, we have
\begin{equation}
\mathfrak d(\Omega) \circ \gamma =\rho (\gamma )\mathfrak d(\Omega),
\qquad \gamma \in \Gamma .
\tag{2.9$_b$}
\end{equation}

We are going to discuss the case $n=2$ in detail. Take
$w=(w_1, w_2) \in \mathfrak H^2$. For $z=(z_1, z_2) \in \mathfrak H^2$, we put
\begin{equation}
F(z)=\int _{w_1}^{z_1} \int _{w_2}^{z_2} \mathfrak d(\Omega),
\tag{2.10}
\end{equation}
a period integral of Eichler--Shimura type.
Let $\mathcal H$ denote the vector space of all $V$-valued holomorphic
functions on $\mathfrak H^2$. For $\varphi \in \mathcal H$ and $\gamma \in \Gamma$,
we define a function $\gamma\varphi$ on $\mathfrak H^2$ by
\begin{equation}
(\gamma\varphi )(z)
=\rho (\gamma )\varphi (\gamma ^{-1}z).
\tag{2.11}
\end{equation}
Then $\mathcal H$ becomes a left $\Gamma$-module. Since

$$
\frac {\partial}{\partial z_1}\frac {\partial}{\partial z_2} (\gamma F-F)=0,
$$
we can write
$$
\gamma F-F=g(\gamma ; z_1)+h(\gamma ; z_2),
$$
where $g(\gamma ; z_1) \in \mathcal H$ (resp. $h(\gamma ; z_2) \in \mathcal H$)
is a function which depends only on $z_1$ (resp. $z_2$) (cf. [Y3], p. 208, Lemma 5.1).
We regard $g$ and $h$ as $1$-cochains in $C^1(\Gamma , \mathcal H)$.
Then clearly we have ($d_1$ in \S1.1 is abbreviated to $d$)
$$
dg(\gamma _1, \gamma _2; z_1)+dh(\gamma _1, \gamma _2; z_2)=0.
$$
Put
$$
f (\Omega)(\gamma _1, \gamma _2)=dg(\gamma _1, \gamma _2; z_1).
$$
We abbreviate $f (\Omega)$ to $f$.
We see that $f (\gamma _1, \gamma _2) \in V$ is a constant.
Furthermore, in $\mathcal H$, $f$ is a coboundary.
Hence $f$ satisfies the cocycle condition
\begin{equation}
\gamma _1f (\gamma _2, \gamma _3)-f (\gamma _1\gamma _2, \gamma _3)
+f (\gamma _1, \gamma_2\gamma _3)-f (\gamma _1, \gamma _2)=0.
\tag{2.12}
\end{equation}
The $2$-cocycle $f$ determines a cohomology class in
$H^2(\Gamma , V)$.

Let us give an explicit formula for $f$.
For $x \in F$, let $x^{\prime}$ denote the conjugate of $x$ over $\mathbf Q$.
For $\gamma = \begin{pmatrix} a & b \\ c & d \end{pmatrix} \in \Gamma$, let
$\gamma ^{\prime}=\begin{pmatrix} a^{\prime} & b^{\prime} \\
c^{\prime} & d^{\prime} \end{pmatrix}$.
We regard $\gamma$, $\gamma ^{\prime} \in SL(2, \mathbf R)$.
Then, for $\gamma \in \Gamma$, we have
\begin{equation*}
\begin{aligned}
& F(\gamma (z))=F(\gamma z_1, \gamma ^{\prime}z_2)
=\int _{w_1}^{\gamma z_1} \int _{w_2}^{\gamma ^{\prime}z_2} \mathfrak d(\Omega) \\
= & \int _{\gamma w_1}^{\gamma z_1}
\int _{\gamma ^{\prime}w_2}^{\gamma ^{\prime}z_2} \mathfrak d(\Omega)
+\int _{\gamma w_1}^{\gamma z_1} \int _{w_2}^{\gamma ^{\prime}w_2} \mathfrak d(\Omega)
+\int _{w_1}^{\gamma w_1} \int _{w_2}^{\gamma ^{\prime}z_2} \mathfrak d(\Omega) \\
= & (\rho _{l_1}(\gamma ) \otimes \rho _{l_2}(\gamma ^{\prime}))F(z)
+\int _{\gamma w_1}^{\gamma z_1} \int _{w_2}^{\gamma ^{\prime}w_2} \mathfrak d(\Omega)
+\int _{w_1}^{\gamma w_1} \int _{w_2}^{\gamma ^{\prime}z_2} \mathfrak d(\Omega) ,
\end{aligned}
\end{equation*}
Substituting $z$ by $\gamma ^{-1}z$ in this formula, we get
\begin{equation}
(\rho _{l_1}(\gamma ) \otimes \rho _{l_2}(\gamma ^{\prime}))F(\gamma ^{-1}z)-F(z)
=-\int _{\gamma w_1}^{z_1} \int _{w_2}^{\gamma ^{\prime}w_2} \mathfrak d(\Omega)
-\int _{w_1}^{\gamma w_1} \int _{w_2}^{z_2} \mathfrak d(\Omega) .
\tag{2.13}
\end{equation}
We may take
\begin{equation}
g(\gamma ; z_1)=
-\int _{\gamma w_1}^{z_1} \int _{w_2}^{\gamma ^{\prime}w_2} \mathfrak d(\Omega),
\tag{2.14}
\end{equation}
\begin{equation}
h(\gamma ; z_2)
=-\int _{w_1}^{\gamma w_1} \int _{w_2}^{z_2} \mathfrak d(\Omega) .
\tag{2.15}
\end{equation}
For $\gamma _1$, $\gamma _2 \in \Gamma$, we have
\begin{equation}
f (\gamma _1, \gamma _2)=(\gamma _1 g)(\gamma _2; z_1)
-g(\gamma _1\gamma _2; z_1)+g(\gamma _1; z_1),
\tag{2.16}
\end{equation}
\begin{equation}
f (\gamma _1, \gamma _2)=- \{(\gamma _1 h)(\gamma _2; z_2)
-h(\gamma _1\gamma _2; z_2)+h(\gamma _1; z_2) \} .
\tag{2.17}
\end{equation}
By (2.14) and (2.16), we have
\begin{equation*}
\begin{aligned}
f (\gamma _1, \gamma _2) = &
(\rho _{l_1}(\gamma _1) \otimes \rho _{l_2}(\gamma _1^{\prime}))
g(\gamma _2; \gamma _1^{-1}z_1)
-g(\gamma _1\gamma _2; z_1)+g(\gamma _1; z_1) \\
= & -(\rho _{l_1}(\gamma _1) \otimes \rho _{l_2}(\gamma _1^{\prime}))
\int _{\gamma _2w_1}^{\gamma _1^{-1}z_1}
\int _{w_2}^{\gamma _2^{\prime}w_2} \mathfrak d(\Omega) \\
& + \int _{\gamma _1\gamma _2w_1}^{z_1}
\int _{w_2}^{\gamma _1^{\prime}\gamma _2^{\prime}w_2} \mathfrak d(\Omega)
-\int _{\gamma _1w_1}^{z_1}
\int _{w_2}^{\gamma _1^{\prime}w_2} \mathfrak d(\Omega) \\
= &- \int _{\gamma _1\gamma _2w_1}^{z_1}
\int _{\gamma _1^{\prime}w_2}^{\gamma _1^{\prime}\gamma _2^{\prime}w_2} \mathfrak d(\Omega)
+ \int _{\gamma _1\gamma _2w_1}^{z_1}
\int _{w_2}^{\gamma _1^{\prime}\gamma _2^{\prime}w_2} \mathfrak d(\Omega)
-\int _{\gamma _1w_1}^{z_1}
\int _{w_2}^{\gamma _1^{\prime}w_2} \mathfrak d(\Omega) \\
= & \int _{\gamma _1\gamma _2w_1}^{z_1}
\int _{w_2}^{\gamma _1^{\prime}w_2} \mathfrak d(\Omega)
-\int _{\gamma _1w_1}^{z_1}
\int _{w_2}^{\gamma _1^{\prime}w_2} \mathfrak d(\Omega) \\
= & \int _{\gamma _1\gamma _2w_1}^{\gamma _1w_1}
\int _{w_2}^{\gamma _1^{\prime}w_2} \mathfrak d(\Omega)
\end{aligned}
\end{equation*}
using (2.9$_b$). Thus we obtain an explicit formula
\begin{equation}
f (\gamma _1, \gamma _2)
=\int _{\gamma _1\gamma _2w_1}^{\gamma _1w_1}
\int _{w_2}^{\gamma _1^{\prime}w_2} \mathfrak d(\Omega).
\tag{2.18}
\end{equation}
By (2.9$_b$), (2.18) can be written as
\begin{equation}
f (\gamma _1, \gamma _2)
=(\rho _{l_1}(\gamma _1) \otimes \rho _{l_2}(\gamma _1^{\prime}))
\int _{\gamma _2w_1}^{w_1}
\int _{\gamma _1^{\prime -1}w_2}^{w_2} \mathfrak d(\Omega).
\tag{2.19}
\end{equation}

Suppose that $w_1$ is replaced by $w_1^{\ast}$, $w_2$ remaining the same.
Then $g(\gamma ; z_1)$ changes to $g(\gamma , z_1)+a(\gamma )$, where
$$
a(\gamma )=\int _{\gamma w_1}^{\gamma w_1^{\ast}}
\int _{w_2}^{\gamma ^{\prime}w_2} \mathfrak d(\Omega).
$$
Hence $f (\gamma _1, \gamma _2)$ changes to
$f (\gamma _1, \gamma _2)+\gamma _1a(\gamma _2)
-a(\gamma _1\gamma _2)+a(\gamma _1)$.
Suppose that $w_2$ is replaced by $w_2^{\ast}$, $w_1$ remaining the same.
Then $h(\gamma ; z_2)$ changes to $h(\gamma , z_2)+b(\gamma )$, where
$$
b(\gamma )=\int _{w_1}^{\gamma w_1}\int _{w_2}^{w_2^{\ast}} \mathfrak d(\Omega).
$$
By (2.17), $f (\gamma _1, \gamma _2)$ changes to
$f (\gamma _1, \gamma _2)-\gamma _1b(\gamma _2)
+b(\gamma _1\gamma _2)-b(\gamma _1)$.
Therefore the cohomology class of $f$ does not depend on
the choice of the \lq\lq base points" $w_1$, $w_2$.

Put $\overline{\Gamma} =\Gamma /(\{ \pm 1_2 \} \cap \Gamma )$.
By (2.18), we see that $f$ can be regarded as a $2$-cocycle
of $\overline{\Gamma}$ taking values in $V$. Depending on the context,
we consider $f$ as a $2$-cocycle on $\overline{\Gamma}$. We see that
the cocycle $f$ is normalized, i.e.,
\begin{equation}
f (1, \gamma )=f (\gamma , 1)=0
\qquad \text {for every} \quad \gamma \in \overline{\Gamma} .
\tag{2.20}
\end{equation}

Now assume that $\Omega$ is a cusp form. Then the cocycle $f =f (\Omega)$ satisfies
the \lq\lq parabolic condition". 
Namely let $q \in \Gamma$ be a parabolic element
and $w^{\ast}=(w_1^{\ast}, w_2^{\ast})$ be the fixed point of $q^{\prime}$.
Since $f$ is a cusp form, we may replace $w_2$ by $w_2^{\ast}$.
\footnote{For every $g \in {\rm SL}(2, F)$, we have the Fourier expansion
$(\Omega |_k \, g)(z)=\sum _{0 \ll \xi \in L} a_g(\xi )\mathbf e_F(\xi z)$
where $L$ is a lattice in $F$. We have the estimate
$\vert a_g(\xi )\vert \leq M\xi ^{k_1/2}\xi ^{\prime k_2/2}$
with a positive constant $M$ depending on $\Omega$ and $g$
(cf. [Sh7], p. 280, Proposition A6.4).
Using this estimate, it is not difficult to check the absolute convergence
of the integral (2.10) defining $F(z)$ when $w_2$ is replaced by $w_2^{\ast}$.}
Let $f ^{\ast}$ be the cocycle obtained from $(w_1, w_2^{\ast})$. We have
$$
f ^{\ast}(\gamma _1, \gamma _2)=f (\gamma _1, \gamma _2)
-\gamma _1b(\gamma _2)+b(\gamma _1\gamma _2)-b(\gamma _1)
$$
with a $1$-cochain $b$ and $f ^{\ast}(q, \gamma )=0$.
Therefore
$$
f (q, \gamma )
=qb(\gamma )-b(q\gamma )+b(q ), \qquad \gamma \in \Gamma ,
$$
i.e., $f (q, \gamma )$ is of the form of coboundary whenever
$q$ is parabolic. Similar argument applies to $f (\gamma ,q)$.

\medskip

{\bf 2.3.} We are going to investigate closely the relation between
the critical values of $L$-function $L(s, \Omega)$ and the cocycle $f (\Omega)$.
Until the end of this subsection, we assume $\Gamma =SL(2, \mathcal O_F)$.
Let $\epsilon$ be the fundamental unit of $F$. We put
$$
\sigma = \begin{pmatrix} 0 & 1 \\ -1 & 0 \end{pmatrix} , \qquad
\mu = \begin{pmatrix} \epsilon & 0 \\ 0 & \epsilon ^{-1} \end{pmatrix} .
$$
We regard $\sigma$ and $\mu$ as elements of $\overline\Gamma$.
Taking $\gamma _1=\gamma _2=\gamma _3=\sigma$ in (2.12),
we obtain
\begin{equation}
\sigma f (\sigma , \sigma )=f (\sigma , \sigma )
\tag{2.21}
\end{equation}
in view of (2.20). As the base points, we choose
$$
w_1=i\epsilon ^{-1}, \qquad w_2=i\infty .
$$
By (2.18), we get
\begin{equation}
f (\sigma , \mu )=f (\sigma , \sigma )
=-\int _{i\epsilon ^{-1}}^{i\epsilon} \int _0^{i\infty} \mathfrak d(\Omega).
\tag{2.22}
\end{equation}
Put
$$
P=\left\{ \begin{pmatrix} u & v \\ 0 & u^{-1} \end{pmatrix} \bigg|
\ u \in E_F, \ v \in \mathcal O_F \right\} \subset \Gamma .
$$
By (2.18), we get
\begin{equation}
f (p, \gamma )=0 \qquad \text {for every} \quad p \in P, \ \gamma \in \Gamma
\tag{2.23}
\end{equation}
since we have $pw_2=w_2$ for $p \in P$. 
Taking $\gamma _1=p \in P$ in (2.12), we obtain
\begin{equation}
f (p\gamma _1, \gamma _2)=pf (\gamma _1, \gamma _2)
\qquad \text {for every} \quad p \in P, \ \gamma _1, \gamma _2 \in \Gamma .
\tag{2.24}
\end{equation}
This is the parabolic condition for $\Gamma ={\rm SL}(2, \mathcal O_F)$
and will play a crucial role in the succeeding sections.

For $0 \leq s \leq l_1$, $0 \leq t \leq l_2$, we put
\begin{equation}
P_{s,t}=\int _{i\epsilon ^{-1}}^{i\epsilon}\int _0^{i\infty}
\Omega(z)z_1^sz_2^t dz_1dz_2.
\tag{2.25}
\end{equation}
The components of $f (\sigma , \sigma )$ are given by
$-P_{s, t}$. The condition
$\sigma f (\sigma , \sigma )=f (\sigma , \sigma )$
is equivalent to
\begin{equation}
P_{s, t}=(-1)^{l_1+l_2-s-t}P_{l_1-s, l_2-t}.
\tag{2.26}
\end{equation}
Put $k_1=l_1+2$, $k_2=l_2+2$. By (2.3), we have
\begin{equation}
l_1 \equiv l_2 \mod 2.
\tag{2.27}
\end{equation}
We assume that $l_1 \geq l_2$. Then we have
$$
k_0=k_1, \qquad k_1^{\prime}=0, \qquad k_2^{\prime}=k_1-k_2.
$$
Since $E_F^2=\langle \epsilon ^2\rangle$, a fundamental domain of
$\mathbf R_+^2/E_F^2$ is given by
$[\epsilon ^{-1}, \epsilon ]\times\mathbf R_+$.
By (2.6), we obtain
\begin{equation}
\int _{\epsilon ^{-1}}^{\epsilon} \int _0^{\infty}
\Omega(iy_1, iy_2)y_1^{s-1}y_2^{s-(k_1-k_2)/2-1}dy_1dy_2
=(2\pi )^{(k_1-k_2)/2}R(s, \Omega)
\tag{2.28}
\end{equation}
when $\Re (s)$ is sufficiently large. We can verify that the integral converges
locally uniformly for $s \in \mathbf C$.
Take $m \in \mathbf Z$ and put $s=m$, $t=m-(k_1-k_2)/2$.
Then $0 \leq s\leq l_1$, $0 \leq t \leq l_2$ hold
if and only if
\begin{equation}
\frac {k_1-k_2}{2} \leq m \leq \frac {k_1+k_2}{2}-2.
\tag{2.29}
\end{equation}
For an integer $m$ in this range, we have
\begin{equation*}
\begin{aligned}
P_{m, m-(k_1-k_2)/2}
= & \int _{i\epsilon ^{-1}}^{i\epsilon}\int _0^{i\infty}
\Omega(z)z_1^mz_2^{m-(k_1-k_2)/2} dz_1dz_2 \\
= \, & i^{2m-(k_1-k_2)/2+2}
\int _{\epsilon ^{-1}}^{\epsilon}\int _0^{\infty}
\Omega(iy_1, iy_2)y_1^my_2^{m-(k_1-k_2)/2}dy_1dy_2.
\end{aligned}
\end{equation*}
Therefore we obtain
\begin{equation}
P_{m, m-(k_1-k_2)/2}=(-1)^{m+1}i^{-(k_1-k_2)/2}
(2\pi )^{(k_1-k_2)/2}R(m+1, \Omega)
\tag{2.30}
\end{equation}
by (2.28). By the functional equation (2.7), this is equal to
$$
(-1)^{m+1}i^{-(k_1-k_2)/2}
(2\pi )^{(k_1-k_2)/2}(-1)^{(k_1+k_2)/2}R(k_1-m-1, \Omega).
$$
Since $k_1-m-2$ satisfies (2.29), we obtain
\begin{equation}
P_{m, m-(k_1-k_2)/2}
=(-1)^{(k_1-k_2)/2}P_{k_1-m-2, (k_1+k_2)/2-m-2}
\tag{2.31}
\end{equation}
using (2.30). We see that (2.31) is consistent with (2.26).
Note that (2.29) is the condition for $L(m+1, \Omega)$ to be a critical value
(cf. [Sh4], (4.14)).

\medskip

{\bf 2.4.} Let $\Omega \in M_k(\Gamma )$ and $f =f (\Omega) \in Z^2(\Gamma , V)$
be the $2$-cocycle attached to $\Omega$ defined by (2.18).
In this subsection, we will write the action of Hecke operators on the cohomology
class of $f (\Omega)$ explicitly.
We denote $f(\Omega )$ also by $f_{\Omega}$.

Let $F$ be a totally real number field of degree $n$ and
$\Gamma$ be a congruence subgroup of ${\rm SL}(2, \mathcal O_F)$.
Let $\varpi$ be a totally positive element of $F$ and let
\begin{equation}
\Gamma \begin{pmatrix} 1 & 0 \\ 0 & \varpi \end{pmatrix}\Gamma
=\sqcup _{i=1}^d \Gamma\beta _i
\tag{2.32}
\end{equation}
be a coset decomposition. Let $\Omega \in M_k(\Gamma )$.
We define the Hecke operator $T(\varpi )$ by

\begin{equation}
\Omega \mid T(\varpi )=N(\varpi )^{k_0/2-1}
\sum _{i=1}^d \Omega |_k \, \beta _i.
\tag{2.33}
\end{equation}
Clearly $T(\varpi )$ does not depend on the choice of
the coset decomposition (2.32). We have $\Omega\vert T(\varpi ) \in M_k(\Gamma )$;
it is a cusp form if $\Omega$ is. By (2.9$_a$), we have
\begin{equation}
\mathfrak d(\Omega \mid T(\varpi ))
= \prod _{\nu =1}^n (\varpi ^{(\nu })^{(k_0+k_{\nu})/2-2}
\sum _{i=1}^d \rho (\beta _i)^{-1}(\mathfrak d(\Omega)\circ \beta _i).
\tag{2.34}
\end{equation}
Put
\begin{equation}
c=\prod _{\nu =1}^n (\varpi ^{(\nu )})^{(k_0+k_{\nu})/2-2}.
\tag{2.35}
\end{equation}
Until the end of this section, we assume that $n=2$.
We define (cf. (2.10))
\begin{equation}
F_{\Omega \mid T(\varpi )}(z)
=\int _{w_1}^{z_1} \int _{w_2}^{z_2} \mathfrak d(\Omega \mid T(\varpi )),
\qquad z=(z_1, z_2).
\tag{2.36}
\end{equation}
By (2.34), we have
\begin{equation*}
\begin{aligned}
F_{\Omega \mid T(\varpi )}(z)
= \, & c\sum _{i=1}^d \beta _i^{-1} \int _{w_1}^{z_1}\int _{w_2}^{z_2}
\mathfrak d(\Omega)\circ \beta _i
=c\sum _{i=1}^d \beta _i^{-1}
\int _{\beta _iw_i}^{\beta _iz_1}\int _{\beta _i^{\prime}w_2}^{\beta _i^{\prime}z_2}
\mathfrak d(\Omega) \\
= \, & c\sum _{i=1}^d \beta _i^{-1}\bigg[
\int _{w_1}^{\beta _iz_1}\int _{w_2}^{\beta _i^{\prime}z_2}\mathfrak d(\Omega)
-\int _{w_1}^{\beta _iz_1}\int _{w_2}^{\beta _i^{\prime}w_2}\mathfrak d(\Omega) \\
& \hskip 1.5em -\int _{w_1}^{\beta _iw_1}\int _{w_2}^{\beta _i^{\prime}z_2}\mathfrak d(\Omega)
+\int _{w_1}^{\beta _iw_1}\int _{w_2}^{\beta _i^{\prime}w_2}\mathfrak d(\Omega) \bigg] .
\end{aligned}
\end{equation*}
Here, for simplicity, we write the action $\rho (\beta _i)^{-1}$ by $\beta _i^{-1}$.
Therefore we obtain
\begin{equation}
\begin{aligned}
F_{\Omega \mid T(\varpi )}(z)
= c \sum _{i=1}^d \beta _i^{-1}\bigg[
& F(\beta _iz)-F(\beta _i(z_1, w_2)) \\
- & F(\beta _i(w_1, z_2))+F(\beta _i(w_1, w_2))\bigg] .
\end{aligned}
\tag{2.37}
\end{equation}
Take $\gamma \in \Gamma$. We have
$$
(\gamma F)(z)-F(z)=g_\Omega(\gamma ; z_1)+h_\Omega(\gamma ; z_2)
$$
with $g_\Omega=g$ (resp. $h_\Omega=h$) defined by (2.14) (resp. (2.15)).
Similarly to this formula, we have a decomposition
\begin{equation}
(\gamma F_{\Omega\mid T(\varpi )})(z)-F_{\Omega\mid T(\varpi )}(z)
=g_{\Omega\mid T(\varpi )}(\gamma ; z_1)+h_{\Omega\mid T(\varpi )}(\gamma ; z_2).
\tag{2.38}
\end{equation}
Here $g_{\Omega\mid T(\varpi )}(\gamma ; z_1)$ (resp. $h_{\Omega\mid T(\varpi )}(\gamma ; z_2)$)
is a $V$-valued holomorphic function on $\mathfrak H^2$ which depends only on
$z_1$ (resp. $z_2$).
If (2.38) holds, then a $2$-cocycle $f _{\Omega \mid T(\varpi )}$
attached to $\Omega \mid T(\varpi )$ is given by (cf. (2.16))
\begin{equation}
f _{\Omega \mid T(\varpi )}(\gamma _1, \gamma _2)
=(\gamma _1g_{\Omega \mid T(\varpi )})(\gamma _2; z_1)
-g_{\Omega \mid T(\varpi )}(\gamma_1\gamma _2; z_1)
+g_{\Omega \mid T(\varpi )}(\gamma _1; z_1).
\tag{2.39}
\end{equation}
Suppose that we have a decomposition
\begin{equation}
c\sum _{i=1}^d (\gamma\beta _i^{-1}F(\beta _i\gamma ^{-1}z) -\beta _i^{-1}F(\beta _iz))
=g_{\Omega\mid T(\varpi )}^{\ast}(\gamma ; z_1)+h_{\Omega\mid T(\varpi )}^{\ast}(\gamma ; z_2).
\tag{2.40}
\end{equation}
Then by (2.37), we see that $g_{\Omega\mid T(\varpi )}(\gamma ; z_1)$ can be taken in the form
$$
g_{\Omega\mid T(\varpi )}(\gamma ; z_1)=g_{\Omega\mid T(\varpi )}^{\ast}(\gamma ; z_1)
+\gamma q(\gamma ^{-1}z_1)-q(z_1)+x(\gamma ).
$$
Here $q(z_1)$ is a $V$-valued holomorphic function which depends only on $z_1$
and does not depend on $\gamma$ and $x(\gamma ) \in V$.
Therefore, if we substract the coboundary which comes from $x(\gamma )$
from the cocycle $f _{\Omega\mid T(\varpi )}$, the resulting cocycle can be
calculated using (2.39) with $g_{\Omega\mid T(\varpi )}^{\ast}(\gamma ; z_1)$ in place of
$g_{\Omega\mid T(\varpi )}(\gamma ; z_1)$.
Now we put
$$
\beta _i\gamma =\delta _i\beta _{j(i)}, \qquad 1 \leq i \leq d, \quad \delta _i \in \Gamma .
$$
Set
\begin{equation}
g_{\Omega\mid T(\varpi )}^{\ast}(\gamma ; z_1)
=c\sum _{i=1}^d \beta _i^{-1}g_\Omega(\delta _i; \beta _iz_1),
\tag{2.41}
\end{equation}
\begin{equation}
h_{\Omega\mid T(\varpi )}^{\ast}(\gamma ; z_2)
=c\sum _{i=1}^d \beta _i^{-1}h_\Omega(\delta _i; \beta _i^{\prime}z_2).
\tag{2.42}
\end{equation}
Since $i \mapsto j(i)$ is a permutation on $d$ letters, we have
\begin{equation*}
\begin{aligned}
& \sum _{i=1}^d (\gamma\beta _{j(i)}^{-1}F(\beta _{j(i)}\gamma ^{-1}z) 
-\beta _i^{-1}F(\beta _iz))
=\sum _{i=1}^d (\beta _i^{-1}\delta _iF(\delta _i^{-1}\beta _iz) 
-\beta _i^{-1}F(\beta _iz)) \\
= & \sum _{i=1}^d \beta _i^{-1}
(g_\Omega(\delta _i; \beta _iz_1)+h_\Omega(\delta _i; \beta _i^{\prime}z_2)).
\end{aligned}
\end{equation*}
Hence we see that (2.40) holds.

For $\gamma _1$, $\gamma _2 \in \Gamma$, we put
\begin{equation}
\beta _i\gamma _1=\delta _i^{(1)}\beta _{j(i)}, \ \delta _i^{(1)} \in \Gamma , \quad
\beta _i\gamma _2=\delta _i^{(2)}\beta _{k(i)}, \ \delta _i^{(2)} \in \Gamma ,
\qquad 1 \leq i \leq d.
\tag{2.43}
\end{equation}
We have
$$
\beta _i\gamma _1\gamma _2=\delta _i^{(1)}\delta _{j(i)}^{(2)}\beta _{k(j(i))}.
$$
Now we calculate $f _{\Omega\mid T(\varpi )}$
using (2.39) with $g_{\Omega\mid T(\varpi )}^{\ast}(\gamma ; z_1)$ in place of
$g_{\Omega\mid T(\varpi )}(\gamma ; z_1)$. Then we have
\begin{equation*}
\begin{aligned}
& f _{\Omega \mid T(\varpi )}(\gamma _1, \gamma _2)
=\gamma _1g_{\Omega \mid T(\varpi )}^{\ast}(\gamma _2; \gamma _1^{-1}z_1)
-g_{\Omega \mid T(\varpi )}^{\ast}(\gamma_1\gamma _2; z_1)
+g_{\Omega \mid T(\varpi )}^{\ast}(\gamma _1; z_1) \\
= \, & c \bigg[ \sum _{i=1}^d \gamma _1\beta _i^{-1}
g_\Omega(\delta _i^{(2)}; \beta _i\gamma _1^{-1}z_1)
-\beta _i^{-1}g_\Omega(\delta _i^{(1)}\delta _{j(i)}^{(2)}; \beta _iz_1)
+\beta _i^{-1}g_\Omega(\delta _i^{(1)}; \beta _iz_1)\bigg] \\
= \, & c \bigg[ \sum _{i=1}^d \beta _i^{-1}\delta _i^{(1)}
g_\Omega(\delta _{j(i)}^{(2)}; \beta _{j(i)}\gamma _1^{-1}z_1)
-\beta _i^{-1}g_\Omega(\delta _i^{(1)}\delta _{j(i)}^{(2)}; \beta _iz_1)
+\beta _i^{-1}g_\Omega(\delta _i^{(1)}; \beta _iz_1)\bigg] \\
= \, & c\sum _{i=1}^d \beta _i^{-1}\bigg[ \delta _i^{(1)}
g_\Omega(\delta _{j(i)}^{(2)}; (\delta _i^{(1)})^{-1}\beta _iz_1)
-g_\Omega(\delta _i^{(1)}\delta _{j(i)}^{(2)}; \beta _iz_1)
+g_\Omega(\delta _i^{(1)}; \beta _iz_1)\bigg] .
\end{aligned}
\end{equation*}
Therefore we obtain an explicit formula
\begin{equation}
f _{\Omega \mid T(\varpi )}(\gamma _1, \gamma _2)
=c\sum _{i=1}^d \beta _i^{-1}f _\Omega(\delta _i^{(1)}, \delta _{j(i)}^{(2)}).
\tag{2.44}
\end{equation}
This formula can be written as
\begin{equation}
f _{\Omega \mid T(\varpi )}(\gamma _1, \gamma _2)
=c\sum _{i=1}^d \beta _i^{-1}
f _\Omega(\beta _i\gamma _1\beta _{j(i)}^{-1},
\beta _{j(i)}\gamma _2\beta _{k(j(i))}^{-1})
\tag{2.45}
\end{equation}
and is consistent with Proposition 1.3.

\medskip

{\bf 2.5.} Assume that the class number of $F$ in the narrow sense is $1$.
Suppose that $\Omega$ is a Hecke eigenform. Then the $L$-function
$L(s, \Omega )$ defined by (2.4) essentially coincides with the Euler product
given in [Sh4] or in Jacquet-Langlands [JL] but there is a subtle diffence;
we are going to explain it briefly for the reader's convenience.

We write $\mathfrak d_F=(\delta )$ with $\delta \gg 0$.
Let $\Omega \in S_{k_1, k_2}(\Gamma )$, $\Gamma ={\rm SL}(2, \mathcal O_F)$
and let
$$
\Omega (z)=\sum _{0 \ll \alpha \in \mathcal O_F}
c(\alpha )\mathbf e_F(\frac {\alpha}{\delta}z)
$$
be the Fourier expansion. We have $a(\alpha /\delta )=c(\alpha )$ (cf. (2.1)).
We set
$$
\Delta =\{ \alpha \in M(2, \mathcal O_F) \mid \det \alpha \gg 0 \} .
$$
Let $\mathfrak m$ be an integral ideal of $F$ and take $m \gg 0$ so that
$\mathfrak m=(m)$. Then we define
$$
T(\mathfrak m)=\sum _{\alpha \in \Delta , \, \det \alpha =m}
\Gamma\alpha\Gamma ,
$$
which is an element of the abstract Hecke ring $\mathcal H(\Gamma , \Delta )$
(cf. [Sh2], p. 54). Let
$T(\mathfrak m)=\sqcup _{i=1}^e \Gamma\beta _i$ be a coset decomposition.
Assume that $k_1 \geq k_2$. We define the action of
$T(\mathfrak m)$ on $\Omega$ by
$$
\Omega \mid T(\mathfrak m)=N(\mathfrak m)^{k_1/2-1}
\sum _{i=1}^e \Omega |_k \beta _i.
$$
Then $\Omega\vert T(\mathfrak m) \in S_k(\Gamma )$; we can verify
easily that it does not depend on the choices of $m$ and $\beta _i$. We have
$$
T(\mathfrak m)=\sqcup _{(d), d \gg 0 , ad=m} \ \sqcup _{b \!\!\! \mod d}
\ \Gamma \begin{pmatrix} a & b \\ 0 & d \end{pmatrix} .
$$
Calculating similarly to [Sh2], p. 79-80, we find that the Fourier expansion
of $\Omega\vert T(\mathfrak m)$ is given by
$$
(\Omega\vert T(\mathfrak m)) (z)=\sum _{0 \ll \alpha \in \mathcal O_F}
c^{\prime}(\alpha )\mathbf e_F(\frac {\alpha}{\delta}z),
$$
\begin{equation}
c^{\prime}(\alpha )=(m^{(2)})^{(k_1-k_2)}
\sum _{(a), a \gg 0, a \mid (m, \alpha )}
N(a)^{k_1-1}(a^{(2)})^{k_2-k_1} c(\frac {m\alpha}{a^2}).
\tag{2.46}
\end{equation}
Now assume that $\Omega$ is a nonzero common eigenfunction for all Hecke
operators $T(\mathfrak m)$. We put
$$
\Omega \mid T(\mathfrak m)=\lambda (\mathfrak m)\Omega .
$$
By (2.46), we find $c(1) \neq 0$. We assume that $\Omega$ is normalized
so that $c(1)=1$. Then we have
$$
\lambda (\mathfrak m)=c(m)(m^{(2)})^{(k_1-k_2)/2}.
$$
Using this formula and (2.46), we obtain, after routine computations, that
\begin{equation}
L(s, \Omega )=(\delta ^{(2)})^{(k_1-k_2)/2}D_F^s
\prod _{\mathfrak p} (1-\lambda (\mathfrak p)N(\mathfrak p)^{-s}
+N(\mathfrak p)^{k_1-1-2s})^{-1}.
\tag{2.47}
\end{equation}
Here $\mathfrak p$ extends over all prime ideals of $F$ and
$D_F=N(\delta )$ is the discriminant of $F$.

When $0 \ll \varpi \in \mathcal O_F$ generates a prime ideal $\mathfrak p$,
we denote $T(\varpi )$ defined by (2.33) also by $T(\mathfrak p)$.

\bigskip

\centerline {\S3. Cohomology of $P$}

\medskip

In this section, we will study cohomology groups of $P$.
Main results are Theorems 3.7 and 3.9 which give the vanishing
of $H^1(P, V)$ and $H^2(P,V)$ when $l_1 \neq l_2$. Hereafter in this paper,
we assume that $[F:{\mathbf Q}]=2$. We also assume
$l_1 \equiv l_2\mod 2$ and $l_2 \leq l_1$.

\medskip

{\bf 3.1.} Put $\Gamma ={\rm PSL}(2, \mathcal O_F)$.
In this section, we define subgroups $P$ and $U$ of $\Gamma$ by
\begin{equation*}
\begin{aligned}
& P=\left\{ \begin{pmatrix} t & 0 \\ u & t^{-1} \end{pmatrix} \bigg|
\ t \in E_F, \ u \in \mathcal O_F \right\}/\{ \pm 1_2 \} , \\
& U=\left\{ \begin{pmatrix} \pm 1 & 0 \\ u & \pm 1 \end{pmatrix} \bigg|
\ u \in \mathcal O_F \right\} /\{ \pm 1_2 \} .
\end{aligned}
\end{equation*}
We write $\mathcal O_F={\mathbf Z}+{\mathbf Z}\omega$.
Let $\epsilon$ be the fundamental unit of $F$ and let
$$
\epsilon ^2=A+B\omega , \qquad \epsilon ^2\omega =C+D\omega .
$$
Then we see that $\epsilon ^2$ is an eigenvalue of
$\begin {pmatrix} A & B \\ C & D \end{pmatrix}$ and that
$\begin {pmatrix} A & B \\ C & D \end{pmatrix} \in {\rm SL}(2, {\mathbf Z})$.
We put
$$
u_1=\begin{pmatrix} 1 & 0 \\ 1 & 1 \end{pmatrix}, \qquad
u_2=\begin{pmatrix} 1 & 0 \\ \omega & 1 \end{pmatrix} \in U, \qquad
t=\begin{pmatrix} \epsilon ^{-1} & 0 \\ 0 & \epsilon \end{pmatrix}.
$$
We have
\begin{equation}
tu_1t^{-1}=u_1^Au_2^B, \qquad tu_2t^{-1}=u_1^Cu_2^D.
\tag{3.1}
\end{equation}
We put
\begin{equation}
\mathcal Z=\{ (U_1, U_2) \in V \times V \mid
(u_1-1)U_2=(u_2-1)U_1 \} .
\tag{3.2}
\end{equation}
It is easy to see that by the mapping
$$
Z^1(U, V) \ni f \longrightarrow (f(u_1), f(u_2)) \in \mathcal Z,
$$
we have an isomorphism $Z^1(U, V) \cong \mathcal Z$.
Put
\begin{equation}
\mathcal B=\{ ((u_1-1)\mathbf b, (u_2-1)\mathbf b) \mid \mathbf b \in V \} .
\tag{3.3}
\end{equation}
Then we have
$B^1(U, V) \cong \mathcal B \subset \mathcal Z$.

We have $V=V_1 \otimes V_2$, $V_1={\mathbf C}^{l_1+1}$, $V_2={\mathbf C}^{l_2+1}$.
Let $\{ \mathbf e_1, \mathbf e_2, \ldots , \mathbf e_{l_1+1} \}$
(resp. $\{ \mathbf e_1^{\prime}, \mathbf e_2^{\prime}, \ldots , \mathbf e_{l_2+1}^{\prime} \}$)
be the standard basis of $V_1$ (resp. $V_2$).

\proclaim Lemma 3.1. We have $\dim V^U=1$ and $V^U$ is spanned by
$\mathbf e_{l_1+1} \otimes \mathbf e_{l_2+1}^{\prime}$.
\par
{\bf Proof.} Put
$$
\widetilde U=\bigg\{ (\begin{pmatrix} 1 & 0  \\ c_1 & 1 \end{pmatrix} ,
\begin{pmatrix} 1 & 0  \\ c_2 & 1 \end{pmatrix}) \, \bigg|
\ c_1, c_2 \in \mathbf C \bigg\} \subset {\rm SL}(2, \mathbf C)^2,
$$
$$
\widetilde U_1=\bigg\{ (\begin{pmatrix} 1 & 0  \\ c_1 & 1 \end{pmatrix} , 1_2) \, \bigg|
\ c_1 \in \mathbf C \bigg\} , \qquad
\widetilde U_2=\bigg\{ (1_2, \begin{pmatrix} 1 & 0  \\ c_2 & 1 \end{pmatrix} ) \, \bigg|
\ c_2 \in \mathbf C \bigg\} .
$$
Since $U$ is Zariski dense in $\widetilde U$, we have
$V^U=V^{\widetilde U}$. We also see easily that
$V_1^{\widetilde U_1}=\mathbf C\mathbf e_{l_1+1}$,
$V_2^{\widetilde U_2}=\mathbf C\mathbf e_{l_2+1}^{\prime}$,
$V^{\widetilde U}=V_1^{\widetilde U_1} \otimes V_2^{\widetilde U_2}$.
This completes the proof.

\proclaim Lemma 3.2. Let
$g=(\begin{pmatrix} 1 & 0  \\ c_1 & 1 \end{pmatrix} ,
\begin{pmatrix} 1 & 0  \\ c_2 & 1 \end{pmatrix}) \in {\rm SL}(2, \mathbf C)^2$.
We assume that $c_1 \neq 0$, $c_2 \neq 0$.
Then the dimension of the subspace of $V$ consisting of all vectors fixed by $g$ is
$l_2+1$. (Note that we have assumed $l_1 \geq l_2$.)
\par
{\bf Proof.} By the definition of the symmetric tensor representation $\rho _l$
of degree $l$, we have
\begin{equation}
\rho _l (\begin{pmatrix} 1 & 0 \\ c & 1 \end{pmatrix})\mathbf e_i
=\sum _{k=i}^{l+1} \binom {k-1}{i-1} c^{k-i}\mathbf e_k.
\tag{3.4}
\end{equation}
Hence for $\rho =\rho _{l_1} \otimes \rho _{l_2}$, we have
\begin{equation*}
\begin{aligned}
& \rho (\begin{pmatrix} 1 & 0  \\ c_1 & 1 \end{pmatrix} ,
\begin{pmatrix} 1 & 0  \\ c_2 & 1 \end{pmatrix})
(\mathbf e_i \otimes \mathbf e_j^{\prime}) \\
= & \sum _{k=i}^{l_1+1} \sum _{l=j}^{l_2+1}
\binom {k-1}{i-1} \binom {l-1}{j-1} c_1^{k-i}c_2^{l-j}
(\mathbf e_k \otimes \mathbf e_l^{\prime}).
\end{aligned}
\tag{3.5}
\end{equation*}
Put
$$
N=\rho (\begin{pmatrix} 1 & 0  \\ c_1 & 1 \end{pmatrix} ,
\begin{pmatrix} 1 & 0  \\ c_2 & 1 \end{pmatrix}), \qquad
n_{kl, ij}=\binom {k-1}{i-1} \binom {l-1}{j-1} c_1^{k-i}c_2^{l-j}.
$$
Then (3.5) can be written as
$$
N(\mathbf e_i \otimes \mathbf e_j^{\prime})
= \sum _{k=i}^{l_1+1} \sum _{l=j}^{l_2+1}
n_{kl, ij}(\mathbf e_k \otimes \mathbf e_l^{\prime}).
$$
The vector
$\sum _{i=1}^{l_1+1} \sum _{j=1}^{l_2+1}
x_{ij}(\mathbf e_i \otimes \mathbf e_j^{\prime})$
is annihilated by $N-\rho (1)$ if and only if
\begin{equation}
\sum _{i=1}^k \sum _{j=1, (i, j) \neq (k,l)}^l n_{kl, ij}x_{ij}=0
\tag{3.6}
\end{equation}
holds for every $1 \leq k \leq l_1+1$, $1 \leq l \leq l_2+1$.
Note that $n_{kl, ij} \neq 0$ if $k \geq l$ and $l \geq j$
since $c_1 \neq 0$, $c_2 \neq 0$; $n_{kl,ij}=0$ if $k<l$ or $l<j$.

First let $l=1$ in (3.6). Then we have
$$
\sum _{i=1}^{k-1} n_{k1, i1}x_{i1}=0, \qquad 2 \leq k \leq l_1+1.
$$
We successively obtain
$$
x_{11}=x_{21}=\cdots =x_{l_1 1}=0.
$$
We are going to show that
\begin{equation}
x_{1l}=x_{2l}=\cdots =x_{l_1+1-l, l}=0, \qquad 1 \leq l \leq l_2+1
\tag{3.7}
\end{equation}
by induction on $l$. Assuming the hypothesis of induction, we obtain
$$
\sum _{i=1}^{k-1} n_{kl, il} x_{il}=0
$$
when $k \leq l_1+2-l$ from (3.6). Hence (3.7) follows.
Now assume $k \geq l_1+3-l$. We write (3.6) as
$$
\sum _{j=1}^{l-1} \sum _{i=1}^k n_{kl,ij} x_{ij}
+\sum _{i=1}^{k-1} n_{kl,il}x_{il}=0.
$$
Using (3.7), we see that this equation is equivalent to
\begin{equation}
\sum _{j=1}^{l-1} \sum _{i=l_1+2-j}^k n_{kl,ij} x_{ij}
+\sum _{i=l_1+2-l}^{k-1} n_{kl,il}x_{il}=0.
\tag{3.8}
\end{equation}
In (3.8), $x_{il}$, $1 \leq i \leq l_1$ are determined by
$x_{\alpha\beta}$, $\beta \leq l-1$ and $x_{l_1+1 l}$
is left undermined. Therefore when $x_{l_1+1 j}$, $1 \leq j \leq l_2+1$
are given, the solution $x_{ij}$ satisfying (3.8) exists and is unique.
This completes the proof.

\proclaim Lemma 3.3. Let
$u=\begin{pmatrix} 1 & 0 \\ c & 1 \end{pmatrix}$, $0 \neq c \in F$.
Then we have
$$ 
\mathbf e_i\otimes \mathbf e_j^{\prime} \in {\rm Im}(u-1), \qquad
{\rm for} \ 1 \leq j \leq l_2+1 \ {\rm if} \ i \geq l_2+3-j.
$$
Here ${\rm Im}(u-1)$ denotes the image of the linear mapping
$V \ni \mathbf v \mapsto (\rho (u)-\rho (1))\mathbf v \in V$.
\par
{\bf Proof.} We have
$$
(u-1)(\mathbf e_i \otimes \mathbf e_{l_2+1}^{\prime})
=(ic\mathbf e_{i+1}+\binom {i+1}{i-1}c^2\mathbf e_{i+2}+\cdots )
\otimes \mathbf e_{l_2+1}^{\prime} .
$$
By descending induction on $i$, we see that
$\mathbf e_i \otimes \mathbf e_{l_2+1}^{\prime} \in {\rm Im}(u-1)$, $i \geq 2$.
This proves our assertion when $j=l_2+1$. When $j=l_2$, we have
$$
(u-1)(\mathbf e_i\otimes\mathbf e_{l_2}^{\prime})
=(ic\mathbf e_{i+1}+\cdots )\otimes\mathbf e_{l_2}^{\prime}
+(\mathbf e_i+ic\mathbf e_{i+1}+\cdots )
\otimes l_2c^{\prime}\mathbf e_{l_2+1}^{\prime}.
$$
The second term belongs to ${\rm Im}(u-1)$ if $i \geq 2$,
and by descending induction on $i$, we can show that
$\mathbf e_i \otimes \mathbf e_{l_2}^{\prime} \in {\rm Im}(u-1)$, $i \geq 3$.
Proceeding similarly by induction on $j$, we see that the assertion holds.

\proclaim Lemma 3.4. We have
$$
{\rm Im}(u_1-1)+{\rm Im} (u_2-1)
=(\oplus _{j=2}^{l_2+1} \mathbf C(\mathbf e_1\otimes \mathbf e_j^{\prime}))
\oplus (\oplus _{i=2}^{l_1+1} \oplus _{j=1}^{l_2+1}
\mathbf C(\mathbf e_i \otimes \mathbf e_j^{\prime})).
$$
In particular, $\dim ({\rm Im}(u_1-1)+{\rm Im} (u_2-1))=\dim V-1$.
\par
{\bf Proof.} Put $W={\rm Im}(u_1-1)+{\rm Im} (u_2-1)$.
Since it is clear that
$\mathbf e_1 \otimes \mathbf e_1^{\prime} \notin {\rm Im}(u_1-1)$,
$\notin {\rm Im}(u_2-1)$, it suffices to show that
$\mathbf e_i \otimes \mathbf e_j^{\prime} \in W$ when $(i, j) \neq (1, 1)$.
We have
$$
(u_2-1)(\mathbf e_1\otimes\mathbf e_{l_2}^{\prime})
=(\omega \mathbf e_2+\cdots )\otimes \mathbf e_{l_2}^{\prime}
+(\mathbf e_1+\omega\mathbf e_2+\cdots )
\otimes l_2\omega ^{\prime}\mathbf e_{l_2+1}^{\prime} .
$$
By Lemma 3.3, we have
$$
\omega (\mathbf e_2\otimes \mathbf e_{l_2}^{\prime})
+l_2\omega ^{\prime}(\mathbf e_1\otimes\mathbf e_{l_2+1}^{\prime})
\in {\rm Im}(u_2-1).
$$
Similarly we have
$$
(\mathbf e_2\otimes \mathbf e_{l_2}^{\prime})
+l_2(\mathbf e_1\otimes\mathbf e_{l_2+1}^{\prime})
\in {\rm Im}(u_1-1).
$$
Since
$\det \begin{pmatrix} 1 & l_2 \\ \omega & l_2\omega ^{\prime} \end{pmatrix} \neq 0$,
we obtain
$$
\mathbf e_1\otimes\mathbf e_{l_2+1}^{\prime} \in W, \qquad
\mathbf e_2\otimes \mathbf e_{l_2}^{\prime} \in W.
$$

We are going to show that
\begin{equation}
\mathbf e_i\otimes\mathbf e_j^{\prime} \in W, \ 1 \leq i \leq l_1+1, j \geq 2, \qquad
\mathbf e_i\otimes \mathbf e_{j-1}^{\prime} \in W, \ 2 \leq i \leq l_1+1, j \geq 2
\tag{$\ast$}
\end{equation}
by descending induction on $j$.
By Lemma 3.3 and by what we have shown, ($\ast$) holds when $j=l_2+1$.
We assume that ($\ast$) holds for $j+1$. We have
$$
(u_2-1)(\mathbf e_i\otimes\mathbf e_{j-1}^{\prime})
=(i\omega\mathbf e_{i+1}+\cdots )\otimes \mathbf e_{j-1}^{\prime}
+(\mathbf e_i+i\omega\mathbf e_{i+1}+\cdots )
\otimes ((j-1)\omega ^{\prime}\mathbf e_j^{\prime}+\cdots ).
$$
Suppose that $i \geq 2$. Then the second term on the right-hand side belongs to
$W$. By induction on $i$, we obtain
$$
\mathbf e_i\otimes\mathbf e_{j-1}^{\prime} \in W, \qquad i \geq 3.
$$
If $i=1$, we obtain
$$
i\omega\mathbf e_2\otimes\mathbf e_{j-1}^{\prime}
+(j-1)\omega ^{\prime}\mathbf e_1\otimes\mathbf e_j^{\prime} \in W.
$$
Considering $u_1-1$, we obtain
$$
i\mathbf e_2\otimes\mathbf e_{j-1}^{\prime}
+(j-1)\mathbf e_1\otimes\mathbf e_j^{\prime} \in W.
$$
Hence ($\ast$) holds for $j$. This completes the proof.

\medskip

By Lemma 3.1, we have
\begin{equation}
\dim \mathcal B=\dim V-1.
\tag{3.9}
\end{equation}
Consider the surjective linear mapping
$$
\mathcal Z \ni (U_1, U_2) \mapsto (u_2-1)U_1
\in {\rm Im}(u_1-1) \cap {\rm Im}(u_2-1).
$$
The kernel of this mapping consists of $(U_1, U_2)$ such that
$U_1 \in {\rm Ker}(u_2-1)$, $U_2 \in {\rm Ker}(u_1-1)$.
Hence by Lemma 3.2, we have
\begin{equation}
\dim \mathcal Z=\dim ({\rm Im}(u_1-1) \cap {\rm Im}(u_2-1))+2l_2+2.
\tag{3.10}
\end{equation}
By (3.9) and (3.10), we obtain
\begin{equation}
\dim H^1(U,V)=\dim ({\rm Im}(u_1-1) \cap {\rm Im}(u_2-1))+2l_2+3-\dim V.
\tag{3.11}
\end{equation}

\proclaim Lemma 3.5. We have $\dim H^1(U,V)=2$.
\par
{\bf Proof.} We have
\begin{equation*}
\begin{aligned}
& \dim ({\rm Im}(u_1-1) \cap {\rm Im}(u_2-1)) \\
= \, & \dim ({\rm Im}(u_1-1))+\dim ({\rm Im}(u_2-1))
-\dim ({\rm Im}(u_1-1)+{\rm Im} (u_2-1)).
\end{aligned}
\end{equation*}
By Lemma 3.2, we have $\dim ({\rm Im}(u_i-1))=\dim V-(l_2+1)$,
$i=1$, $2$. Then by Lemma 3.4, we get
$$
\dim ({\rm Im}(u_1-1) \cap {\rm Im}(u_2-1))
=\dim V-2l_2-1.
$$
The assertion follows from (3.11).

\medskip

{\bf 3.2.} In this subsection, we will prove the following theorems.

\proclaim Theorem 3.6. The eigenvalues of the action of $t$ on
$H^1(U, V)$ are $\epsilon ^{l_1+2}(\epsilon ^{\prime})^{-l_2}$
and $\epsilon ^{-l_1-2}(\epsilon ^{\prime})^{l_2}$.
In particular, $H^1(U, V)^{P/U}=0$.
\par

\proclaim Theorem 3.7.
$$
\dim H^1(P,V)=
\begin{cases}
0 \qquad {\rm if} \quad l_1\neq l_2 \ \ {\rm or} \ \ N(\epsilon )^{l_1}=-1, \\
1 \qquad {\rm if} \quad l_1=l_2 \ \ {\rm and} \ \ N(\epsilon )^{l_1}=1.
\end{cases}
$$
Here $N(\epsilon )$ denotes the norm of $\epsilon$.
\par

Taking $G=P$, $N=U$, $M=V$ in (1.9), we obtain the exact sequence
\begin{equation*}
\begin{CD}
0 @> >> & H^1(P/U, V^U) @> >> H^1(P,V)
@> >> H^1(U,V)^{P/U} @> >> 0,
\end{CD}
\end{equation*}
since $P/U \cong \mathbf Z$.
Therefore Theorem 3.7 follows immediately from Theorem 3.6,
since $\dim H^1(P/U, V^U)$ is easily seen to be equal to $0$ (resp. $1$) if
$l_1\neq l_2$ or $N(\epsilon )^{l_1}=-1$
(resp. if $l_1=l_2$ and $N(\epsilon )^{l_1}=1$), in view of Lemma 3.1.

\medskip

{\bf Proof of Theorem 3.6.} First we recall the following fact on the action of
$t$ on $H^q(U, V)$ (cf. (1.3)). Let $f \in Z^q(U, V)$ and let $\bar f \in H^q(U,V)$ be
the cohomology class represented by $f$. Put
$$
g(n_1, n_2, \ldots n_q)=t^{-1}f(tn_1t^{-1}, tn_2t^{-1}, \ldots , tn_qt^{-1}), \qquad
n_i \in U, \ 1 \leq i \leq q.
$$
Then $g \in Z^q(U,V)$ and $\bar f \mapsto \bar g$ is the action of $t$.

As in the proof of Lemma 3.4, let
$$
W={\rm Im}(u_1-1)+{\rm Im} (u_2-1)
=(\oplus _{j=2}^{l_2+1} \mathbf C(\mathbf e_1\otimes \mathbf e_j^{\prime}))
\oplus (\oplus _{i=2}^{l_1+1} \oplus _{j=1}^{l_2+1}
\mathbf C(\mathbf e_i \otimes \mathbf e_j^{\prime})).
$$
We have
$$
V=\mathbf C(\mathbf e_1\otimes\mathbf e_1^{\prime}) \oplus W.
$$
We may assume that $l_1>0$ since our assertion is clearly true if
$l_1=l_2=0$.

Put $\mathbf t_1=\mathbf e_1\otimes\mathbf e_{l_2+1}^{\prime}$.
Let us show that for
$$
\mathbf t_2=\omega (\mathbf e_1\otimes\mathbf e_{l_2+1}^{\prime})
+\sum _{i=2}^{l_1+1} x_i(\mathbf e_i\otimes\mathbf e_{l_2+1}^{\prime})
$$
with suitably chosen $x_i \in \mathbf C$, we have
\begin{equation}
(u_2-1)\mathbf t_1=(u_1-1)\mathbf t_2.
\tag{3.12}
\end{equation}
To this end, for $i \geq 1$, put
$$
W_i=\oplus _{k=i}^{l_1+1} \mathbf C
(\mathbf e_k\otimes\mathbf e_{l_2+1}^{\prime}).
$$
We have
\begin{equation*}
\begin{aligned}
& (u_2-1)(\mathbf e_1\otimes\mathbf e_{l_2+1}^{\prime})
=(\omega\mathbf e_2+\omega ^2\mathbf e_3+\cdots )
\otimes \mathbf e_{l_2+1}^{\prime}, \\
& (u_1-1)(\mathbf e_i\otimes\mathbf e_{l_2+1}^{\prime})
=(i\mathbf e_{i+1}+\binom {i+1}{i-1}\mathbf e_{i+2}+\cdots )
\otimes \mathbf e_{l_2+1}^{\prime}.
\end{aligned}
\end{equation*}
We see that
$$
(u_2-1)\mathbf t_1 \equiv (u_1-1)\mathbf t_2 \mod W_3.
$$
For $x_2=(\omega ^2-\omega )/2$, we have
$$
(u_2-1)\mathbf t_1 \equiv (u_1-1)\mathbf t_2 \mod W_4.
$$
In this way, we can determine $x_i$ successively so that (3.12) holds.
Let $f_1 \in Z^1(U, V)$ be the $1$-cocycle which corresponds to
the point $(\mathbf t_1, \mathbf t_2) \in \mathcal Z$.

Put $\mathbf t_3=\mathbf e_{l_1+1}\otimes\mathbf e_1^{\prime}$.
Similarly to the above, we can show that for
$$
\mathbf t_4=\omega ^{\prime}(\mathbf e_{l_1+1}\otimes\mathbf e_1^{\prime})
+\sum _{j=2}^{l_2+1} y_j(\mathbf e_{l_1+1}\otimes\mathbf e_j^{\prime}),
$$
the relation
\begin{equation}
(u_2-1)\mathbf t_3=(u_1-1)\mathbf t_4
\tag{3.13}
\end{equation}
holds when $y_j$ are suitably chosen.
Let $f_2 \in Z^1(U, V)$ be the $1$-cocycle which corresponds to
the point $(\mathbf t_3, \mathbf t_4) \in \mathcal Z$.

Let $\bar f_i$ be the class of $f_i$ in $H^1(U, V)$, $i=1$, $2$.
Let us show that $\{ \bar f_1, \bar f_2 \}$ gives a basis of $H^1(U, V)$.
To this end, assume that $\alpha f_1+\beta f_2 \in B^1(U, V)$
for $\alpha$, $\beta \in \mathbf C$. Then there exists $\mathbf b \in V$
such that
\begin{equation}
\alpha\mathbf t_1+\beta\mathbf t_3=(u_1-1)\mathbf b,
\tag{i}
\end{equation}
\begin{equation}
\alpha\mathbf t_2+\beta\mathbf t_4=(u_2-1)\mathbf b
\tag{ii}
\end{equation}
hold. Put
$$
\mathbf b=\sum _{i=1}^{l_1+1}\sum _{j=1}^{l_2+1}
x_{ij}(\mathbf e_i\otimes\mathbf e_j^{\prime}).
$$
On the left-hand side of (i), the coefficient of the tensor
$\mathbf e_1\otimes\mathbf e_{l_2+1}^{\prime}$ is $\alpha$
and the coefficients of $\mathbf e_1\otimes\mathbf e_j^{\prime}$ are
$0$ for $1 \leq j \leq l_2$. We have
$$
(u_1-1)(\mathbf e_1\otimes\mathbf e_j^{\prime})
=j(\mathbf e_1\otimes\mathbf e_{j+1}^{\prime})
+\sum _{l=j+2}^{l_2+1} z_l(\mathbf e_1\otimes\mathbf e_l^{\prime})
+A,
$$
where $z_l \in \mathbf Z$ and $A$ is a term which does not contain
$\mathbf e_1\otimes\mathbf e_l^{\prime}$. Therefore we have
$x_{11}=\cdots =x_{1 l_2-1}=0$. By comparing the coefficients of the tensor
$\mathbf e_1\otimes\mathbf e_{l_2+1}^{\prime}$ on the both sides of (i),
we obtain
$$
\alpha =l_2x_{1 l_2}.
$$
By comparing the coefficients of the tensor
$\mathbf e_1\otimes\mathbf e_{l_2+1}^{\prime}$ on the both sides of (ii),
we get
$$
\alpha\omega=l_2\omega ^{\prime}x_{1 l_2}.
$$
Hence we obtain $x_{1 l_2}=0$, $\alpha =0$.
Similarly by comparing the coefficients of the tensor
$\mathbf e_{l_1+1}\otimes\mathbf e_1^{\prime}$
for the both sides of (i) and (ii), we obtain $\beta =0$.

Let $f_1^{\prime}$ be the image of $f_1$ under the action of $t$
and let $(U_1^{\prime}, U_2^{\prime}) \in \mathcal Z$ be the point
corresponding to $f_1^{\prime}$. Then we have
\begin{equation*}
U_1^{\prime}=f_1^{\prime}(u_1)
=t^{-1}f_1(tu_1t^{-1})=t^{-1}f_1(u_1^Au_2^B)
=t^{-1}[u_1^Af_1(u_2^B)+f_1(u_1^A)],
\end{equation*}
\begin{equation*}
U_2^{\prime}=f_1^{\prime}(u_2)
=t^{-1}f_1(tu_2t^{-1})=t^{-1}f_1(u_1^Cu_2^D)
=t^{-1}[u_1^Cf_1(u_2^D)+f_1(u_1^C)].
\end{equation*}
For $i=1$, $2$, we have $f_1(u_i)=\mathbf t_i$ and
\begin{equation}
f_1(u_i^n)=(1+u_1+\cdots +u_1^{n-1})\mathbf t_i \qquad
{\rm if} \quad n>0.
\tag{3.14}
\end{equation}
\begin{equation}
f_1(u_i^{-n})=-(u_1^{-1}+\cdots +u_1^{-n})\mathbf t_i \qquad
{\rm if} \quad n>0.
\tag{3.15}
\end{equation}
From these formulas, we see easily that the coefficient of
$\mathbf e_1\otimes\mathbf e_{l_2+1}^{\prime}$ in $tU_1^{\prime}$ is
$A+B\omega$. Hence the coefficient of
$\mathbf e_1\otimes\mathbf e_{l_2+1}^{\prime}$ in $U_1^{\prime}$ is
$\epsilon ^{l_1}(\epsilon ^{\prime})^{-l_2}(A+B\omega )
=\epsilon ^{l_1+2}(\epsilon ^{\prime})^{-l_2}$. Similarly we see that
the coefficient of
$\mathbf e_1\otimes\mathbf e_{l_2+1}^{\prime}$ in $U_2^{\prime}$ is
$\omega\epsilon ^{l_1+2}(\epsilon ^{\prime})^{-l_2}$. 

Now let
$$
f_1^{\prime}\equiv \gamma f_1+\delta f_2 \mod B^1(U, V)
$$
with $\gamma$, $\delta \in \mathbf C$. Then there exists $\mathbf c \in V$
such that
\begin{equation}
\gamma\mathbf t_1+\delta\mathbf t_3-U_1^{\prime}
=(u_1-1)\mathbf c,
\tag{iii}
\end{equation}
\begin{equation}
\gamma\mathbf t_2+\delta\mathbf t_4-U_2^{\prime}
=(u_2-1)\mathbf c.
\tag{iv}
\end{equation}
Put $\mathbf c=\sum _{i=1}^{l_1+1}\sum _{j=1}^{l_2+1} y_{ij}
(\mathbf e_i\otimes\mathbf e_j^{\prime})$. Comparing the coefficients of
$\mathbf e_1\otimes\mathbf e_{l_2+1}^{\prime}$ on the both sides of (iii),
we obtain
$$
\gamma -\epsilon ^{l_1+2}(\epsilon ^{\prime})^{-l_2}=l_2y_{1l_2}.
$$
Comparing the coefficients of
$\mathbf e_1\otimes\mathbf e_{l_2+1}^{\prime}$ on the both sides of (iv),
we obtain
$$
(\gamma -\epsilon ^{l_1+2}(\epsilon ^{\prime})^{-l_2})\omega 
=l_2y_{1l_2}\omega ^{\prime}.
$$
From these two formulas, we obtain
$y_{1l_2}=0$, $\gamma =\epsilon ^{l_1+2}(\epsilon ^{\prime})^{-l_2}$.
Similarly comparing the coefficients of
$\mathbf e_{l_1+1}\otimes\mathbf e_1^{\prime}$ on the both sides of
(iii) and (iv), we obtain $\delta =0$. Thus we have shown
\begin{equation}
f_1^{\prime} \equiv \epsilon ^{l_1+2}(\epsilon ^{\prime})^{-l_2}f_1
\mod B^1(U, V).
\tag{3.16}
\end{equation}

Next let $f_2^{\prime}$ be the image of $f_2$ under the action of $t$
and let $(U_3^{\prime}, U_4^{\prime})$ be the point of $\mathcal Z$
corresponding to $f_2^{\prime}$. Here
$U_3^{\prime}=f_2^{\prime}(u_1)$, $U_4^{\prime}=f_2^{\prime}(u_2)$.
Then we have
\begin{equation*}
U_3^{\prime}=f_2^{\prime}(u_1)
=t^{-1}f_2(tu_1t^{-1})=t^{-1}f_2(u_1^Au_2^B)
=t^{-1}[u_1^Af_2(u_2^B)+f_2(u_1^A)],
\end{equation*}
\begin{equation*}
U_4^{\prime}=f_2^{\prime}(u_2)
=t^{-1}f_2(tu_2t^{-1})=t^{-1}f_2(u_1^Cu_2^D)
=t^{-1}[u_1^Cf_2(u_2^D)+f_2(u_1^C)].
\end{equation*}
The coefficient of $\mathbf e_{l_1+1}\otimes\mathbf e_1^{\prime}$ in
$tU_3^{\prime}$ is $A+B\omega ^{\prime}=\epsilon ^{-2}$.
The coefficient of $\mathbf e_{l_1+1}\otimes\mathbf e_1^{\prime}$ in
$tU_4^{\prime}$ is $C+D\omega ^{\prime}=\epsilon ^{-2}\omega ^{\prime}$.
By a similar argument to the above, we obtain
\begin{equation}
f_2^{\prime} \equiv \epsilon ^{-l_1-2}(\epsilon ^{\prime})^{l_2}f_2
\mod B^1(U, V).
\tag{3.17}
\end{equation}
This completes the proof of Theorem 3.6.

\medskip

{\bf 3.3.} In this subsection, we will prove the following theorems.

\proclaim Theorem 3.8. We have $\dim H^2(U,V)=1$ and $t$ acts on it
as the multiplication by $\epsilon ^{l_1}(\epsilon ^{\prime})^{l_2}$.
\par

\proclaim Theorem 3.9.
We have $H^2(P, V)=0$ except for the case
$l_1=l_2$ and $N(\epsilon )^{l_1}=1$.
If $l_1=l_2$ and $N(\epsilon )^{l_1}=1$, then we have $\dim H^2(P,V)=1$.
\par

First we will prove the part of Theorem 3.8 concerning the dimension.

\proclaim Lemma 3.10. We have $\dim H^2(U,V)=1$.
\par
{\bf Proof.} Let $U_1$ be the subgroup of $U$ generated by $u_1$.
We have the exact sequence 
\begin{equation}
\begin{CD}
0 @> >> & U_1 @> >> U @> >> U_2 @> >> 0
\end{CD}
\tag{3.18}
\end{equation}
and the associated spectral sequence (cf. (1.8))
\begin{equation}
E_2^{p,q}=H^p(U_2, H^q(U_1,V)) \Longrightarrow H^n(U,V).
\tag{3.19}
\end{equation}
Let $E^n=H^n(U,V)$ and $\{ F^i \}$ denote the filtration on $E^n$
induced by (3.19). We have
$F^p(E^n)/F^{p+1}(E^n)\cong E_{\infty}^{p,n-p}$.
Since $U_2\cong\mathbf Z$, we have $E_2^{2,q}=E_{\infty}^{2,q}=0$.
Since $F^3(E^2)=0$, we get $F^2(E^2)=0$. Since $U_1\cong\mathbf Z$,
we have $E_2^{p,2}=E_{\infty}^{p,2}=0$. Hence we get $E^2/F^1(E^2)=0$.
We have $F^1(E^2)/F^2(E^2) \cong E_{\infty}^{1,1}$.
Therefore it is sufficient to show that $\dim E_{\infty}^{1,1}=1$.

We consider
$$
E_2^{1,1}=H^1(U_2, H^1(U_1, V)).
$$
The map $Z^1(U_1,V) \ni f \mapsto f(u_1) \in V$ induces the isomorphism
\begin{equation}
H^1(U_1, V) \cong V/{\rm Im}(u_1-1).
\tag{3.20}
\end{equation}
The action of $u \in U_2$ on the right-hand side of (3.20) is given by
$$
V/{\rm Im}(u_1-1) \ni v \hskip -0.5em \mod {} {\rm Im}(u_1-1) \longrightarrow
u^{-1}v \hskip -0.5em \mod {} {\rm Im}(u_1-1) \in V/{\rm Im}(u_1-1).
$$
Since $\bar u_2=u_2 \mod U_1$ is a generator of $U_2$, we have
$$
H^1(U_2, H^1(U_1, V)) \cong (V/{\rm Im}(u_1-1))/{\rm Im}(\bar u_2-1)
\cong V/({\rm Im}(u_1-1)+{\rm Im}(u_2-1)).
$$
By Lemma 3.4, we obtain
$$
\dim H^1(U_2, H^1(U_1, V))=\dim E_2^{1,1}=1.
$$
Since $E_2^{3,0}=0$, we have $E_{\infty}^{1,1}=E_2^{1,1}$.
This completes the proof.

\medskip

{\bf Proof of Theorem 3.8.} We set $\tau =u_1$, $\eta =u_2$.
Let $\mathcal F$ be the free group on two free generators
$\widetilde\tau$ and $\widetilde\eta$ and let
$\pi :\mathcal F \longrightarrow U$ be the surjective homomorphism
such that
$$
\pi (\widetilde\tau )=\tau , \qquad \pi (\widetilde\eta )=\eta .
$$
Let $R$ be the kernel of $\pi$. For $a$, $b \in \mathcal F$,
let $[a, b]=aba^{-1}b^{-1}$ be the commutator of $a$ and $b$.
We see easily that
$$
R=\langle x[\widetilde\tau , \widetilde\eta ]x^{-1} \mid x \in \mathcal F \rangle ,
\qquad R=[\mathcal F, \mathcal F].
$$
We have the isomorphism (cf. (1.11))
\begin{equation}
H^2(U, V) \cong H^1(R, V)^U/{\rm Im}(H^1(\mathcal F, V)).
\tag{3.21}
\end{equation}
We have
\begin{equation}
H^1(R,V)^U=\{ \varphi \in {\rm Hom}(R, V) \mid
\varphi (grg^{-1})=g\varphi (r), \ \ g \in \mathcal F, \ r \in R \} .
\tag{3.22}
\end{equation}
Hence $\varphi \in H^1(R, V)^U$ is completely determined by
$\varphi ([\widetilde\tau , \widetilde\eta ])$.
For $b \in H^1(\mathcal F, V)$, we have
$$
b([\widetilde\tau , \widetilde\eta ])
=(1-\eta )b(\widetilde\tau )+(\tau -1)b(\widetilde\eta ).
$$
Let $W$ be the suspace of $V$ as in the proof of Lemma 3.4.
For $\varphi \in {\rm Im}(H^1(\mathcal F, V))$,
the formula above shows that
$\varphi ([\widetilde\tau , \widetilde\eta ])$ can take an arbitrary vector in $W$.
In particular, it follows that $\dim H^2(U, V) \leq 1$. Since
$\dim H^2(U, V)=1$ by Lemma 3.10, we see that there exists
$\varphi _1\in H^1(R, V)^U$ such that
$\varphi  _1([\widetilde\tau , \widetilde\eta ])=\mathbf e_1\otimes\mathbf e_1^{\prime}$.
This $\varphi _1$ corresponds to a generator of $H^2(U, V)$.

Let $f \in Z^2(U, V)$. For $g \in \mathcal F$, we put $\bar g=\pi (g)$.
There exists $a \in C^1(\mathcal F, V)$ such that (cf. (1.12))
\begin{equation}
f(\bar g_1, \bar g_2)=g_1a(g_2)+a(g_1)-a(g_1g_2), \qquad g_1, g_2 \in \mathcal F.
\tag{3.23}
\end{equation}
The corresponding element $\varphi \in H^1(R,V)^U$ to $f$ is obtained as
the restriction of $a$ to $R$.
Now let $\xi$ be an automorphism of $\mathcal F$. Since $\xi$ stabilizes
$R=[\mathcal F, \mathcal F]$, $\xi$ induces an automorphism of $U=\mathcal F/R$,
which we denote by $\bar\xi$. We have
$$
\bar\xi (\bar g)=\overline {\xi (g)}, \qquad g \in \mathcal F.
$$
From (3.23), we obtain
\begin{equation}
f(\bar\xi (\bar g_1), \bar\xi (\bar g_2))
=\xi (g_1)a(\xi (g_2))+a(\xi (g_1))-a(\xi (g_1)\xi (g_2)), \ \ g_1, g_2 \in \mathcal F.
\tag{3.24}
\end{equation}

\proclaim Lemma 3.11. For
$\gamma =\begin{pmatrix} a & b \\ c & d \end{pmatrix} \in {\rm SL}(2, \mathbf Z)$,
let $\xi (\gamma )$ be the automorphism of $U$ defined by
$\xi (\gamma )(\tau )=\tau ^a\eta ^c$, $\xi (\gamma )(\eta )=\tau ^b\eta ^d$.
Then there exists an automorphism $\widetilde{\xi (\gamma )}$ of $\mathcal F$ such that
$\xi (\gamma )=\overline{\widetilde{\xi (\gamma )}}$. Moreover $\widetilde{\xi (\gamma )}$
can be taken so that
\begin{equation}
\varphi (\widetilde{\xi (\gamma )}(g))
\equiv \varphi (g) \mod W
\tag{3.25}
\end{equation}
holds for every $\varphi \in H^1(R, V)^U$ and every $g \in [\mathcal F, \mathcal F]$.
\par
{\bf Proof.} For $\gamma _1$, $\gamma _2 \in {\rm SL}(2, \mathbf Z)$, we have
$\xi (\gamma _1\gamma _2)=\xi (\gamma _1)\xi (\gamma _2)$.
For two automorphisms $\xi _1$, $\xi _2$ of $\mathcal F$, we have
$\overline{\xi _1\xi _2}=\bar\xi _1\bar\xi _2$.
Therefore to show the first assertion, it is sufficient to verify it for generators
$\gamma _1=\begin{pmatrix} 1 & 1 \\ 0 & 1 \end{pmatrix}$,
$\gamma _2=\begin{pmatrix} 0 & 1 \\ -1 & 0 \end{pmatrix}$ of ${\rm SL}(2, \mathbf Z)$.
Clearly the formulas
$\widetilde\xi (\gamma _1)(\widetilde\tau )=\widetilde\tau$,
$\widetilde\xi (\gamma _1)(\widetilde\eta )=\widetilde\tau\widetilde\eta$,
$\widetilde\xi (\gamma _2)(\widetilde\tau )=\widetilde\eta ^{-1}$,
$\widetilde\xi (\gamma _2)(\widetilde\eta )=\widetilde\tau$
define automorphisms $\widetilde\xi (\gamma _1)$ and $\widetilde\xi (\gamma _2)$
of $\mathcal F$ satisfying the requirements.

To show the latter assertion, we first note that
\begin{equation}
u\mathbf v\equiv \mathbf v \mod W \quad
\text {for every $u \in U$ and every $\mathbf v \in V$}.
\tag{3.26}
\end{equation}
Let $\varphi \in H^1(R, V)^U$.
Since $\widetilde\xi (\gamma )$ can be taken from the subgroup of ${\rm Aut}(\mathcal F)$
generated by $\widetilde\xi (\gamma _1)$ and $\widetilde\xi (\gamma _2)$, it is sufficient
to show (3.25) for these generators. Moreover since
$\varphi (x [\widetilde\tau , \widetilde\eta ]x^{-1})=x\varphi ([\widetilde\tau , \widetilde\eta ])$
for $x \in \mathcal F$, it is enough to verify (3.25) for $g=[\widetilde\tau , \widetilde\eta ]$
in view of (3.26). For $\widetilde\xi (\gamma _1)$, we have
$$
\varphi (\widetilde\xi (\gamma _1)([\widetilde\tau , \widetilde\eta ]))
=\varphi (\widetilde\tau [\widetilde\tau , \widetilde\eta ]\widetilde\tau ^{-1})
=\tau\varphi ([\widetilde\tau , \widetilde\eta ]) \equiv
\varphi ([\widetilde\tau , \widetilde\eta ]) \mod W
$$
by (3.26). For $\widetilde\xi (\gamma _2)$, we can check (3.25) similarly since
$\widetilde\xi (\gamma _2)([\widetilde\tau , \widetilde\eta ])
=\widetilde\eta ^{-1}[\widetilde\tau , \widetilde\eta ]\widetilde\eta$.
This completes the proof of Lemma 3.11.

Applying Lemma 3.11 to
$\gamma =\begin{pmatrix} A & B \\ C & D \end{pmatrix}$,
we see that there exists an automorphism $\xi _t$ of $\mathcal F$ such that (cf. (3.1))
$$
\bar\xi _t(u)=tut^{-1}, \qquad u \in U.
$$
Under the action of $t$, $f$ is transformed to the $2$-cocycle $f^{\prime} \in Z^2(U,V)$
where
$$
f^{\prime}(h_1, h_2)=t^{-1}f(th_1t^{-1}, th_2t^{-1}), \qquad h_1, h_2 \in U.
$$
By (3.24), we obtain
\begin{equation}
\begin{aligned}
& t^{-1}f(t\bar g_1t^{-1}, t\bar g_2t^{-1}) \\
= \, & g_1t^{-1}a(\xi _t(g_2))+t^{-1}a(\xi _t(g_1))-t^{-1}a(\xi _t(g_1)\xi _t(g_2)),
\qquad g_1, g_2 \in \mathcal F.
\end{aligned}
\tag{3.27}
\end{equation}
This formula shows that $1$-cochain $a^{\prime} \in C^1(\mathcal F, V)$ which splits
$f^{\prime}$ is given by
$$
a^{\prime}(g)=t^{-1}a(\xi _t(g)), \qquad g \in \mathcal F.
$$

Now suppose that $f$ (resp. $f^{\prime}$) $\in Z^2(U, V)$ corresponds to
$\varphi$ (resp. $\varphi ^{\prime}$) $\in H^1(R, V)^U$. We have
\begin{equation}
\varphi ^{\prime}([\widetilde\tau , \widetilde\eta ])
=t^{-1}\varphi (\xi _t([\widetilde\tau , \widetilde\eta ])).
\tag{3.28}
\end{equation}
We may assume that $\varphi =\varphi _1$, i.e., $\varphi ([\widetilde\tau , \widetilde\eta ])
=\mathbf e_1\otimes\mathbf e_1^{\prime}$. Then by (3.25), we obtain
$$
\varphi ^{\prime}([\widetilde\tau , \widetilde\eta ])
\equiv t^{-1}\varphi ([\widetilde\tau , \widetilde\eta ])
\equiv \epsilon ^{l_1}(\epsilon ^{\prime})^{l_2}
\varphi ([\widetilde\tau , \widetilde\eta ]) \mod W.
$$
This completes the proof of Theorem 3.8.

\medskip

{\bf Proof of Theorem 3.9.} Set $T=P/U$. Then $T$ is generated by $t \mod U$.
We consider the spectral sequence
\begin{equation}
E_2^{p,q}=H^p(T, H^q(U,V)) \Longrightarrow H^n(P,V).
\tag{3.29}
\end{equation}
Let $E^n=H^n(P,V)$ and $\{ F^i \}$ denote the filtration induced by (3.29).
Since $T \cong \mathbf Z$, we have $E_2^{p,q}=0$ for $p \geq 2$, $q \geq 0$.
Hence $F^2(E^2)/F^3(E^2)\cong E_{\infty}^{2,0}=0$.
Since $F^3(E^2)=0$, we obtain $F^2(E^2)=0$. By Theorem 3.6, we have
$E_2^{1,1}=H^1(T, H^1(U, V))=0$. Hence we have
$F^1(E^2)/F^2(E^2)\cong E_{\infty}^{1,1}=0$. Therefore we obtain
\begin{equation}
\dim H^2(P, V)=\dim E^2/F^1(E^2)=\dim E_{\infty}^{0,2}.
\tag{3.30}
\end{equation}

Now assume $l_1\neq l_2$ or $N(\epsilon )^{l_1} \neq 1$.
By Theorem 3.8, we have $H^2(U, V)^T=0$. Hence we get
$E_2^{0,2}=E_{\infty}^{0,2}=0$. Next assume that
$l_1=l_2$ and $N(\epsilon )^{l_1}=1$. By Theorem 3.8, we have
$\dim E_2^{0,2}=\dim H^2(U, V)^T=1$. We clearly have
$E_2^{0,2}\cong E_{\infty}^{0,2}$. This completes the proof.

\bigskip

\centerline{\S4. On the parabolic condition}

\medskip

In this section (in particular subsection 4.1), we will show that
it is possible to deduce information on critical values of $L$-functions,
once we know a corresponding $2$-cocycle which satisfies the parabolic condition.

From this section until the end of the paper, we define subgroups of $\Gamma$ by
\begin{equation*}
\begin{aligned}
& P=\left\{ \begin{pmatrix} u & v \\ 0 & u^{-1} \end{pmatrix} \bigg|
\ u \in E_F, \ v \in \mathcal O_F \right\}/\{ \pm 1_2 \} , \\
& U=\left\{ \begin{pmatrix} \pm 1 & v \\ 0 & \pm 1 \end{pmatrix} \bigg|
\ v \in \mathcal O_F \right\} /\{ \pm 1_2 \}
\end{aligned}
\end{equation*}
restoring the notation to that of \S2. We see that Theorems 3.7, 3.9
and the fact $H^1(P,V)^{P/U}=0$ stated in Theorem 3.6 are valid,
considering the isomorphism $P \ni p \mapsto {}^tp^{-1} \in {}^tP$ and noting that
$g \mapsto \rho (g)$ and $g \mapsto \rho ({}^t g^{-1})$ are equivalent as
representations of ${\rm SL}(2, \mathbf C)^2$.

\medskip

{\bf 4.1.} Let $V_1$ (resp. $V_2$) be the representation space of $\rho _{l_1}$
(resp. $\rho _{l_2}$). We take a basis
$\{ \mathbf e_1, \mathbf e_2, \ldots , \mathbf e_{l_1+1} \}$ of $V_1$
so that $\rho _{l_1}(\begin{pmatrix} a & 0 \\ 0 & 1 \end{pmatrix})\mathbf e_i
=a^{l_1+1-i}\mathbf e_i$.
Similarly we take a basis
$\{ \mathbf e_1^{\prime}, \mathbf e_2^{\prime}, \ldots , \mathbf e_{l_2+1}^{\prime} \}$
of $V_2$ so that
$\rho _{l_2}(\begin{pmatrix} a & 0 \\ 0 & 1 \end{pmatrix})\mathbf e_i^{\prime}
=a^{l_2+1-i}\mathbf e_i^{\prime}$.
We assume that $l_1 \geq l_2$, $l_1 \equiv l_2 \mod 2$.
We put $k_1=l_1+2$, $k_2=l_2+2$, $k=(k_1, k_2)$.
Let $\Omega \in S_k(\Gamma )$.
We assume that $l_1$ is even if $N(\epsilon )=-1$. (This assumption is
(A) in \S1.)

We recall the formulas:
\begin{equation}
f (\gamma _1, \gamma _2)
=\int _{\gamma _1\gamma _2w_1}^{\gamma _1w_1}
\int _{w_2}^{\gamma _1^{\prime}w_2} \mathfrak d(\Omega),
\qquad w_1=i\epsilon ^{-1}, w_2=i\infty ,
\tag{4.1}
\end{equation}

\begin{equation}
f(\sigma , \mu )
=-\int _{i\epsilon ^{-1}}^{i\epsilon} \int _0^{i\infty} \mathfrak d(\Omega ).
\tag{4.2}
\end{equation}
The formula (2.30) shows that the coefficients of
$\mathbf e_i \otimes \mathbf e_{i-(l_1-l_2)/2}^{\prime}$ in $f(\sigma , \mu )$,
$(l_1-l_2)/2+1 \leq i \leq (l_1+l_2)/2+1$ are related to the critical values
of $L(s, \Omega)$.

The parabolic condition on the cocycle $f$ is
\begin{equation}
f(p\gamma _1, \gamma _2)=pf(\gamma _1, \gamma _2)
\qquad \text {for every} \quad p \in P, \ \gamma _1, \gamma _2 \in \Gamma .
\tag{4.3}
\end{equation}

Now suppose that we add the coboundary of $b \in C^1(\Gamma , V)$
$$
b(\gamma _1\gamma _2)-\gamma _1b(\gamma _2)-b(\gamma _1)
$$
to $f$. We assume that the resulting $2$-cocycle is normalized and still satisfies
the parabolic condition (4.3). Then we see easily that $b$ must satisfy the condition
\begin{equation}
b(p\gamma )=pb(\gamma )+b(p), \qquad p \in P, \ \gamma \in \Gamma .
\tag{4.4}
\end{equation}
Put $A=f(\sigma , \mu )$. After adding the cobounday of $b$, $A$ changes to
$A+b(\sigma\mu )-\sigma b(\mu )-b(\sigma )$. By (4.4), we have
$$
b(\sigma\mu )=b(\mu ^{-1}\sigma )=\mu ^{-1}b(\sigma )+b(\mu ^{-1}), \qquad
b(\mu ^{-1})=-\mu ^{-1}b(\mu ).
$$
Therefore $A$ changes to
$$
A+(\mu ^{-1}-1)b(\sigma )-(\sigma +\mu ^{-1})b(\mu ).
$$
By (4.4), we have $b \vert P \in Z^1(P, V)$.
Suppose that $l_1\neq l_2$. By Theorem 3.7, we have
$$
b(\mu )=(\mu -1)\mathbf b, \qquad \mathbf b \in V.
$$
Since $(\sigma +\mu ^{-1})(\mu -1)=(\mu ^{-1}-1)(\sigma -1)$,
we see that $A$ changes to
$$
A+(\mu ^{-1}-1)[b(\sigma )+(1-\sigma )\mathbf b].
$$
This formula shows that the components of $A$ related to the critical values
do not change by adding a coboundary, since
$\mu ^{-1}(\mathbf e_i\otimes\mathbf e_{i-(l_1-l_2)/2})
=\mathbf e_i\otimes\mathbf e_{i-(l_1-l_2)/2})$.
Next suppose that $l_1=l_2$.
By Theorem 3.7 and by the exact sequence below it, we have
$$
b(\mu )=(\mu -1)\mathbf b+\mathbf b_0, \qquad \mathbf b \in V, \quad \mathbf b_0 \in V^U.
$$
Hence $A$ changes to
$$
A+(\mu ^{-1}-1)[b(\sigma )+(1-\sigma )\mathbf b] -(\sigma +\mu ^{-1})\mathbf b_0.
$$
Since $\mathbf b_0 \in V^U$, this formula shows that the components of $A$
related to the critical values do not change except for two critical values
$L(1, \Omega )$ and $L(l_1+1, \Omega )$ at the edges.

\medskip

{\bf 4.2.} Let $\bar Z^2(\Gamma , V)$ be the subgroup of $Z^2(\Gamma , V)$ consisting of
normalized $2$-cocycles. Put
$$
\bar B^2(\Gamma , V)= \{ f=db \mid b \in C^1(\Gamma , V), \ b(1)=0 \} .
$$
Then we have
$$
\bar Z^2(\Gamma , V)\cap B^2(\Gamma , V)=\bar B^2(\Gamma , V)
$$
and therefore
$$
\bar Z^2(\Gamma , V)/\bar B^2(\Gamma , V) \subset  Z^2(\Gamma , V)/B^2(\Gamma , V).
$$
Since every $2$-cocycle can be normalized by adding a coboundary, we have
$$
H^2(\Gamma , V)=\bar Z^2(\Gamma , V)/\bar B^2(\Gamma , V).
$$
Put
\begin{equation}
Z_{\rm P}^2(\Gamma , V)= \{ f \in \bar Z^2(\Gamma , V) \mid
\text {$f$ satisfies the parabolic condition (4.3)} \} ,
\tag{4.5}
\end{equation}
\begin{equation*}
\begin{aligned}
B_{\rm P}^2(\Gamma , V)= & \{ f \in \bar B^2(\Gamma , V) \mid
f=db, b \in C^1(\Gamma , V), \\
& b(p\gamma )=pb(\gamma )+b(p), \quad
p \in P, \gamma \in \Gamma \} .
\end{aligned}
\tag{4.6}
\end{equation*}
An element of $Z_P^2(\Gamma , V)$ is called a
{\it normalized parabolic $2$-cocycle}. The next lemma can easily be verified.

\proclaim Lemma 4.1. We have
$$
Z_{\rm P}^2(\Gamma , V)\cap \bar B^2(\Gamma , V)=B_{\rm P}^2(\Gamma ,V).
$$
\par

By Lemma 4.1, we have
$$
Z_{\rm P}^2(\Gamma , V)/B_{\rm P}^2(\Gamma , V) \subset 
\bar Z^2(\Gamma , V)/\bar B^2(\Gamma , V)=H^2(\Gamma , V).
$$
We define the parabolic part $H_{\rm P}^2(\Gamma , V)$ of $H^2(\Gamma , V)$ by
\begin{equation}
H_{\rm P}^2(\Gamma , V)=Z_{\rm P}^2(\Gamma , V)/B_{\rm P}^2(\Gamma , V).
\tag{4.7}
\end{equation}

{\bf 4.3.} As another application of Theorem 3.7, we are going to show
the nonvanishing of the cohomology class attached to a Hecke eigenform.

\proclaim Lemma 4.2. Assume $l_1$ is even if $N(\epsilon )=-1$.
Let $f \in Z_P^2(\Gamma , V)$ be a normalized
parabolic $2$-cocycle. For $(l_1-l_2)/2+1 \leq i \leq (l_1+l_2)/2+1$, let $c_i$ be
the coefficient of
$\mathbf e_i \otimes \mathbf e_{i-(l_1-l_2)/2}^{\prime}$ in $f(\sigma , \mu )$.
Assume that $c_i \neq 0$ for some $i$ if $l_1 \neq l_2$ and that $c_i \neq 0$
for some $i \neq 1, l_1+1$ if $l_1=l_2$. Then the cohomology class of
$f$ is non-trivial.
\par
{\bf Proof.} Suppose that the cohomology class of $f$ is trivial.
Then there exists $b \in C^1(\Gamma , V)$ such that
$$
f(\gamma _1, \gamma _2)=\gamma _1b(\gamma _2)+b(\gamma _1)
-b(\gamma _1\gamma _2), \qquad \gamma _1, \gamma _2 \in \Gamma .
$$
From $f(1, 1)=0$, we get $b(1)=0$. Since $f$ satisfies the parabolic condition,
we have
$$
p\gamma _1b(\gamma _2)+b(p\gamma _1)-b(p\gamma _1\gamma _2)
=p\gamma _1b(\gamma _2)+pb(\gamma _1)-pb(\gamma _1\gamma _2)
$$
for $p \in P$. Taking $\gamma _2=\gamma _1^{-1}$ and
writing $\gamma _1$ as $\gamma$, we find
$$
b(p\gamma )=pb(\gamma )+b(p), \qquad p \in P, \ \gamma \in \Gamma .
$$
Now
\begin{equation*}
\begin{aligned}
f(\sigma , \mu ) & =\sigma b(\mu )+b(\sigma )-b(\sigma\mu )
=\sigma b(\mu )+b(\sigma )-b(\mu ^{-1}\sigma ) \\
& =\sigma b(\mu )+b(\sigma )-\mu ^{-1}b(\sigma )-b(\mu ^{-1}).
\end{aligned}
\end{equation*}
Since $b(\mu ^{-1})=-\mu ^{-1}b(\mu )$, we obtain
\begin{equation}
f(\sigma , \mu )=(1-\mu ^{-1})b(\sigma )+(\sigma +\mu ^{-1})b(\mu ).
\tag{4.8}
\end{equation}

First we consider the case $l_1 \neq l_2$. Since $b \vert P \in Z^1(P, V)$
and $H^1(P, V)=0$ (Theorem 3.7), there exists $\mathbf b \in V$
such that $b(\mu )=(\mu -1)\mathbf b$. Then we have
$$
f(\sigma , \mu )=(1-\mu ^{-1})[b(\sigma )+(1-\sigma )\mathbf b].
$$
We have $\mu ^{-1}(\mathbf e_i \otimes \mathbf e_{i-(l_1-l_2)/2}^{\prime})
=N(\epsilon )^{l_1}(\mathbf e_i \otimes \mathbf e_{i-(l_1-l_2)/2}^{\prime})$.
Hence $c_i$ vanishes for all $i$. This is a contradiction and the proof is complete
in this case.

Next we consider the case $l_1=l_2$. By Theorem 3.7, there exist
$\mathbf b \in V$ and $\mathbf b_0 \in V^U$ such that
$$
b(\mu )=(\mu -1)\mathbf b +\mathbf b_0.
$$
Then we have
$$
f(\sigma , \mu )=(1-\mu ^{-1})[b(\sigma )+(1-\sigma )\mathbf b]
+(\sigma +\mu ^{-1})\mathbf b_0.
$$
Since $\mathbf b_0 \in V^U$, this formula shows that
$c_i=0$ if $i \neq 1, l_1+1$. This is a contradiction and completes the proof.

\proclaim Proposition 4.3. Let $k=(k_1, k_2)$, $k_1 \geq k_2$,
$k_1 \equiv k_2 \equiv 0 \mod 2$. Let $\Omega \in S_k(\Gamma )$ and
let $f=f(\Omega )$ be the normalized parabolic $2$-cocycle attached to
$\Omega$ (cf. (4.1)). We assume that the class number of $F$
in the narrow sense is $1$ and that $\Omega$ is a nonzero Hecke eigenform.
If $k_1 \neq k_2$, we assume $k_2 \geq 4$. If $k_1=k_2$, we assume
$k_2 \geq 6$. Then the cohomology class of $f$ in $H^2(\Gamma , V)$
is non-trivial.
\par
{\bf Proof.} Let $k_1=l_1+2$, $k_2=l_2+2$.
By (2.30), we see that the coefficient $c_i$ of
$\mathbf e_i \otimes \mathbf e_{i-(l_1-l_2)/2}^{\prime}$ in $f(\sigma , \mu )$
is $L(l_1+2-i, \Omega )$ times a nonzero constant for
$(l_1-l_2)/2 +1\leq i \leq (l_1+l_2)/2+1$. It is well known that
$L(s, \Omega ) \neq 0$ for $\Re (s) \geq (k_1+1)/2$ (cf. [Sh4], Proposition 4.16).
For $i=(l_1-l_2)/2+1$, $c_i$ is nonzero times $L((k_1+k_2)/2-1, \Omega )$.
Since $(k_1+k_2)/2-1 \geq (k_1+1)/2$ if $k_2 \geq 3$,
our assertion follows from Lemma 4.2 if $k_1 \neq k_2$.
Assume $k_1=k_2$. For $i=2$, $c_i$ is nonzero times $L(k_1-2, \Omega )$.
Since $k_1-2 \geq (k_1+1)/2$ if $k_1 \geq 5$,
our assertion in this case also follows from Lemma 4.2.

\medskip

{\bf 4.4.} With a free group $\mathcal F$ with finitely many generators,
we write $\Gamma =\mathcal F/R$.
Let $\pi : \mathcal F \longrightarrow \Gamma$ be the canonical homomorphism
with ${\rm Ker}(\pi )=R$. For $g \in \mathcal F$, we put
$\pi (g)=\bar g$. We regard $V$ as an $\mathcal F$-module by
$gv=\bar g v$, $g \in \mathcal F$, $v \in V$. By (1.11), we have
\begin{equation}
H^2(\Gamma ,V) \cong H^1(R,V)^{\Gamma}/{\rm Im}(H^1(\mathcal F,V)).
\tag{4.9}
\end{equation}

We are going to examine the part of the right-hand side of (4.9)
which corresponds to $H_{\rm P}^2(\Gamma , V)$.
Put $\mathcal P=\pi ^{-1}(P)$. Let $f \in Z_{\rm P}^2(\Gamma , V)$.
Take a $1$-cochain $a \in C^1(\mathcal F, V)$ which satisfies (1.12).
Then we have
$$
f(\bar p \bar g_1, \bar g_2)=pg_1a(g_2)+a(pg_1)-a(pg_1g_2),
\qquad p \in \mathcal P, \ g_1, g_2 \in \mathcal F.
$$
By the parabolic condition on $f$, this is equal to
$$
p(g_1a(g_2)+a(g_1)-a(g_1g_2)).
$$
Hence we have
$$
a(pg_1g_2)-a(pg_1)=pa(g_1g_2)-pa(g_1),
\qquad p \in \mathcal P, \ g_1, g_2 \in \mathcal F.
$$
Taking $g_1=g_2^{-1}=g$, we obtain
\begin{equation}
a(pg)=pa(g)+a(p), \qquad p \in \mathcal P, \ g \in \mathcal F.
\tag{4.10}
\end{equation}
Conversely if $a$ satisfies (4.10), then $f$ satisfies the parabolic condition.

Let $\varphi =a\vert R$. We note that $a$ satisfies (1.13) and
$\varphi \in H^1(R, V)^{\Gamma}$. For every $s \in P$, we take an element
$\widetilde s \in \mathcal P$ such that $\pi (\widetilde s)=s$.
We fix the choice of $\widetilde s$. Then we write $a(\widetilde s)$ as
$\widetilde a(s)$. By (1.13), we have
\begin{equation}
a(\widetilde sr)=s\varphi (r)+\widetilde a(s), \qquad s \in P, \ r \in R.
\tag{4.11}
\end{equation}
Now for $s_1$, $s_2 \in P$ and
$r_1$, $r_2 \in R$, we have
\begin{equation*}
\begin{aligned}
& a(\widetilde s_1r_1\widetilde s_2r_2)
=a((\widetilde {s_1s_2})(\widetilde {s_1s_2})^{-1}
\widetilde s_1\widetilde s_2\widetilde s_2^{-1}r_1\widetilde s_2r_2) \\
=  \, & s_1s_2\varphi ((\widetilde {s_1s_2})^{-1}
\widetilde s_1\widetilde s_2\widetilde s_2^{-1}r_1\widetilde s_2r_2)
+\widetilde a(s_1s_2) \\
= \, & s_1s_2[\varphi (\widetilde s_2^{-1}r_1\widetilde s_2)+\varphi (r_2)
+\varphi ((\widetilde {s_1s_2})^{-1}\widetilde s_1\widetilde s_2)]
+\widetilde a(s_1s_2) \\
= \, & s_1\varphi (r_1)+s_1s_2\varphi (r_2)
+\varphi (\widetilde s_1\widetilde s_2(\widetilde {s_1s_2})^{-1})
+\widetilde a(s_1s_2),
\end{aligned}
\end{equation*}
using (1.13), (1.14) and (4.11). On the other hand, by (4.10), we have
\begin{equation*}
\begin{aligned}
& a(\widetilde s_1r_1\widetilde s_2r_2)
=s_1a(\widetilde s_2r_2)+a(\widetilde s_1r_1) \\
= \, & s_1(s_2\varphi (r_2)+\widetilde a(s_2))
+s_1\varphi (r_1)+\widetilde a(s_1).
\end{aligned}
\end{equation*}
Comparing two results, we obtain
\begin{equation}
\varphi (\widetilde s_1\widetilde s_2(\widetilde {s_1s_2})^{-1})
=s_1\widetilde a(s_2)+\widetilde a(s_1)-\widetilde a(s_1s_2).
\tag{4.12}
\end{equation}
The condition (4.12) can be interpreted as follows.
The group extension
\begin{equation*}
\begin{CD}
1 @> >> R @> >> \mathcal P @> >> P @> >> 0.
\end{CD}
\end{equation*}
defines the factor set
\begin{equation}
(s_1, s_2) \longrightarrow \widetilde s_1\widetilde s_2(\widetilde {s_1s_2})^{-1}
\tag{4.13}
\end{equation}
of $P$ taking values in $R$. Mapping this factor set by $\varphi$, we obtain
a $2$-cocycle of $P$ taking values in $V$ (cf. Lemma 1.4). Then (4.12) means that
this $2$-cocycle splits. 

Next we consider the condition (4.10) on a double coset
$\mathcal P\delta R$, where $\delta$ is an arbitrary element of $\mathcal F$.
Since $R$ is a normal subgroup of $\mathcal P$, we have
$$
\mathcal P\delta R=\sqcup _{\widetilde s \in \mathcal P/R}
\ \widetilde s\delta R.
$$
We assume that $\mathcal P\delta R \neq \mathcal PR$.
We write $a(\widetilde s\delta )$ as $\widetilde a(s\delta )$.
By (1.13), we have
\begin{equation}
a(\widetilde s\delta r)=s\delta \varphi (r)+\widetilde a(s\delta ),
\qquad s \in P, \ r \in R.
\tag{4.14}
\end{equation}
Now for $s_1$, $s_2 \in P$ and
$r_1$, $r_2 \in R$, we have
\begin{equation*}
\begin{aligned}
& a(\widetilde s_1r_1\widetilde s_2\delta r_2)
=a((\widetilde {s_1s_2})\delta \delta ^{-1}(\widetilde {s_1s_2})^{-1}
\widetilde s_1\widetilde s_2\delta\delta ^{-1}\widetilde s_2^{-1}r_1\widetilde s_2\delta r_2) \\
=  \, & s_1s_2\delta [\varphi (\delta ^{-1}(\widetilde {s_1s_2})^{-1}
\widetilde s_1\widetilde s_2\delta ^{-1})
+\varphi (\delta ^{-1}\widetilde s_2^{-1}r_1\widetilde s_2\delta ) +\varphi (r_2)]
 +\widetilde a(s_1s_2\delta )\\
=  \, & s_1s_2 \varphi ((\widetilde {s_1s_2})^{-1}\widetilde s_1\widetilde s_2)
+s_1\varphi (r_1) +s_1s_2\delta \varphi (r_2)+\widetilde a(s_1s_2\delta )\\
=  \, & \varphi (\widetilde s_1\widetilde s_2(\widetilde {s_1s_2})^{-1})
+s_1\varphi (r_1) +s_1s_2\delta \varphi (r_2)+\widetilde a(s_1s_2\delta )
\end{aligned}
\end{equation*}
using (1.13), (1.14) and (4.12).
On the other hand, we have, using (4.10),
\begin{equation*}
\begin{aligned}
& a(\widetilde s_1r_1\widetilde s_2\delta r_2)
=s_1a(\widetilde s_2\delta r_2)+a(\widetilde s_1r_1) \\
=  \, & s_1[s_2\delta \varphi (r_2)+a(\widetilde s_2\delta )]
+a(\widetilde s_1r_1) \\
=  \, & s_1s_2\delta \varphi (r_2)+s_1\widetilde a(s_2\delta )
+s_1\varphi (r_1) +\widetilde a(s_1).
\end{aligned}
\end{equation*}
Comparing two formulas, we obtain
\begin{equation}
\varphi (\widetilde s_1\widetilde s_2(\widetilde {s_1s_2})^{-1})
=s_1\widetilde a(s_2\delta )+\widetilde a(s_1)-\widetilde a(s_1s_2\delta ).
\tag{4.15}
\end{equation}

Thus we have shown the following:
For $f \in Z_{\rm P}^2(\Gamma , V)$, take $a \in C^1(\mathcal F, V)$
which satisfies (1.12). Define $\varphi$ and $\widetilde a$ as above.
Then (4.12) holds. Conversely take
$\varphi \in H^1(R, V)^{\Gamma}$. Suppose that the $2$-cocycle
$$
(s_1, s_2) \longrightarrow \varphi (\widetilde s_1\widetilde s_2(\widetilde {s_1s_2})^{-1})
$$
of $P$ taking values in $V$ splits. There exists $\widetilde a \in C^1(P,V)$
with which (4.12) holds. Take a double coset decomposition
$\mathcal F=\sqcup _{\delta} \mathcal P\delta R$. 
We put $\widetilde a(s\delta )=\widetilde a(s)$.
Then we define $a \in C^1(\mathcal F, V)$ by the formula
$$
a(\widetilde s\delta r)=s\delta \varphi (r)+\widetilde a(s\delta ), \qquad r \in R.
$$
Then $a$ satisfies (4.10). Therefore the $2$-cocycle $f$ determined by (1.12)
belongs to $Z_{\rm P}^2(\Gamma , V)$.
Thus we have proved the following proposition.

\proclaim Proposition 4.4. On the right-hand side of (4.9),
the subgroup which corresponds to $H_{\rm P}^2(\Gamma , V)$
consists of the class of $\varphi \in H^1(R, V)^{\Gamma}$ for which the $2$-cocycle 
$(s_1, s_2) \mapsto \varphi (\widetilde s_1\widetilde s_2(\widetilde {s_1s_2})^{-1})$
of $P$ taking values in $V$ splits.
\par

By Theorem 3.9, we have $H^2(P, V)=0$ if $l_1 \neq l_2$.
Hence the next proposition follows.

\proclaim Proposition 4.5. If $l_1 \neq l_2$, then we have
$H^2(\Gamma , V)=H_{\rm P}^2(\Gamma , V)$.
\par

It is known that there are no holomorphic Eisenstein series of weight
$(k_1, k_2)$ if $k_1 \neq k_2$ ([Sh6], Proposition 2.1).
We can interpret this proposition as the cohomological counter part of this fact.

\medskip

{\bf Remark 4.6.} In view of the results of Matsushima-Shimura [MS],
Hida [Hi1], [Hi2] and Harder [Ha], we should be able to prove that
$\dim H_P^2(\Gamma , V)=4\dim S_{l_1+2, l_2+2}(\Gamma )$.
The author does not work out the details yet.
The parabolic cohomology group is also discussed in [Hi2].

\bigskip
\newpage

\centerline{\S5. Decompositions of $H^2(\Gamma , V)$}

\medskip

{\bf 5.1.} Let $F$ be a real quadratic field and let $\Gamma ={\rm PSL}(2, {\mathcal O}_F)$.
We define elements $\sigma$, $\mu$, $\tau$ and $\eta$ of $\Gamma$ by
$$
\sigma =\begin{pmatrix} 0 & 1 \\ -1 & 0 \end{pmatrix}, \quad
\mu =\begin{pmatrix} \epsilon & 0 \\ 0 & \epsilon ^{-1} \end{pmatrix}, \quad
\tau =\begin{pmatrix} 1 & 1 \\ 0 & 1 \end{pmatrix}, \quad
\eta =\begin{pmatrix} 1 & \omega \\ 0 & 1 \end{pmatrix}.
$$
Here we choose an $\omega$ so that
$\mathcal O_F=\mathbf Z+\mathbf Z\omega$.
Let $\mathcal F$ be the free group on four letters
$\widetilde\sigma$, $\widetilde\mu$, $\widetilde\tau$, $\widetilde\eta$.
Let $\pi : \mathcal F \longrightarrow \Gamma$ be the homomorphism
such that
$$
\pi (\widetilde\sigma )=\sigma , \quad \pi (\widetilde\mu )=\mu , \quad
\pi (\widetilde\tau )=\tau , \quad \pi (\widetilde\eta )=\eta .
$$
By Vaser$\rm {\check{s}}$tein [V],  $\pi$ is surjective.
Let $R$ be the kernel of $\pi$.
For $\gamma \in \Gamma$, we choose a $\widetilde\gamma \in \mathcal F$
so that $\pi (\widetilde\gamma )=\gamma$. For $\gamma =\sigma$, $\mu$,
$\tau$ and $\eta$, we choose $\widetilde\gamma$ so that the notation to be
consistent. We choose $\widetilde 1=1$.
For other $\gamma$, we will specify the choice of $\widetilde\gamma$ later
(cf. (5.2) and \S6.2).

Let $f \in Z^2(\Gamma , V)$ be a normalized $2$-cocycle.
There exists $a \in C^1(\mathcal F, V)$ which satisfies
$$
f(\gamma _1, \gamma _2)=\gamma _1a(\widetilde\gamma _2)
+a(\widetilde\gamma _1)-a(\widetilde\gamma _1\widetilde\gamma _2).
$$
A corresponding element $\varphi \in H^1(R,V)^{\Gamma}$ to $f$ is given by
$\varphi =a\vert R$. As was shown in \S1.5, adding a cobounday to $f$,
we may assume that $f \in Z^2(\Gamma , V)$ is given by
\begin{equation}
f(\gamma _1, \gamma _2)=
-\varphi(\widetilde\gamma _1\widetilde\gamma _2(\widetilde{\gamma _1\gamma _2})^{-1}),
\qquad \gamma _1, \gamma _2 \in \Gamma .
\tag{5.1}
\end{equation}

Let $\mathcal F_P$ be the subgroup of $\mathcal F$ generated by
$\widetilde\mu$, $\widetilde\tau$ and $\widetilde\eta$.
Let $\pi _P$ be the restriction of $\pi$ to $\mathcal F_P$
and let $R_P$ be the kernel of $\pi _P$.
We see that $R_P$ is generated by the elements corresponding to
the relations (iv), (v), (vi) of Appendix and their conjugates.
Suppose that $f$ satisfies the parabolic condition (4.3).
Then, by (4.12), we see that we may assume that $\varphi \vert R_P=0$
in addition to (5.1), adding a cobounday to $f$ if necessary.

Conversely assume that $\varphi \vert R_P=0$.
Take a complete set of representatives $\Delta$ for $P\backslash\Gamma$
and fix it. We have
$$
\Gamma =\sqcup _{\delta \in \Delta} P\delta .
$$
For $\gamma =p\delta$, $p \in P$, $\delta \in \Delta$, we define
\begin{equation}
\widetilde\gamma =\widetilde p\widetilde\delta .
\tag{5.2}
\end{equation}
In (5.1), write $\gamma _1=p_1\delta _1$, $p_1 \in P$, $\delta _1\in \Delta$,
$\gamma _1\gamma _2=p_2\delta _2$, $p_2 \in P$, $\delta _2 \in \Delta$.
Let $p \in P$. Then we have
$$
\widetilde {p\gamma _1}=\widetilde {pp_1}\widetilde\delta _1
=\widetilde{pp_1}(\widetilde p\widetilde p_1)^{-1}\widetilde p\widetilde\gamma _1,
\qquad \widetilde{p\gamma _1\gamma _2}=\widetilde{pp_2}\widetilde p_2^{-1}
\widetilde{\gamma _1\gamma _2}.
$$
Hence, by (5.1), we have
\begin{equation*}
\begin{aligned}
f(p\gamma _1, \gamma _2)
& =-\varphi (\widetilde{pp_1}(\widetilde p\widetilde p_1)^{-1}\widetilde p
\widetilde\gamma _1\widetilde\gamma _2
\{ \widetilde {pp_2}(\widetilde p\widetilde p_2)^{-1}\widetilde p
\widetilde {\gamma _1\gamma _2} \}^{-1}) \\
& =-\varphi (\widetilde p\widetilde\gamma _1\widetilde\gamma _2
(\widetilde{\gamma _1\gamma _2})^{-1}\widetilde p^{-1})
=-p\varphi (\widetilde\gamma _1\widetilde\gamma _2
(\widetilde{\gamma _1\gamma _2})^{-1})
=pf(\gamma _1, \gamma _2).
\end{aligned}
\end{equation*}
Therefore $f$ satisfies the parabolic condition (4.3).

The value $f(\sigma , \mu )$ of the cocycle is related to the critical values of
the $L$-function. By (5.1), we have
$$
f(\sigma , \mu )=-\varphi (\widetilde\sigma \widetilde\mu(\widetilde{\sigma\mu})^{-1})
=-\varphi (\widetilde\sigma \widetilde\mu
(\widetilde{\mu ^{-1}\sigma})^{-1}).
$$
We assume that $\sigma \in \Delta$. Then we have
$$
f(\sigma , \mu )=-\varphi (\widetilde\sigma \widetilde\mu
\widetilde\sigma ^{-1}(\widetilde{\mu^{-1}})^{-1}),
$$
since $\widetilde{\mu ^{-1}\sigma}=\widetilde{\mu ^{-1}}\widetilde\sigma$.
As $\widetilde{\mu ^{-1}}\widetilde\mu \in R_P$, we have
\begin{equation*}
\begin{aligned}
f(\sigma , \mu ) & =-\varphi (\widetilde\sigma \widetilde\mu \widetilde\sigma ^{-1}
\widetilde\mu )
=-\varphi (\widetilde\sigma \widetilde\mu \widetilde\sigma ^{-2}
\widetilde\sigma \widetilde\mu )
=-\varphi (\widetilde\sigma \widetilde\mu \widetilde\sigma ^{-2}
(\widetilde\sigma \widetilde\mu )^{-1}
\widetilde\sigma \widetilde\mu \widetilde\sigma \widetilde\mu ) \\
& = -\sigma\mu \varphi (\widetilde\sigma ^{-2})
-\varphi (\widetilde\sigma \widetilde\mu \widetilde\sigma \widetilde\mu ).
\end{aligned}
\end{equation*}
Therefore we obtain
\begin{equation}
f(\sigma , \mu )=-\varphi ((\widetilde\sigma \widetilde\mu )^2)
+\sigma\mu \varphi (\widetilde\sigma ^2).
\tag{5.3}
\end{equation}

{\bf 5.2.} Let us consider the action of Hecke operators.
Let $\varpi$ be a totally positive element of $F$. Let
$$
\Gamma \begin{pmatrix} 1 & 0 \\ 0 & \varpi \end{pmatrix} \Gamma
=\sqcup _{i=1}^d \Gamma \beta _i
$$
be a coset decomposition. We put
$$
c=\prod _{\nu =1}^2 (\varpi ^{(\nu )})^{(k_0+k_{\nu})/2 -2}.
$$
Let $f \in Z^2(\Gamma , V)$ and put $g=cT(\varpi )f$. The explicit
form of $g$ is given as follows (cf. Proposition 1.3 and (2.45)). Let
$$
\beta _i\gamma _1=\delta _i^{(1)}\beta _{j(i)}, \quad \delta _i^{(1)} \in \Gamma , \qquad
\beta _i\gamma _2=\delta _i^{(2)}\beta _{k(i)}, \quad \delta _i^{(2)} \in \Gamma ,
$$
for $1 \leq i \leq d$. Here $j$ and $k$ are permutations on $d$ letters. Then
\begin{equation}
g(\gamma _1, \gamma _2)=c\sum _{i=1}^d \beta _i^{-1}
f(\beta _i\gamma _1\beta _{j(i)}^{-1}, \beta _{j(i)}\gamma _2\beta _{k(j(i))}^{-1}).
\tag{5.4}
\end{equation}
We assume that $f \in Z^2_P(\Gamma , V)$ and that it is given by (5.1) with
$\varphi \in H^1(R, V)^{\Gamma}$ satisfying $\varphi \vert R_P=0$. Then we have
\begin{equation}
g(\gamma _1, \gamma _2)=-c\sum _{i=1}^d \beta _i^{-1}
\varphi (\widetilde {\beta _i\gamma _1\beta _{j(i)}^{-1}}
\widetilde{\beta _{j(i)}\gamma _2\beta _{k(j(i))}^{-1}}
(\widetilde{\beta _i\gamma _1\gamma _2\beta _{k(j(i))}^{-1}})^{-1}).
\tag{5.5}
\end{equation}
Let $\psi \in H^1(R, V)^{\Gamma}$ be a corresponding element to $g$.
We are going to give an explicit form of $\psi$.
There exists $b \in C^1(\mathcal F, V)$ such that
$$
g(\bar x_1, \bar x_2)=x_1b(x_2)+b(x_1)-b(x_1x_2), \qquad
x_1, x_2 \in \mathcal F
$$
and $\psi$ is given as the restriction of $b$ to $R$. Here $\bar x=\pi (x)$, $x \in \mathcal F$.
We assume that $(\varpi )$ is a prime ideal.
Then $d=N(\varpi )+1$ and $\{ \beta _i \}$ can be taken as
$$
\left\{ \begin{pmatrix} 1 & u \\ 0 & \varpi \end{pmatrix} , u \mod \varpi , \quad
\begin{pmatrix} \varpi & 0 \\ 0 & 1 \end{pmatrix} \right\} .
$$
Take $p \in P$ and let $\beta _ip\beta _{j(i)}^{-1} \in \Gamma$ for
$1\leq i \leq d$. Then we see easily that
\begin{equation}
\beta _i p\beta _{j(i)}^{-1} \in P, \qquad 1 \leq i \leq d.
\tag{5.6}
\end{equation}
By (5.5), (5.6) and $\varphi\vert R_P=0$, we find
\begin{equation}
g(p_1, p_2)=0, \qquad p_1, p_2 \in P.
\tag{5.7}
\end{equation}
We have
$$
b(x_1x_2)=x_1b(x_2)+b(x_1)-g(\bar x_1, \bar x_2), \qquad
x_1, x_2 \in \mathcal F
$$
and we can use this formula to determine the value $b(x)$, $x \in \mathcal F$
by the induction on the length of the element $x$. As the initial conditions,
we may assume that
$$
b(\widetilde\mu )=0, \quad b(\widetilde\tau )=0, \quad b(\widetilde\eta )=0, \quad
b(\widetilde\sigma )=0.
$$
Then, by (5.7), we see that
\begin{equation}
b\vert \mathcal F_P=0.
\tag{5.8}
\end{equation}
The next Proposition is a special case of Proposition 1.5.

\proclaim Proposition 5.1. Suppose $\gamma _j\in \Gamma$ are
given for $1 \leq j \leq m$. For every $j$, we define $p_j \in S_d$ by
$$
\beta _i\gamma _j\beta _{p_j(i)}^{-1} \in \Gamma , \qquad 1 \leq i \leq d.
$$
We define $q_j \in S_d$ inductively by
$$
q_1=p_1, \qquad q_k=p_kq_{k-1}, \quad 2 \leq k \leq m.
$$
We assume that $\gamma _j \in P$ or $\gamma _j=\sigma$ for every $j$.
Then we have
\begin{equation*}
\begin{aligned}
& b(\widetilde\gamma _1\widetilde\gamma _2\cdots\widetilde\gamma _m) \\
= \, & c\sum _{i=1}^d \beta _i^{-1}
\varphi (\widetilde{\beta _i\gamma _1\beta _{q_1(i)}^{-1}}
\widetilde{\beta _{q_1(i)}\gamma _2\beta _{q_2(i)}^{-1}}\cdots
\widetilde{\beta _{q_{m-1}(i)}\gamma _m\beta _{q_m(i)}^{-1}}
(\widetilde{\beta _i\gamma _1\gamma _2\cdots \gamma _m\beta _{q_m(i)}^{-1}})^{-1}).
\end{aligned}
\tag{5.9}
\end{equation*}
\par

\medskip

{\bf 5.3.} For the practical computation, it is convenient to decompose
$H^2(\Gamma , V)$ into a direct sum of subspaces under the action
of the automorphisms of $\Gamma$. We put
$$
Z=\left\{ \begin{pmatrix} u & 0 \\ 0 & u \end{pmatrix} \bigg| \ u \in E_F \right\} ,
$$
which is the center of ${\rm GL}(2, {\mathcal O}_F)$. Then we have
$$
Z\cdot {\rm SL}(2, {\mathcal O}_F)/Z \cong
{\rm SL}(2, {\mathcal O}_F)/\{ \pm\begin{pmatrix} 1 & 0 \\ 0 & 1 \end{pmatrix} \}
={\rm PSL}(2, {\mathcal O}_F)=\Gamma .
$$
By this isomorphism, we regard $\Gamma$ as a subgroup of
${\rm PGL}(2, {\mathcal O}_F)={\rm GL}(2, {\mathcal O}_F)/Z$.
Hereafter we assume that $l_1$ and $l_2$ are even. When $l$ is even,
we define a representation $\rho _l^{\prime}$ of ${\rm GL}(2, {\mathbf C})$ by
$$
\rho _l^{\prime}(g)=\rho _l(g)\det (g)^{-l/2},
\qquad g \in {\rm GL}(2, {\mathbf C}).
$$
Then $\rho _l^{\prime}$ is trivial on the center.
We put $\rho ^{\prime}=\rho _{l_1}^{\prime}\otimes\rho _{l_2}^{\prime}$.
By $gv=\rho ^{\prime}(g)v$, $g \in {\rm GL}(2, {\mathcal O}_F)$,
$v \in V$, we regard $V$ as a left ${\rm GL}(2, {\mathcal O}_F)$-module.
Since $\rho ^{\prime}(z)={\rm id}$, $z \in Z$, we can regard $V$ as
a ${\rm PGL}(2, {\mathcal O}_F)$-module.
Since $\rho ^{\prime}\vert\Gamma =\rho\vert\Gamma$,
the $\Gamma$-module structure of $V$ is the same as before.

We have
$$
{\rm PGL}(2, {\mathcal O_F})/{\rm PSL}(2, {\mathcal O}_F) \cong
E_F/E_F^2 \cong {\mathbf Z/2\mathbf Z}\oplus {\mathbf Z/2\mathbf Z}.
$$
By conjugation, ${\rm PGL}(2, {\mathcal O}_F)$ acts on
$H^2(\Gamma , V)$ and it decomposes into a direct sum
of four subspaces. We put
$$
\nu =\begin{pmatrix} \epsilon & 0 \\ 0 & 1 \end{pmatrix} , \qquad
\delta =\begin{pmatrix} -1 & 0 \\ 0 & 1 \end{pmatrix} .
$$
We see that ${\rm PGL}(2, {\mathcal O}_F)$ is generated by
$\nu$ and $\delta$ over ${\rm PSL}(2, {\mathcal O}_F)$.
We first examine the action of $\nu$. For $f \in Z^2(\Gamma , V)$,
define $\widetilde ef \in Z^2(\Gamma , V)$ by (cf. (1.3))
\begin{equation}
\widetilde ef(\gamma _1, \gamma _2)
=\nu ^{-1}f(\nu\gamma _1\nu ^{-1}, \nu\gamma _2\nu ^{-1}),
\qquad \gamma _1, \gamma _2 \in \Gamma .
\tag{5.10}
\end{equation}
Then $\widetilde e$ induces an automorphism $e$ of $H^2(\Gamma , V)$.
Since $\nu ^2=\mu$, $\widetilde e^2$ is obtained from the inner automorphism
by $\mu$. Hence $e^2=1$. By (5.10), we see that
$\widetilde ef$ is a parabolic cocycle if $f$ is parabolic.
Therefore, by the action of $e$, we have the decompositions
$$
H^2(\Gamma , V)=H^2(\Gamma , V)^+ \oplus H^2(\Gamma , V)^-, \quad
H_P^2(\Gamma , V)=H_P^2(\Gamma , V)^+ \oplus H_P^2(\Gamma , V)^-.
$$
Here we put
$$
H^2(\Gamma , V)^{\pm}=\{ c \in H^2(\Gamma , V) \mid ec=\pm c \} , \
H_P^2(\Gamma , V)^{\pm}=\{ c \in H_P^2(\Gamma , V) \mid ec=\pm c \} .
$$
Explicitly the decomposition is given by
$$
f=\frac {1}{2}\big[(1+\widetilde e)f+(1-\widetilde e)f\big],
\qquad f \in Z^2(\Gamma , V).
$$

\proclaim Proposition 5.2. Let $k=(k_1, k_2)$, $k_1 \geq k_2$,
$k_1$ and $k_2$ are even. Let $\Omega \in S_k(\Gamma )$ and
let $f=f(\Omega )$ be the normalized parabolic $2$-cocycle attached to
$\Omega$ by (4.1). We assume that the class number of $F$
in the narrow sense is $1$ and that $\Omega$ is a nonzero Hecke eigenform.

(1) If $k_1 \neq k_2$, we assume $k_2 \geq 6$. If $k_1=k_2$, we assume
$k_2 \geq 8$. Then the cohomology class of $(1+\widetilde e)f$ in
$H^2(\Gamma , V)$ is non-trivial.

(2) If $k_1 \neq k_2$, we assume $k_2 \geq 4$. If $k_1=k_2$, we assume
$k_2 \geq 6$. Then the cohomology class of $(1-\widetilde e)f$ in
$H^2(\Gamma , V)$ is non-trivial.

\par
{\bf Proof.} We apply Lemma 4.2 in a similar way to the proof of
Proposition 4.3. We use the same notation as there.
By (5.10), we have
$$
(\widetilde ef)(\sigma , \mu )=\nu ^{-1}f(\nu\sigma\nu ^{-1}, \nu\mu\nu ^{-1})
=\nu ^{-1}f(\mu\sigma , \mu )=\nu ^{-1}\mu f(\sigma , \mu )
=\nu f(\sigma , \mu ).
$$
We have
$$
\nu (\mathbf e_i \otimes \mathbf e_{i-(l_1-l_2)/2}^{\prime})
=N(\epsilon )^{l_1/2+1-i}(\mathbf e_i \otimes \mathbf e_{i-(l_1-l_2)/2}^{\prime})
=N(\epsilon )^{k_1/2-i}(\mathbf e_i \otimes \mathbf e_{i-(l_1-l_2)/2}^{\prime}).
$$
By the assumption, we have $N(\epsilon )=-1$. The range of $i$ is
$\frac {k_1}{2} -\frac {l_2}{2} \leq i \leq \frac {k_1}{2} +\frac {l_2}{2}$.
We see that $L(l_1+2-i, \Omega )$ is non-vanishing if $i \neq k_1/2$.
To conclude the non-vanishing of the cohomology class of
$(1+\widetilde e)f$, it suffices to find an even integer $j$ such that
$0<j\leq l_2/2$ if $k_1 \neq k_2$ and $0<j\leq l_2/2-1$ if $k_1=k_2$.
Such a $j$ exists under the condition stated in (1).
To conclude the non-vanishing of the cohomology class of
$(1-\widetilde e)f$, it suffices to find an odd integer $j$ such that
$0<j\leq l_2/2$ if $k_1 \neq k_2$ and $0<j\leq l_2/2-1$ if $k_1=k_2$.
Such a $j$ exists under the condition stated in (2). This completes the proof.

\medskip

We put
$$
\overline\Gamma ^{\ast}=\{ \gamma \in {\rm GL}(2, {\mathcal O}_F) \mid
\det (\gamma )=\epsilon ^n, n \in {\mathbf Z} \} , \qquad
\Gamma ^{\ast}=Z\overline\Gamma ^{\ast}/Z.
$$
Then $\Gamma ^{\ast}$ is generated by $\nu$ over $\Gamma$ and we have
$[\Gamma ^{\ast}:\Gamma ]=2$.
Let
$$
{\rm Res}: H^2(\Gamma ^{\ast}, V) \longrightarrow H^2(\Gamma , V), \qquad
T: H^2(\Gamma , V) \longrightarrow H^2(\Gamma ^{\ast}, V)
$$
be the restriction map and the transfer map respectively.

\proclaim Proposition 5.3. We have

(1) ${\rm Res}(H^2(\Gamma ^{\ast}, V))=H^2(\Gamma , V)^+$.

(2) $T(H^2(\Gamma , V)^+)=H^2(\Gamma ^{\ast}, V)$.

(3) ${\rm Ker}(T)=H^2(\Gamma , V)^-$.
\par
{\bf Proof.} It is clear that
${\rm Res}(H^2(\Gamma ^{\ast}, V))\subseteq H^2(\Gamma , V)^+$.
Let $f \in Z^2(\Gamma , V)$ and take a coset decomposition
$\Gamma ^{\ast}=\Gamma \sqcup \nu ^{-1}\Gamma$. Then by
Proposition 1.2, we have
$$
\widetilde T(f)(\gamma _1, \gamma _2)
=f(\gamma _1\nu ^a, \nu ^{-a}\gamma _2\nu ^b)
+\nu ^{-1}f(\nu\gamma _1\nu ^c, \nu ^{-c}\gamma _2\nu ^d)
$$
for $\gamma _1$, $\gamma _2 \in \Gamma ^{\ast}$.
Here $\widetilde T(f)$ denotes a cocycle which represents the transfer of
the class of $f$;
$a$, $b$, $c$, $d=0$ or $-1$ and they are determined so that
all arguments on the right-hand side belong to $\Gamma$.
In particular, if $\gamma _1$, $\gamma _2 \in \Gamma$, then we have
$$
\widetilde T(f)(\gamma _1, \gamma _2)=f(\gamma _1, \gamma _2)
+\nu ^{-1}f(\nu\gamma _1\nu ^{-1}, \nu \gamma _2\nu ^{-1}).
$$
Therefore we have
$$
{\rm Res}\circ T=1+e
$$
on $H^2(\Gamma , V)$. This formula combined with (1.6) shows that
$T\circ {\rm Res}$ and ${\rm Res}\circ T$ are the multiplication by $2$
on $H^2(\Gamma ^{\ast}, V)$ and $H^2(\Gamma , V)^+$ respectively.
Hence (1) and (2) follow. Then (3) follows since
$H^2(\Gamma , V)^- \subset {\rm Ker}(T)$ and $T\vert H^2(\Gamma ,V)^+$
is injective. This completes the proof.

\medskip

{\bf 5.4.} We have
$$
H^2(\Gamma , V)\cong H^1(R, V)^{\Gamma}/{\rm Im}(H^1(\mathcal F, V)).
$$
Let us consider the action of $e$ on the right-hand side under this isomorphism.
We use the same notation as in \S5.1.
Let $\xi$ be the automorphism of $\Gamma$ defined by
$\xi (\gamma )=\nu\gamma\nu ^{-1}$, $\gamma \in \Gamma$. Put
$$
\epsilon =A+B\omega , \qquad \epsilon\omega =C+D\omega .
$$
Then we have
$\begin {pmatrix} A & B \\ C & D \end{pmatrix} \in {\rm GL}(2, {\mathbf Z})$.
We have
$$
\nu\sigma\nu ^{-1}=\sigma\mu ^{-1}, \quad \nu\mu\nu ^{-1}=\mu , \quad
\nu\tau\nu ^{-1}=\tau ^A\eta ^B, \quad \nu\eta\nu ^{-1}=\tau ^C\eta ^D.
$$
Using Lemma 3.11, we can check that there exists an automorphism $\widetilde\xi$
of $\mathcal F$ which satisfies
\begin{equation}
\pi (\widetilde\xi (g))=\xi (\pi (g)), \qquad g \in \mathcal F.
\tag{5.11}
\end{equation}
Now let $f \in Z^2(\Gamma , V)$ and take $a \in C^1(\mathcal F, V)$ so that
$$
f(\pi (g_1), \pi (g_2))=g_1a(g_2)+a(g_1)-a(g_1g_2), \qquad
g_1, g_2 \in \mathcal F.
$$
Then we have
\begin{equation*}
\begin{aligned}
& (\widetilde ef)(\pi (g_1), \pi (g_2))
=\nu ^{-1}f(\xi (\pi (g_1)), \xi (\pi (g_2)))
=\nu ^{-1}f(\pi (\widetilde\xi (g_1)), \pi (\widetilde\xi (g_2)) \\
= & \, g_1\nu ^{-1}a(\widetilde\xi (g_2))
+\nu ^{-1}a(\widetilde\xi (g_1))
-\nu ^{-1}a(\widetilde\xi (g_1g_2))
\end{aligned}
\end{equation*}
for $g_1$, $g_2 \in \mathcal F$. Put
$$
a^{\prime}(g)=\nu ^{-1}a(\widetilde\xi (g)), \qquad g \in \mathcal F.
$$
Then we have
$$
(\widetilde ef)(\pi (g_1), \pi (g_2))
=g_1a^{\prime}(g_2)+a^{\prime}(g_1)-a^{\prime}(g_1g_2), \qquad
g_1, g_2 \in \mathcal F.
$$
Thus we obtain the following proposition.

\proclaim Proposition 5.4. Let $f \in Z^2(\Gamma , V)$ and let
$\varphi \in H^1(R, V)^{\Gamma}$ be a corresponding element.
Then a corresponding element $\psi$ of $H^1(R, V)^{\Gamma}$ to
$\widetilde ef$ is given by
$$
\psi (r)=\nu ^{-1}\varphi (\widetilde\xi (r)), \qquad r \in R.
$$

We can check easily that the map $\varphi \longrightarrow \psi$
induces a map from \newline $H^1(R, V)^{\Gamma}/{\rm Im}(H^1(\mathcal F, V))$
to itself and gives an automorphism of order $2$.

\medskip

{\bf 5.5.} For the actual computation, the cohomology group
$H^2(\Gamma ^{\ast}, V)$ is easier to handle than $H^2(\Gamma , V)$.
By the action of $\delta$, we can further decompose
$H^2(\Gamma ^{\ast}, V)$ so that
$$
H^2(\Gamma ^{\ast}, V)
=H^2(\Gamma ^{\ast}, V)^+\oplus H^2(\Gamma ^{\ast}, V)^-.
$$
Let $\widetilde d$ (resp. $d$) denote the action of $\delta$ on 
$Z^2(\Gamma ^{\ast}, V)$ (resp. $H^2(\Gamma ^{\ast}, V)$).

\proclaim Proposition 5.5. Let $k=(k_1, k_2)$, $k_1 \geq k_2$,
$k_1$ and $k_2$ are even. Let $\Omega \in S_k(\Gamma )$ and
let $f=f(\Omega )$ be the normalized parabolic $2$-cocycle attached to
$\Omega$ by (4.1). We assume that the class number of $F$
in the narrow sense is $1$ and that $\Omega$ is a nonzero Hecke eigenform.
Take $f^{\ast} \in Z^2(\Gamma ^{\ast}, V)$ so that
$f^{\ast}\vert\Gamma =(1+\widetilde e)f$.
\footnote{$f^{\ast}=\widetilde T(f)$ satisfies this condition.}
If $k_1 \neq k_2$, we assume $k_2 \geq 6$. If $k_1=k_2$, we assume
$k_2 \geq 8$. Then the cohomology class of $(1+\widetilde d)f^{\ast}$ in
$H^2(\Gamma ^{\ast}, V)$ is non-trivial.
\par
{\bf Proof.} The proof is similar to that of Proposition 5.2.
We consider the restriction $f_0$ of $(1+\widetilde d)f^{\ast}$ to $\Gamma$.
Since $\delta$ commutes with $\sigma$ and $\mu$, we find
$$
f_0(\sigma , \mu )=(1+\delta )(1+ \nu )f(\sigma , \mu ).
$$
We have
$$
\delta (\mathbf e_i \otimes \mathbf e_{i-(l_1-l_2)/2}^{\prime})
=\mathbf e_i \otimes \mathbf e_{i-(l_1-l_2)/2}^{\prime}.
$$
Hence the assertion follows from Lemma 4.2 in the same way as Proposition 5.2.

Until the end of this subsection,
we assume that $\sigma$, $\nu$ and $\tau$ generate $\Gamma ^{\ast}$.
(This assumption is satisfied if
$\mathcal O_F={\mathbf Z}+{\mathbf Z}\epsilon$.)
Let $\mathcal F^{\ast}$ be the free group on three letters
$\widetilde\sigma$, $\widetilde\nu$ and $\widetilde\tau$.
We define a surjective homomorphism $\pi ^{\ast}$ of $\mathcal F^{\ast}$ onto
$\Gamma ^{\ast}$ by
$$
\pi ^{\ast}(\widetilde\sigma )=\sigma , \qquad
\pi ^{\ast}(\widetilde\nu )=\nu , \qquad \pi ^{\ast}(\widetilde\tau )=\tau
$$
and let $R^{\ast}$ be the kernel of $\pi ^{\ast}$. We see that
$\delta$ commutes with $\sigma$ and $\nu$ and
$\delta\tau\delta ^{-1}=\tau ^{-1}$.
We can define an automorphism 
$x \mapsto x_{\delta}$ of $\mathcal F^{\ast}$ by
$(\widetilde\sigma )_{\delta}=\widetilde\sigma$,
$(\widetilde\nu )_{\delta}=\widetilde\nu$,
$(\widetilde\tau )_{\delta}=\widetilde\tau ^{-1}$. Then we have
$$
\pi ^{\ast}(x_{\delta})=\delta\pi ^{\ast}(x)\delta ^{-1},
\qquad x \in \mathcal F^{\ast}.
$$
The following proposition can be shown in a similar manner to Proposition 5.4.

\proclaim Proposition 5.6. Let $f \in Z^2(\Gamma ^{\ast}, V)$ and let
$\varphi \in H^1(R^{\ast}, V)^{\Gamma ^{\ast}}$ be a corresponding element.
Then a corresponding element $\psi$ of $H^1(R^{\ast}, V)^{\Gamma ^{\ast}}$ to
$\widetilde df$ is given by
$$
\psi (r)=\delta ^{-1}\varphi (r_{\delta}), \qquad r \in R^{\ast}.
$$
\par

Let $\varphi \in H^1(R^{\ast}, V)^{\Gamma ^{\ast}}$. We define
$\varphi _{\delta} \in H^1(R^{\ast}, V)^{\Gamma ^{\ast}}$ by the formula
\begin{equation}
\varphi _{\delta}(r)=\delta ^{-1}\varphi (r_{\delta}).
\tag{5.12}
\end{equation}
Then we can check easily that $(\varphi _{\delta})_{\delta}=\varphi$ and
$H^1(R^{\ast}, V)^{\Gamma ^{\ast}}$ decomposes into
a direct sum of $\pm 1$ eigenspaces under the action of $\delta$:
\begin{equation}
H^1(R^{\ast}, V)^{\Gamma ^{\ast}}
=H^1(R^{\ast}, V)^{\Gamma ^{\ast},+} \oplus
H^1(R^{\ast}, V)^{\Gamma ^{\ast},-} .
\tag{5.13}
\end{equation}

{\bf 5.6.} Let $l_1$ and $l_2$ be nonnegative even integers.
We assume that $l_1 \geq l_2$.
Let $\Omega \in S_{l_1+2,l_2+2}(\Gamma )$.
Define $L(s, \Omega )$ and $R(s, \Omega )$ by (2.4) and (2.5) respectively.
The functional equation is (cf. (2.7))
$$
R(s, \Omega )=(-1)^{(l_1+l_2)/2}R(l_1+2-s, \Omega ).
$$
For an integer $m$, $L(m, \Omega )$ is a critical value if and only if
\begin{equation}
\frac {l_1-l_2}{2}+1 \leq m \leq \frac {l_1+l_2}{2}+1.
\tag{5.14}
\end{equation}
The central critical value is $L(l_1/2+1, \Omega )$ which vanishes
if $(l_1+l_2)/2$ is odd. By (2.30), we have
\begin{equation}
R(m, \Omega )=(-1)^m i^{(l_1-l_2)/2}(2\pi )^{(l_2-l_1)/2}
P_{m-1, m-1-(l_1-l_2)/2}.
\tag{5.15}
\end{equation}
Here $P_{s,t}$ denotes the period integral given by (2.25).
Let $f =f(\Omega )\in Z_P^2(\Gamma , V)$ be the parabolic
$2$-cocycle defined by (4.1). Then we have
$$
f (\sigma, \mu)
=-\int _{i\epsilon ^{-1}}^{i\epsilon}
\int _0^{i\infty} \mathfrak d(\Omega )
$$
and $-P_{m-1, m-1-(l_1-l_2)/2}$ is equal to the coefficient of
${\mathbf e}_{l_1+2-m} \otimes {\mathbf e}^{\prime}_{(l_1+l_2)/2+2-m}$
in $f(\sigma , \mu )$.

Using the operator $\widetilde e$ (cf. (5.10)), we define
$$
f^+=(1+\widetilde e)f, \qquad f^{-}=(1-\widetilde e)f.
$$
We have $f^{\pm} \in Z_P^2(\Gamma , V)$.
As was shown in the proof of Proposition 5.2, we have
\begin{equation}
f^+(\sigma , \mu )=(1+\nu )f(\sigma , \mu ), \qquad
f^-(\sigma , \mu )=(1-\nu )f(\sigma , \mu ).
\tag{5.16}
\end{equation}
We have
\begin{equation}
\nu ({\mathbf e}_{l_1+2-m} \otimes {\mathbf e}^{\prime}_{(l_1+l_2)/2+2-m})
=N(\epsilon )^{m-1-l_1/2}
{\mathbf e}_{l_1+2-m} \otimes {\mathbf e}^{\prime}_{(l_1+l_2)/2+2-m}.
\tag{5.17}
\end{equation}
Assume $N(\epsilon )=-1$.
Suppose that $l_1/2$ is even. By (5.17), we see that $f^+(\sigma , \mu )$
contains information on $R(m, \Omega )$ for odd $m$ and
$f^-(\sigma , \mu )$ contains information on $R(m, \Omega )$ for even $m$.
If $l_1/2$ is odd, then $f^+(\sigma , \mu )$
contains information on $R(m, \Omega )$ for even $m$ and
$f^-(\sigma , \mu )$ contains information on $R(m, \Omega )$ for odd $m$.

To treat $f^-$ efficiently, we will need more techniques which will be
explained in the next section.

\bigskip

\centerline{\S6. Numerical examples I}

\medskip

{\bf 6.1.} In this section, we assume that $F={\mathbf Q}(\sqrt {5})$.
(The formulas (6.1) $\sim$ (6.6) and those given in \S6.5 are valid
for any real quadratic field.) We use the notation of \S5.
The elements $\sigma$, $\nu$ and $\tau$ of $\Gamma ^{\ast}$ satisfy the relations

\begin{equation}
\sigma ^2=1.
\tag{i$^{\prime}$}
\end{equation}
\begin{equation}
(\sigma\tau ) ^3=1.
\tag{ii$^{\prime}$}
\end{equation}
\begin{equation}
(\sigma\nu ) ^2=1.
\tag{iii$^{\prime}$}
\end{equation}
\begin{equation}
\tau\nu\tau\nu ^{-1}=\nu\tau\nu ^{-1}\tau .
\tag{iv$^{\prime}$}
\end{equation}
\begin{equation}
\nu ^2\tau\nu ^{-2}=\tau\nu\tau\nu ^{-1}.
\tag{v$^{\prime}$}
\end{equation}

\proclaim Theorem 6.1. The fundamental relations satisfied by the generators
$\sigma$, $\nu$, $\tau$ of $\Gamma ^{\ast}$ are (i$^{\prime}$) $\sim$ (v$^{\prime}$).
\par

This theorem follows from Theorem A.1. We sketch a proof.
We have $\mu =\nu ^2$, $\eta =\nu\tau\nu ^{-1}$.
Then we can check easily that the relations (i) $\sim$ (vii) in Theorem A.1 follow
from (i$^{\prime}$) $\sim$ (v$^{\prime}$). Suppose that
\begin{equation}
u_1u_2\cdots u_m=1
\tag{$\ast$}
\end{equation}
is a relation. Here $u_i$ is one of $\sigma$, $\nu$, $\nu ^{-1}$, $\tau$, $\tau ^{-1}$.
In ($\ast$), we substitute $\nu ^{-1}$ by $\mu ^{-1}\nu$. Then we obtain a relation
\begin{equation}
v_1v_2\cdots v_n=1.
\tag{$\ast\ast$}
\end{equation}
Here $v_i$ is one of $\sigma$, $\nu$, $\mu ^{-1}$, $\tau$, $\tau ^{-1}$.
The number of $v_i$ such that $v_i=\nu$ is even. If this number is $0$,
then ($\ast\ast$) is the relation among the elements $\sigma$, $\mu$ and $\tau$.
If this number is positive, then in ($\ast\ast$), a term of the form
$\nu X\nu$ is contained, where $X$ is an expression which contains
only $\sigma$, $\tau$ and $\mu$. We may replace $\nu X\nu$ by
$\nu X\nu ^{-1}\mu$. By the relations
$$
\nu\sigma\nu ^{-1}=\sigma\nu ^{-2}=\sigma\mu ^{-1}, \qquad
\nu\tau\nu ^{-1}=\eta ,
$$
$\nu X\nu ^{-1}$ is transformed to an expression which contains
only $\sigma$, $\mu$, $\tau$, $\eta$ and their inverses.
Repeating this procedure, ($\ast\ast$) can be reduced to
a relation among the elements $\sigma$, $\mu$, $\tau$ and $\eta$.
By Theorem A.1, this relation follows from the fundamental relations
(i) $\sim$ (vii). Since (i) $\sim$ (vii) follow from
(i$^{\prime}$) $\sim$ (v$^{\prime}$), our assertion is proved.

\medskip

Let $\mathcal F^{\ast}$ be the free group on three letters
$\widetilde\sigma$, $\widetilde\nu$, $\widetilde\tau$.
We define a surjective homomorphism
$\pi ^{\ast}: \mathcal F^{\ast} \longrightarrow \Gamma ^{\ast}$
by  $\pi ^{\ast}(\widetilde\sigma )=\sigma$, $\pi ^{\ast}(\widetilde\nu )=\nu$,
$\pi ^{\ast}(\widetilde\tau )=\tau$. Let $R^{\ast}$ be the kernel of $\pi ^{\ast}$.
We have $\Gamma ^{\ast}=\mathcal F^{\ast}/R^{\ast}$.
By Theorem 6.1, $R^{\ast}$ is generated by the elements
\begin{equation}
\widetilde\sigma ^2,
\tag{i$^{\ast}$}
\end{equation}
\begin{equation}
(\widetilde\sigma\widetilde\tau ) ^3,
\tag{ii$^{\ast}$}
\end{equation}
\begin{equation}
(\widetilde\sigma\widetilde\nu )^2,
\tag{iii$^{\ast}$}
\end{equation}
\begin{equation}
\widetilde\tau\widetilde\nu\widetilde\tau\widetilde\nu ^{-1}
(\widetilde\nu\widetilde\tau\widetilde\nu ^{-1}\widetilde\tau )^{-1},
\tag{iv$^{\ast}$}
\end{equation}
\begin{equation}
\widetilde\nu ^2\widetilde\tau\widetilde\nu ^{-2}
(\widetilde\tau\widetilde\nu\widetilde\tau\widetilde\nu ^{-1})^{-1}
\tag{v$^{\ast}$}
\end{equation}
and their conjugates.

Let $P^{\ast}$ be the subgroup of $\Gamma ^{\ast}$ consisting of elements
which can be represented by upper triangular matrices.
Let $\mathcal F_{P^{\ast}}$ be the subgroup of $\mathcal F^{\ast}$
generated by $\widetilde \nu$ and $\widetilde\tau$.
Then $\pi ^{\ast}\vert \mathcal F_{P^{\ast}}: 
\mathcal F_{P^{\ast}} \longrightarrow P^{\ast}$ is surjective.
Let $R_{P^{\ast}}$ be the kernel of this homomorphism.
We see that $R_{P^{\ast}}$ is generated by (iv$^{\ast}$) and (v$^{\ast}$)
and their conjugates.

We have $[\mathcal F^{\ast}:(\pi ^{\ast})^{-1}(\Gamma )]=2$.
The following lemma can be proved easily by
applying the method of Reidemeister--Schreier (cf. Schreier [Sc], Suzuki [Su], \S6).

\proclaim Lemma 6.2. The group $(\pi ^{\ast})^{-1}(\Gamma )$ is the free group
on five elements
$\widetilde\sigma$, $\widetilde\nu ^2$, $\widetilde\tau$,
$\widetilde\nu\widetilde\sigma\widetilde\nu ^{-1}$ and
$\widetilde\nu\widetilde\tau\widetilde\nu ^{-1}$.
\par

We put $\widetilde\nu ^2=\widetilde\mu$,
$\widetilde\nu\widetilde\tau\widetilde\nu ^{-1}=\widetilde\eta$.
Let $\mathcal F$ be the free group on four elements
$\widetilde\sigma$, $\widetilde\mu$, $\widetilde\tau$ and $\widetilde\eta$.
Then our notation becomes consistent with that given in the beginning of section 5.
We have $\mathcal FR^{\ast}=(\pi ^{\ast})^{-1}(\Gamma )$.

\medskip

{\bf 6.2.} For every $\gamma \in \Gamma ^{\ast}$, we choose
$\widetilde\gamma \in \mathcal F^{\ast}$ so that
$\pi ^{\ast}(\widetilde\gamma )=\gamma$. For explicit calculations,
it is necessary to specify the choice of $\widetilde\gamma$.
First let $p \in P$. We can write $p=\mu ^a\tau ^b\eta ^c$ and
this expression is unique. We put
$\widetilde p=\widetilde\mu ^a\widetilde\tau ^b\widetilde\eta ^c$.
Next let $p \in P^{\ast}$. We have $p \in P$ or $p=\nu p_1$ with $p_1 \in P$.
In the latter case, we put $\widetilde p=\widetilde\nu\widetilde p_1$.

Let $\Delta$ be a complete set of representatives for $P\backslash \Gamma$
as in \S5.1. Then $\Delta$ is also a complete set of representatives for
$P^{\ast}\backslash \Gamma ^{\ast}$. For $\gamma \in \Gamma ^{\ast}$,
we write $\gamma =p\delta$ with $p \in P^{\ast}$, $\delta \in \Delta$
and put $\widetilde\gamma =\widetilde p\widetilde\delta$.
Our task is to specify the choice of $\Delta$ and define $\widetilde\delta$
for $\delta \in \Delta$. To specify $\Delta$ is equivalent to
choose one element from every coset $P\gamma$, $\gamma \in \Gamma$.
Let $\gamma =\begin{pmatrix} a & b \\ c & d \end{pmatrix}$.

(1) The case where $P\gamma =P$. We take $1$ as the representative.
We take the identity element of $\mathcal F$ as $\widetilde 1$.

(2) The case where $c \in E_F$. We can take an element of the form
$\begin{pmatrix} 0 & -1 \\ 1 & d \end{pmatrix}$ as the representative.
We define
$$
\widetilde {\begin{pmatrix} 0 & -1 \\ 1 & d \end{pmatrix}}
=\widetilde\sigma \widetilde {\begin{pmatrix} 1 & d \\ 0 & 1 \end{pmatrix}}.
$$

(3) 
\footnote{In this paper, this step will be used for the actual calculations
only in the case $a \in E_F$. Since it will become necessary in future
calculations, we write one (tentative) algorithm explicitly.}
The case where $c \neq 0$ and $c \notin E_F$.
We note that $\mathcal O_F$ is a Euclidean ring with respect to
the absolute value of the norm (cf. [HW], Theorem 247, p. 213):
For every $x$, $y \in \mathcal O_F$, $x \neq 0$, there exist
$q$, $r \in \mathcal O_F$ such that
$$
y=qx+r, \qquad \vert N(r)\vert <\vert N(x)\vert .
$$
We have
$$
\begin{pmatrix} u & 0 \\ 0 & u^{-1}\end{pmatrix}
\begin{pmatrix} a & b \\ c & d \end{pmatrix}
=\begin{pmatrix} ua & ub \\ u^{-1}c & u^{-1}d \end{pmatrix}, \quad
\begin{pmatrix} 1 & t \\ 0 & 1 \end{pmatrix}
\begin{pmatrix} a & b \\ c & d \end{pmatrix}
=\begin{pmatrix} a+tc & b+td \\ c & d \end{pmatrix}.
$$
First mulplying $\gamma$ on the left by
$\begin{pmatrix} u & 0 \\ 0 & u^{-1}\end{pmatrix}$, $u \in E_F$,
we normalize $c$ so that
$$
c \gg 0, \qquad 1 \leq c^{\prime}/c < \epsilon ^2.
$$
Next mulplying $\gamma$ on the left by
$\begin{pmatrix} 1 & t \\ 0 & 1 \end{pmatrix}$, $t \in \mathcal O_F$,
we may assume that $\vert N(a)\vert <\vert N(c)\vert$
by the Euclidean algorithm. However to specify the choice of
$t$ is not necessarily easy. In other words, there can be many
choices of such $a$'s. We make the preference order of the choice of
$a$ as follows. Put $a=\alpha +\beta\epsilon$, $\alpha$, $\beta \in {\mathbf Z}$.

1. $\vert\alpha\vert +\vert\beta\vert$ is minimum.
2. $\vert\alpha\vert$ is minimum.
3. $\vert\beta\vert$ is minimum.
4. $\alpha \geq 0$. 5. $\beta \geq 0$.

We define $\widetilde\delta$ for $\delta \in \Delta$ as follows.
We put $\delta =\begin{pmatrix} a & b \\ c & d \end{pmatrix}$
and proceed by induction on $\vert N(c)\vert$.
The case $\vert N(c)\vert =0$ or $1$ is settled by (1) and (2).
By our choice of $\Delta$, we have $\vert N(a)\vert <\vert N(c)\vert$.
Put $\sigma ^{-1}\delta =p_1\delta _1$, $p_1 \in P$, $\delta _1 \in \Delta$,
$\delta _1=\begin{pmatrix} a_1 & b_1 \\ c_1 & d_1 \end{pmatrix}$.
We have $\vert N(c_1)\vert =\vert N(a)\vert <\vert N(c)\vert$.
We define $\widetilde\delta =\widetilde\sigma\widetilde p_1\widetilde\delta _1$.

\medskip

{\bf 6.3.} Let $f \in Z_P^2(\Gamma , V)$ be a normalized parabolic $2$-cocycle.
We first consider $f^+$ (cf. \S5.6). We put $f^{\ast}=\widetilde T(f)$.
Then $f^{\ast} \in Z^2(\Gamma ^{\ast}, V)$ and
$f^{\ast}\vert \Gamma=f^+$ (cf. \S5.3). We can verify easily
the parabolic condition
\begin{equation}
f^{\ast}(p\gamma _1, \gamma _2)=pf^{\ast}(\gamma _1, \gamma _2),
\qquad p \in P^{\ast}, \ \gamma _1, \gamma _2 \in \Gamma ^{\ast}.
\tag{6.1}
\end{equation}
We have
\begin{equation}
H^2(\Gamma ^{\ast}, V) \cong
H^1(R^{\ast}, V)^{\Gamma ^{\ast}}/{\rm Im}(H^1(\mathcal F^{\ast}, V)).
\tag{6.2}
\end{equation}
Let $\varphi \in H^1(R^{\ast}, V)^{\Gamma ^{\ast}}$ be a corresponding
element to $f^{\ast}$. We recall that $\varphi$ is obtained in the following way.
There exists $a \in C^1(\mathcal F^{\ast}, V)$ such that
\begin{equation}
a(g_1g_2) = g_1a(g_2)+a(g_1)-f^{\ast}(\pi ^{\ast}(g_1), \pi ^{\ast}(g_2)),
\qquad g_1, g_2 \in \mathcal F^{\ast}.
\tag{6.3}
\end{equation}
Then $\varphi =a\vert R^{\ast}$. We may regard (6.3) as a rule for
determining the value $a(g)$ according to the length of a word
$g \in \mathcal F^{\ast}$. We can take
$a(\widetilde\sigma )=a(\widetilde\nu )=a(\widetilde\tau )=0$.
Then we have $a\vert \mathcal F_{P^{\ast}}=0$,
since (6.1) yields $f^{\ast}(p, \gamma )=0$, $p \in P^{\ast}$,
$\gamma \in \Gamma ^{\ast}$. In particular, we have
\begin{equation}
\varphi \vert R_{P^{\ast}}=0.
\tag{6.4}
\end{equation}
As shown in \S1.5, we may assume that
\begin{equation}
f^{\ast}(\gamma _1, \gamma _2)=
-\varphi(\widetilde\gamma _1\widetilde\gamma _2(\widetilde{\gamma _1\gamma _2})^{-1})
\tag{6.5}
\end{equation}
adding a coboundary to $f^{\ast}$. By (6.4), we can check that
$f^{\ast}$ satisfies (6.1) in the same way as in \S5.1.
We have (cf. (5.3))
$$
f^{\ast}(\sigma , \mu )=-\varphi ((\widetilde\sigma \widetilde\mu )^2)
+\sigma\mu \varphi (\widetilde\sigma ^2).
$$
We have
\begin{equation*}
\begin{aligned}
& \varphi ((\widetilde\sigma \widetilde\mu )^2)
=\varphi (\widetilde\sigma\widetilde\nu ^2\widetilde\sigma\widetilde\nu ^2)
=\varphi (\widetilde\sigma\widetilde\nu\widetilde\sigma\widetilde\nu
\widetilde\nu ^{-1}\widetilde\sigma ^{-1}\widetilde\nu
\widetilde\sigma\widetilde\nu ^2) \\
= \, & \varphi (\widetilde\sigma\widetilde\nu\widetilde\sigma\widetilde\nu )
+\varphi (\widetilde\nu ^{-1}\widetilde\sigma ^{-2}\widetilde\nu )
+\varphi (\widetilde\nu ^{-1}\widetilde\sigma\widetilde\nu
\widetilde\sigma\widetilde\nu ^2)
= (1+\nu ^{-1})\varphi ((\widetilde\sigma \widetilde\nu )^2)
-\nu ^{-1}\varphi (\widetilde\sigma ^2).
\end{aligned}
\end{equation*}
Therefore we obtain
\begin{equation}
f^{\ast}(\sigma , \mu )
=-(1+\nu ^{-1})\varphi ((\widetilde\sigma \widetilde\nu )^2)
+(\sigma\mu +\nu ^{-1})\varphi (\widetilde\sigma ^2).
\tag{6.6}
\end{equation}

Clearly $\varphi$ is determined by its values on the elements (i$^{\ast}$) $\sim$
(v$^{\ast}$). By (6.4), $\varphi$ takes the value $0$ on the elements
(iv$^{\ast}$) and (v$^{\ast}$). We have
$\sigma\varphi (\widetilde\sigma ^2)=\varphi (\widetilde\sigma ^2)$.
Take $h \in H^1(\mathcal F^{\ast}, V)$ so that 
$h(\widetilde\sigma )=-\varphi (\widetilde\sigma ^2)/2$,
$h(\widetilde\nu )=0$, $h(\widetilde\tau )=0$. Adding
$h\vert R^{\ast}$ to $\varphi$, we may assume that
$\varphi (\widetilde\sigma ^2)=0$; $\varphi$ still satisfies (6.4).

We analyze the process of adding $h\vert R^{\ast}$ to
$\varphi$ in more detail. For $S$, $T$, $U \in V$, we can find
$h \in H^1(\mathcal F^{\ast}, V)$ such that
$$
h(\widetilde\sigma )=S, \qquad h(\widetilde\tau )=T, \qquad
h(\widetilde\nu )=U.
$$
We find easily that the conditions that $h$ vanishes on the elements
(iv$^{\ast}$) and (v$^{\ast}$) are
\begin{equation}
(1+\tau\nu -\nu -\nu\tau\nu ^{-1})T+(\tau -1)(1-\nu\tau\nu ^{-1})U=0,
\tag{6.7}
\end{equation}
\begin{equation}
(\nu ^2-1-\tau\nu )T+(1+\nu -\nu ^2\tau\nu ^{-1}-\tau )U=0
\tag{6.8}
\end{equation}
respectively. We have
\begin{equation}
h(\widetilde\sigma ^2)=(1+\sigma )S.
\tag{6.9}
\end{equation}
We put
$$
A=\varphi ((\widetilde\sigma\widetilde\nu )^2), \qquad
B=\varphi ((\widetilde\sigma\widetilde\tau )^3).
$$
We note that
\begin{equation}
\sigma\nu A=A, \qquad \sigma\tau B=B.
\tag{6.10}
\end{equation}
Our objective is to determine $A$ explicitly.

\medskip

{\bf 6.4.} Let us consider the Hecke operators.
We put $g^{\ast}=T(\varpi )f^{\ast}$ where $g^{\ast}$
is defined by (5.4) with $\Gamma ^{\ast}$ in place of $\Gamma$.
Let $\psi \in H^1(R^{\ast}, V)^{\Gamma ^{\ast}}$ be a corresponding
element to $g^{\ast}$. We see that Proposition 5.1 remains valid with
$\Gamma ^{\ast}$ and $P^{\ast}$ in place of $\Gamma$ and $P$.
In particular we may assume that $\psi$ is given by the formula
\begin{equation*}
\begin{aligned}
& \psi(\widetilde\gamma _1\widetilde\gamma _2\cdots\widetilde\gamma _m) \\
= \, & c\sum _{i=1}^d \beta _i^{-1}
\varphi (\widetilde{\beta _i\gamma _1\beta _{q_1(i)}^{-1}}
\widetilde{\beta _{q_1(i)}\gamma _2\beta _{q_2(i)}^{-1}}\cdots
\widetilde{\beta _{q_{m-1}(i)}\gamma _m\beta _{q_m(i)}^{-1}}
(\widetilde{\beta _i\gamma _1\gamma _2\cdots \gamma _m\beta _{q_m(i)}^{-1}})^{-1}).
\end{aligned}
\tag{6.11}
\end{equation*}
Here $\gamma _j=\sigma$ or $\gamma _j \in P^{\ast}$ and
$\gamma _1\gamma _2\cdots \gamma _m=1$.

\medskip

{\bf Example 6.3.}  Let us consider $T(2)$. We may take
\begin{equation*}
\begin{aligned}
& \beta _1=\begin{pmatrix} 1 & 0 \\ 0 & 2 \end{pmatrix}, \qquad
\beta _2=\begin{pmatrix} 1 & 1 \\ 0 & 2 \end{pmatrix}, \qquad
\beta _3=\begin{pmatrix} 1 & \epsilon \\ 0 & 2 \end{pmatrix}, \qquad \\
& \beta _4=\begin{pmatrix} 1 & \epsilon ^2\\ 0 & 2 \end{pmatrix}, \qquad
\beta _5=\begin{pmatrix} 2 & 0 \\ 0 & 1 \end{pmatrix}.
\end{aligned}
\end{equation*}
By (6.11), we find
$$
\psi ((\widetilde\sigma\widetilde\tau )^3)
=c(\beta _3^{-1}Z_3+\beta _4^{-1}Z_4),
$$
where
\begin{equation}
Z_3=\varphi ((\widetilde
{\begin{pmatrix} \epsilon & -\epsilon ^2 \\ 2 & -\epsilon ^2 \end{pmatrix}}
\widetilde\tau )^3), \qquad
Z_4=\varphi ((\widetilde
{\begin{pmatrix} \epsilon ^2 & -\epsilon ^2 \\ 2 & -\epsilon \end{pmatrix}})^3).
\tag{6.12}
\end{equation}
We have
$$
\widetilde{\begin{pmatrix} \epsilon & -\epsilon ^2 \\ 2 & -\epsilon ^2 \end{pmatrix}}
=\widetilde\sigma
\widetilde{\begin{pmatrix} \epsilon ^{-1}& 0 \\ 0 & \epsilon \end{pmatrix}}
\widetilde{\begin{pmatrix} 1 & \epsilon \\ 0 & 1 \end{pmatrix}}^{-2}
\widetilde\sigma
\widetilde{\begin{pmatrix} 1 & \epsilon \\ 0 & 1 \end{pmatrix}}^{-1}.
$$
Hence, using (6.4), we have
$$
Z_3=\varphi ((\widetilde\sigma
\widetilde{\begin{pmatrix} \epsilon ^{-1}& -2 \\ 0 & \epsilon \end{pmatrix}}
\widetilde\sigma
\widetilde{\begin{pmatrix} 1 & -\epsilon ^{-1}\\ 0 & 1 \end{pmatrix}})^3).
$$
Similarly we obtain
$$
Z_4=\varphi ((\widetilde\sigma
\widetilde{\begin{pmatrix} \epsilon ^{-2}& -2 \\ 0 & \epsilon ^2\end{pmatrix}}
\widetilde\sigma
\widetilde{\begin{pmatrix} 1 & -1 \\ 0 & 1 \end{pmatrix}})^3).
$$

\medskip

{\bf 6.5.} In general, every element $r$ of $R^{\ast}$ can be written as
$$
r=\widetilde\sigma \widetilde p_1\widetilde\sigma \widetilde p_2
\cdots\widetilde\sigma \widetilde p_m
$$
with $p_i \in P^{\ast}$, $1 \leq i \leq m$ such that
$\sigma p_1\sigma p_2\cdots \sigma p_m=1$. We call such an element
an {\it $m$ terms relation}. Theorem 6.1 assures us that
$\varphi (r)$ can be expressed by $A$ and $B$.
The following formulas can be proved easily.
\begin{equation}
\varphi ((\widetilde\sigma\widetilde\nu ^n)^2)
=(1+\nu ^{-1}+\cdots +\nu ^{1-n})A, \qquad n \geq 1,
\tag{6.13a}
\end{equation}
\begin{equation}
\varphi ((\widetilde\sigma\widetilde\nu ^{-n})^2)
=-(\nu+\nu ^2+\cdots +\nu ^n)A, \qquad n \geq 1,
\tag{6.13b}
\end{equation}
For $t \in E_F$, we put
$$
B(t)=\varphi (
\widetilde\sigma\widetilde{\begin{pmatrix} 1 & t \\ 0 & 1 \end{pmatrix}}
\widetilde\sigma\widetilde{\begin{pmatrix} 1 & t^{-1} \\ 0 & 1 \end{pmatrix}}
\widetilde\sigma\widetilde{\begin{pmatrix} 1 & t \\ 0 & 1 \end{pmatrix}}
\widetilde{\begin{pmatrix} t & 0 \\ 0 & t^{-1} \end{pmatrix}} ).
$$
Then we have $B(1)=B$,
\begin{equation}
B(-t)= -\sigma\begin{pmatrix} t & 0 \\ 0 & t^{-1} \end{pmatrix} B(t)
-\begin{pmatrix} t^{-1} & 0 \\ 0 & t \end{pmatrix}
\varphi ((\widetilde\sigma
\widetilde{\begin{pmatrix} t^{-1} & 0 \\ 0 & t \end{pmatrix}})^2),
\tag{6.14}
\end{equation}
\begin{equation*}
\begin{aligned}
B(\epsilon t)= & \nu ^{-1} B(t) \\
+ &\bigg[ 1+\sigma\begin{pmatrix} 1 & \epsilon t \\ 0 & 1 \end{pmatrix}
\sigma\begin{pmatrix} 1 & \epsilon ^{-1}t^{-1} \\ 0 & 1 \end{pmatrix}
-\sigma\begin{pmatrix} 1 & \epsilon t\\ 0 & 1 \end{pmatrix}\sigma \bigg]A,
\end{aligned}
\tag{6.15}
\end{equation*}
\begin{equation}
B(t)= \sigma\begin{pmatrix} 1 & t \\ 0 & 1 \end{pmatrix} B(t^{-1})
+\varphi ((\widetilde\sigma
\widetilde{\begin{pmatrix} t & 0 \\ 0 & t^{-1} \end{pmatrix}})^2).
\tag{6.16}
\end{equation}
By these formulas, we can express $B(t)$ in terms of $A$ and $B$ explicitly.
Using $B(t)$, we have an explicit formula for $\varphi (r)$ for
a three terms relation $r$:
\begin{equation*}
\begin{aligned}
& \varphi (
\widetilde\sigma\widetilde{\begin{pmatrix} u_1 & x_1 \\ 0 & 1 \end{pmatrix}}
\widetilde\sigma\widetilde{\begin{pmatrix} u_2 & x_2 \\ 0 & 1 \end{pmatrix}}
\widetilde\sigma\widetilde{\begin{pmatrix} u_3 & x_3 \\ 0 & 1 \end{pmatrix}}) \\
= & \begin{pmatrix} u_1^{-1} & 0 \\ 0 & 1 \end{pmatrix} B(u_1^{-1}x_1)
+\varphi ((\widetilde\sigma
\widetilde{\begin{pmatrix} u_1 & 0 \\ 0 & 1 \end{pmatrix}})^2) \\
+ & \begin{pmatrix} u_3^{-1} & -u_3^{-1}x_3 \\ 0 & 1 \end{pmatrix}\sigma
\varphi ((\widetilde\sigma
\widetilde{\begin{pmatrix} u_2 & 0 \\ 0 & 1 \end{pmatrix}})^2).
\end{aligned}
\tag{6.17}
\end{equation*}

For an $m$ terms relation $r \in R^{\ast}$, $m \geq 4$, we may write
$p_i=\begin{pmatrix} u_i & x_i \\ 0 & 1 \end{pmatrix}$,
$u_i \in E_F$, $x_i \in \mathcal O_F$, $ 1 \leq i \leq m$.
We see that $\varphi (r)$ reduces to an $(m-2)$ terms relation if
$x_i=0$ for some $i$. If $x_i \in E_F$ for some $i$,
$\varphi (r)$ reduces to an $(m-1)$ terms relation. For example,
if $x_1 \in E_F$ and $m\geq 4$, we have

\begin{equation*}
\begin{aligned}
& \varphi (
\widetilde\sigma\widetilde{\begin{pmatrix} u_1 & x_1 \\ 0 & 1 \end{pmatrix}}
\widetilde\sigma\widetilde{\begin{pmatrix} u_2 & x_2 \\ 0 & 1 \end{pmatrix}}
\widetilde\sigma\widetilde{\begin{pmatrix} u_3 & x_3 \\ 0 & 1 \end{pmatrix}}
\widetilde\sigma\cdots 
\widetilde\sigma\widetilde{\begin{pmatrix} u_m & x_m \\ 0 & 1 \end{pmatrix}} )\\
= \, & \begin{pmatrix} u_1^{-1}u^{-1} & -u_1^{-1} \\ 0 & u \end{pmatrix}
 \varphi (
\widetilde\sigma
\widetilde{\begin{pmatrix} 1 & -u^{-1} \\ 0 & 1 \end{pmatrix}}
\widetilde{\begin{pmatrix} u_2 & x_2 \\ 0 & 1 \end{pmatrix}}
\widetilde\sigma\widetilde{\begin{pmatrix} u_3 & x_3 \\ 0 & 1 \end{pmatrix}}
\widetilde\sigma \\
& \cdots 
\widetilde\sigma\widetilde{\begin{pmatrix} u_m & x_m \\ 0 & 1 \end{pmatrix}} 
\widetilde{\begin{pmatrix} u_1^{-1}u^{-1} & -u_1^{-1} \\ 0 & u \end{pmatrix}})\\
+ & \begin{pmatrix} u_1^{-1} & 0 \\ 0 & 1 \end{pmatrix} B(u)
+\varphi ((\widetilde\sigma
\widetilde{\begin{pmatrix} u_1 & 0 \\ 0 & 1 \end{pmatrix}})^2).
\end{aligned}
\tag{6.18}
\end{equation*}
Here $u=u_1^{-1}x_1$. For a general $m$ terms relation $r$,
the explicit reduction of $\varphi (r)$ to $A$ and $B$ is a highly non-trivial problem.
The author has an idea on a heuristic algorithm to solve this problem,
but it will not be discussed in this paper.
For our present purposes, the formulas (6.13a) $\sim$ (6.18) are sufficient.

\medskip

{\bf 6.6.} For actual computations, it is convenient to use the decomposition (5.13).
Proposition 5.5 shows that we will lose little information by
assuming $\varphi \in H^1(R^{\ast}, V)^{\Gamma ^{\ast},+}$, so we do assume this.
Then we have
\begin{equation*}
\begin{aligned}
& -\varphi ((\widetilde\sigma\widetilde\tau )^3)
=\varphi (\widetilde\tau ^{-1}\widetilde\sigma\widetilde\tau ^{-1}\widetilde\sigma
\widetilde\tau ^{-1}\widetilde\sigma )
=\tau ^{-1}\varphi ((\widetilde\sigma\widetilde\tau ^{-1})^3) \\
= \, & \tau ^{-1}\varphi (((\widetilde\sigma\widetilde\tau )^3)_{\delta})
=\tau ^{-1}\delta\varphi ((\widetilde\sigma\widetilde\tau )^3).
\end{aligned}
\end{equation*}
Hence
$$
(\delta\tau +1)B=0.
$$
Similarly we obtain
$$
(\delta -1)A=0.
$$

Now we are ready to state explicit numerical examples.
First by numerical computations, we have verified:

\medskip

{\bf Fact 1.} Suppose $0 \leq l_2\leq l_1 \leq 20$.
Then adding $h \vert R^{\ast}$, $h \in H^1(\mathcal F^{\ast}, V)$
to $\varphi$ (keeping $\varphi$ in the plus space under the action of $\delta$),
we may assume $B=0$.

\medskip

Therefore our task is to find constraints on
$A=\varphi ((\widetilde\sigma\widetilde\nu )^2)$. Note that $(\sigma\nu -1)A=0$.
We put
$x=\begin{pmatrix} \epsilon & -\epsilon ^2 \\ 2 & -\epsilon ^2 \end{pmatrix}\tau$ and
\begin{equation}
Z_A^+=\{ \mathbf v \in V \mid (\sigma\nu -1)\mathbf v=0, \ (\delta -1)\mathbf v=0, \
xZ_3=Z_3 \} .
\tag{6.19}
\end{equation}
Here some explanation is called for on the meaning of $xZ_3=Z_3$.
First note that $Z_3$ is defined by (6.12); clearly we must have $xZ_3=Z_3$.
Using the formulas (6.13a) $\sim$ (6.18), we see that
$Z_3$ can be expressed by $A$. Therefore $xZ_3=Z_3$ gives
a constraint on $A$. We define a linear mapping
\begin{equation}
\zeta ^+: Z_A^+ \longrightarrow {\mathbf C}^{l_2+1}
\tag{6.20}
\end{equation}
as follows. Let $\mathbf v \in Z_A^+$. We let the coefficient of
${\mathbf e} _{l_1+2-m}\otimes {\mathbf e}^{\prime}_{(l_1+l_2)/2+2-m}$
in $(1+\nu ^{-1})\mathbf v$ be equal to the $(l_1+l_2)/2+2-m$-th coefficient of
$\zeta ^+(\mathbf v)$, for $(l_1-l_2)/2+1 \leq m \leq (l_1+l_2)/2+1$ (cf. (6.6)).

Suppose that $\varphi$ as above corresponds to a (nonzero) Hecke eigenform
$\Omega \in S_{l_1+2, l_2+2}(\Gamma )$. Suppose that
$l_1$ and $l_2$ are in the range of Fact 1. Then $\zeta ^+(A) \neq 0$
if $l_2 \geq 4$ in the case $l_1\neq l_2$, if $l_2 \geq 6$ in the case
$l_1=l_2$ by Proposition 5.5.

\medskip

{\bf Example 6.4.} We take $l_1=8$, $l_2=4$. Then $\dim S_{10,6}(\Gamma )=1$.
We find $\zeta ^+(Z_A^+)$ is one dimensional and consists of scalar multiples of
${}^t (4, 0, 1, 0, 4)$. Hence we obtain
$$
R(7, \Omega )/R(5, \Omega )=4, \qquad \Omega \in S_{10,6}(\Gamma ).
$$
My computer calculates this example in six seconds.

\medskip

{\bf Example 6.5.} In the same way as in Example 6.4, we obtain the following
numerical values.
$$
R(9, \Omega )/R(7, \Omega )=6, \qquad \Omega \in S_{14,6}(\Gamma ).
$$
$$
R(6, \Omega )/R(4, \Omega )=\frac{25}{6}, \qquad \Omega \in S_{8,8}(\Gamma ).
$$
$$
R(8, \Omega )/R(6, \Omega )=7, \qquad \Omega \in S_{12,8}(\Gamma ).
$$
$$
R(10, \Omega )/R(8, \Omega )=\frac{720}{11}, \qquad \Omega \in S_{12,10}(\Gamma ).
$$
The spaces of cusp forms appearing in this example are all one dimensional.

\medskip

{\bf 6.7.} To deal with the case where $\dim S_{l_1+2, l_2+2}(\Gamma )>1$,
it is necessary to use the action of Hecke operators.
To this end, we consider the contribution of $H^1(\mathcal F^{\ast}, V)$ to $Z_A^+$.
Take $h \in H^1(\mathcal F^{\ast}, V)$ and put
$$
h(\widetilde\sigma )=S, \qquad h(\widetilde\nu )=U, \qquad
h(\widetilde\tau )=T.
$$
We require that $h\vert R^{\ast}$ vanishes on the elements
(i$^{\ast}$), (ii$^{\ast}$), (iv$^{\ast}$), (v$^{\ast}$).
These conditions are equivalent to
\begin{equation}
(\sigma +1)S=0,
\tag{6.21}
\end{equation}
\begin{equation}
\{ (\sigma\tau )^2+\sigma\tau +1\} (\sigma T+S)=0
\tag{6.22}
\end{equation}
and (6.7), (6.8). We have
$$
h((\widetilde\sigma\widetilde\nu )^2)=(\sigma\nu +1)(\sigma U+S).
$$
We also require that
\begin{equation}
(\delta -1)(\sigma\nu +1)(\sigma U+S)=0.
\tag{6.23}
\end{equation}
Let $B_A^+$ be the subspace of $V$ generated by
$(\sigma\nu +1)(\sigma U+S)$ when $S$, $T$, $U$ extend over vectors of $V$
satifying the relations (6.7), (6.8), (6.21), (6.22) and (6.23). 
We have $B_A^+ \subset Z_A^+$. As shown in \S4.1, we have
\begin{equation}
\zeta ^+(B_A^+)=\{ 0\} \quad {\rm if} \ l_1 \neq l_2, \qquad
\dim \zeta ^+(B_A^+) \leq 1\quad {\rm if} \ l_1=l_2.
\tag{6.24}
\end{equation}
By Proposition 5.5, we have
$$
\dim Z_A^+/B_A^+ \geq \dim S_{l_1+2.l_2+2}(\Gamma ) \qquad
{\rm if} \ l_2 \geq 4, \  l_1 \neq l_2 \ {\rm or \ if} \ \  l_1=l_2,  \ l_2 \geq 6.
$$
Now by numerical computations, we have verified:

\medskip

{\bf Fact 2.} Suppose $0 \leq l_2\leq l_1 \leq 20$.
Then $\dim S_{l_1+2, l_2+1}(\Gamma )=\dim Z_A^+/B_A^+$.

\medskip

This fact means that the constraints posed on
$A=\varphi ((\widetilde\sigma\widetilde\nu )^2)$ is enough.

\medskip

{\bf Example 6.6.} We take $l_1=12$, $l_2=8$.
We have $\dim S_{14,10}(\Gamma )=2$.
Moreover we have $\zeta ^+(Z_A^+)=2$ in this case. Hence
$\zeta ^+$ gives an isomorphism of $Z_A^+/B_A^+$ into ${\mathbf C}^{l_2+1}$.
Calculating the action of $T(2)$ on $Z_A^+/B_A^+$ using (6.11),
we find that the eigenvalues are $-2560\pm 960\sqrt{106}$.
Take an eigenvector in $Z_A^+/B_A^+$ and map it by $\zeta ^+$.
Then we find
$$
R(11, \Omega )/R(7, \Omega )=1616-76\sqrt{106}, \qquad
R(9, \Omega )/R(7, \Omega )=\frac {58}{3}-\frac {5}{6}\sqrt{106}
$$
if $0 \neq \Omega \in S_{14,10}(\Gamma )$ satisfies
$\Omega\vert T(2)=(-2560+960\sqrt{106})\Omega$.
If $0 \neq \Omega \in S_{14,10}(\Gamma )$ satisfies
$\Omega\vert T(2)=(-2560-960\sqrt{106})\Omega$, then we have
$$
R(11, \Omega )/R(7, \Omega )=1616+76\sqrt{106}, \qquad
R(9, \Omega )/R(7, \Omega )=\frac {58}{3}+\frac {5}{6}\sqrt{106}.
$$

{\bf Remark 6.7.} The relation $\dim \zeta ^+(Z_A^+)=\dim S_{l_1+2,l_2+2}(\Gamma )$
is rather accidental in the above example. It holds in many cases but we have
$\dim S_{l_1+2,l_2+2}(\Gamma )>\dim \zeta ^+(Z_A^+)$ in general.
Even in the general case, we can obtain ratios of $L$-values
by finding an eigenvector of Hecke operators in $Z_A^+/B_A^+$
and mapping it by $\zeta ^+$.

\medskip

{\bf 6.8.} We next consider the $2$-cocycle $f^-$ (cf. \S 5.6).
The technique of calculation is basically same as for $f^+$, but this case is
somewhat more complicated. Put
\begin{equation}
H^1(R^{\ast}, V)^{\Gamma}
=\{ \varphi \in {\rm Hom}(R^{\ast}, V) \mid
\varphi (grg^{-1})=g\varphi (r), \ g \in \mathcal F, r \in R^{\ast} \} .
\tag{6.25}
\end{equation}
Let $\varphi \in H^1(R^{\ast}, V)^{\Gamma}$. We put
$$
(e\varphi )(r)=\nu ^{-1}\varphi (\widetilde\nu r\widetilde\nu ^{-1}),
\qquad r \in R^{\ast}.
$$
Then we can verify easily that
$$
e\varphi \in H^1(R^{\ast}, V)^{\Gamma}, \qquad e^2\varphi =\varphi .
$$
Therefore $H^1(R^{\ast}, V)^{\Gamma}$ decomposes as
\begin{equation}
H^1(R^{\ast}, V)^{\Gamma}
=H^1(R^{\ast}, V)^{\Gamma, +} \oplus H^1(R^{\ast}, V)^{\Gamma ,-},
\tag{6.26}
\end{equation}
where, for $\epsilon =\pm 1$,
$$
H^1(R^{\ast}, V)^{\Gamma ,\epsilon}
=\{ \varphi \in {\rm Hom}(R^{\ast}, V)^{\Gamma} \mid
\varphi (\widetilde\nu r\widetilde\nu ^{-1})=\epsilon\nu\varphi (r), r \in R^{\ast} \} .
$$

First we take an arbitrary normalized $2$-cocycle $f \in Z^2(\Gamma , V)$.
Since $\mathcal FR^{\ast}$ is a free group, there exists
$a \in C^1(\mathcal FR^{\ast}, V)$ such that
\begin{equation}
f(\pi ^{\ast}(g_1), \pi ^{\ast}(g_2))
=g_1a(g_2)+a(g_1)-a(g_1g_2), \qquad g_1, g_2 \in \mathcal FR^{\ast}.
\tag{6.27}
\end{equation}
As shown in \S1.4, we have
$$
a(gr)=ga(r)+a(g), \qquad
a(grg^{-1})=ga(r), \qquad g \in \mathcal FR^{\ast}, \ r \in R^{\ast}.
$$
Put $\varphi =a\vert R^{\ast}$. Then the above formulas imply
$\varphi \in H^1(R^{\ast}, V)^{\Gamma}$. From the isomorphism
$\Gamma \cong \mathcal FR^{\ast}/R^{\ast} (\cong \mathcal F/\mathcal F\cap R^{\ast}
=\mathcal F/R)$, we obtain
\begin{equation}
H^2(\Gamma , V) \cong H^1(R^{\ast}, V)^{\Gamma}
/{\rm Im}(H^1(\mathcal FR^{\ast}, V))
\tag{6.28}
\end{equation}
and the procedure $f \mapsto \varphi$ described above gives an explicit form
of the isomorphism (6.28). We consider the decomposition of $H^2(\Gamma , V)$
under the action of $\nu$ (cf. the formula below (5.10)). Then we have
\begin{equation}
H^2(\Gamma , V)^{\pm} \cong H^1(R^{\ast}, V)^{\Gamma , \pm}
/({\rm Im}(H^1(\mathcal FR^{\ast}, V))\cap H^1(R^{\ast}, V)^{\Gamma , \pm})).
\tag{6.29}
\end{equation}

\medskip

{\bf 6.9.} Now we consider the $2$-cocycle $f^-$. Let
$\varphi \in H^1(R^{\ast}, V)^{\Gamma , -}$ be a corresponding element.
As for $f^+$, we may assume that
\begin{equation}
\varphi \vert R_{P^{\ast}}=0,
\tag{6.30}
\end{equation}
\begin{equation}
f^-(\gamma _1, \gamma _2)=
-\varphi(\widetilde\gamma _1\widetilde\gamma _2(\widetilde{\gamma _1\gamma _2})^{-1})
\tag{6.31}
\end{equation}
adding a coboundary to $f^-$. We put
$$
A=\varphi ((\widetilde\sigma\widetilde\nu )^2), \qquad
B=\varphi ((\widetilde\sigma\widetilde\tau )^3).
$$
The formulas (6.13a) $\sim$ (6.18) hold with
the following modifications.
\begin{equation}
\varphi ((\widetilde\sigma\widetilde\nu ^n)^2)
=(1-\nu ^{-1}+\nu ^{-2}+\cdots +(-1)^{1-n}\nu ^{1-n})A, \qquad n \geq 1,
\tag{6.13a$^-$}
\end{equation}
\begin{equation}
\varphi ((\widetilde\sigma\widetilde\nu ^{-n})^2)
=(\nu-\nu ^2+\nu ^3-\cdots +(-1)^{1-n}\nu ^n)A, \qquad n \geq 1.
\tag{6.13b$^-$}
\end{equation}
We define $B(t)$, $t \in E_F$ by the same formula as before.
In (6.15), the term $\nu ^{-1}B(t)$ should be replaced by
$-\nu ^{-1}B(t)$; (6.14) and (6.16) hold without any change.
For $u =\pm \epsilon ^n\in E_F$, we define
$\epsilon _0(u)=(-1)^n$. On the right-hand side of (6.17),
the first term should be multiplied by $\epsilon _0(u_1)$
and the third term should be multiplied by $\epsilon _0(u_3)$.
On the right-hand side of (6.18), both of the first and the second term
should be multiplied by $\epsilon _0(u_1)$.

We may and do assume that $f^-$ belongs to the plus subspace of
$H^2(\Gamma , V)^-$ under the action of $\delta$. Then we have
$$
(\delta -1)A=0, \qquad (\delta\tau +1)B=0.
$$
By numerical computations, we have verified

\medskip

{\bf Fact 3.} Suppose $0 \leq l_2\leq l_1 \leq 20$.
Then adding $h \vert R^{\ast}$ for $h \in H^1(\mathcal FR^{\ast}, V)$
such that $h \vert R^{\ast} \in H^1(R^{\ast}, V)^{\Gamma , -}$
to $\varphi$ (keeping $\varphi$ in the plus space under the action of $\delta$),
we may assume $B=0$.

\medskip

Therefore our task is to find constraints on
$A=\varphi ((\widetilde\sigma\widetilde\nu )^2)$.
Note that $(\sigma\nu +1)A=0$. We put
$x=\begin{pmatrix} \epsilon & -\epsilon ^2 \\ 2 & -\epsilon ^2 \end{pmatrix}\tau$
and
\begin{equation}
Z_A^-=\{ \mathbf v \in V \mid (\sigma\nu +1)\mathbf v=0, \ (\delta -1)\mathbf v=0, \
xZ_3=Z_3 \} .
\tag{6.32}
\end{equation}
Here the meaning of the constraint $xZ_3=Z_3$ is the same as for $Z_A^+$.
We define a linear mapping
$$
\zeta ^-: Z_A^- \longrightarrow {\mathbf C}^{l_2+1}
$$
as follows. Let $\mathbf v \in Z_A^-$. We let the coefficient of
${\mathbf e} _{l_1+2-m}\otimes {\mathbf e}^{\prime}_{(l_1+l_2)/2+2-m}$
in $(1-\nu ^{-1})\mathbf v$ be equal to the $(l_1+l_2)/2+2-m$-th coefficient of
$\zeta ^-(\mathbf v)$, for $(l_1-l_2)/2+1 \leq m \leq (l_1+l_2)/2+1$ (cf. (6.6)).

\medskip

{\bf Example 6.8.} We take $l_1=8$, $l_2=6$. Then $\dim S_{10,8}(\Gamma )=1$.
We find $\zeta ^-(Z_A^-)$ is one dimensional and consists of scalar multiples of
\newline ${}^t (2, 0, 7/90, 0, -7/90,0,-2)$. Hence we obtain
$$
R(8, \Omega )/R(6, \Omega )=\frac{180}{7}, \qquad \Omega \in S_{10,8}(\Gamma ).
$$

\medskip

{\bf Example 6.9.} In the same way as in Example 6.8, we obtain the following
numerical values.
$$
R(9, \Omega )/R(7, \Omega )=\frac{70}{3}, \qquad \Omega \in S_{12,8}(\Gamma ).
$$
$$
R(9, \Omega )/R(7, \Omega )=42, \qquad \Omega \in S_{12,10}(\Gamma ).
$$
The spaces of cusp forms appearing in this example are all one dimensional.

\medskip

{\bf 6.10.} To treat the case where $\dim S_{l_1+2, l_2+2}(\Gamma )>1$,
it is necessary to consider Hecke operators.

First let us write down
${\rm Im}(H^1(\mathcal FR^{\ast}, V))\cap H^1(R^{\ast}, V)^{\Gamma , \pm}$
which appears on the right-hand side of (6.29), explicitly. Take
$h \in Z^1(\mathcal FR^{\ast}, V)$. We put
\begin{equation}
(e_0h)(x)=\nu ^{-1}h(\widetilde\nu x\widetilde\nu ^{-1}),
\qquad x \in \mathcal FR^{\ast}.
\tag{6.33}
\end{equation}
We can check easily that $e_0h \in Z^1(\mathcal FR^{\ast}, V)$ and that
$$
(e_0^2h)(x)=h(x)+(\nu ^{-2}-x\nu ^{-2})h(\widetilde\nu ^2),
\qquad x \in \mathcal FR^{\ast}.
$$
If we restrict $h$ to $R^{\ast}$, then the action $e_0$ coincides with the action
of $e$ defined in \S6.8. We have $(e_0^2h)\vert R^{\ast}=h\vert R^{\ast}$.
We put
$$
h^{\pm}=h\pm e_0h.
$$
A general element of
${\rm Im}(H^1(\mathcal FR^{\ast}, V))\cap H^1(R^{\ast}, V)^{\Gamma , \pm}$
can be obtained as $h^{\pm}\vert R^{\ast}$ from a general element
$h \in Z^1(\mathcal FR^{\ast}, V)$.

Let $Z^1(\mathcal FR^{\ast}, V)^{\pm}$ be the subgroup of
$Z^1(\mathcal FR^{\ast}, V)$ consisting of all elements whose restrictions
to $R^{\ast}$ belong to $H^1(R^{\ast}, V)^{\Gamma , \pm}$.
Take $\epsilon _1=\pm 1$ and put $h^{\pm}=h+\epsilon _1e_0h$.
For the free generators $\widetilde\sigma$, $\widetilde\tau$, $\widetilde\nu ^2$,
$\widetilde\nu\widetilde\sigma\widetilde\nu ^{-1}$,
$\widetilde\nu\widetilde\tau\widetilde\nu ^{-1}$ of $\mathcal FR^{\ast}$, we put
$$
h(\widetilde\sigma )=S_1, \quad h(\widetilde\tau )=T_1, \quad
h(\widetilde\nu ^2)=U, \quad h(\widetilde\nu\widetilde\sigma\widetilde\nu ^{-1})=V_1,
\quad h(\widetilde\nu\widetilde\tau\widetilde\nu ^{-1})=W_1.
$$
Then we find
$$
h^{\pm}(\widetilde\sigma )=S_1+\epsilon _1\nu ^{-1}V_1,
$$
$$
h^{\pm}(\widetilde\tau )=T_1+\epsilon _1\nu ^{-1}W_1,
$$
$$
h^{\pm}(\widetilde\nu ^2 )=(1+\epsilon _1\nu ^{-1})U,
$$
$$
h^{\pm}(\widetilde\nu\widetilde\sigma\widetilde\nu ^{-1})
=V_1+\epsilon _1\nu S_1+\epsilon _1\nu ^{-1}(1-\nu ^4\sigma )U,
$$
$$
h^{\pm}(\widetilde\nu\widetilde\tau\widetilde\nu ^{-1})
=W_1+\epsilon _1\nu T_1+\epsilon _1(\nu ^{-1}-\nu\tau\nu ^{-2} )U.
$$
Fix $\epsilon _1=\pm 1$ and put
\begin{equation}
h^{\pm}(\widetilde\sigma )=S, \quad h^{\pm}(\widetilde\tau )=T.
\tag{6.34}
\end{equation}
Then $V_1$ and $W_1$ are eliminated and we obtain
\begin{equation}
h^{\pm}(\widetilde\nu ^2 )=(1+\epsilon _1\nu ^{-1})U,
\tag{6.35}
\end{equation}
\begin{equation}
h^{\pm}(\widetilde\nu\widetilde\sigma\widetilde\nu ^{-1})
=\epsilon _1\nu S+\epsilon _1\nu ^{-1}(1-\nu ^4\sigma )U,
\tag{6.36}
\end{equation}
\begin{equation}
h^{\pm}(\widetilde\nu\widetilde\tau\widetilde\nu ^{-1})
=\epsilon _1\nu T+\epsilon _1(\nu ^{-1}-\nu\tau\nu ^{-2} )U.
\tag{6.37}
\end{equation}
Clearly $S$, $T$ and $U$ can take arbirary three vectors of $V$.
The formulas (6.34) $\sim$ (6.37) describe a general element of
$Z^1(\mathcal FR^{\ast} , V)^{\pm}$.
The conditions for $h^{\pm}$ to vanish on the elements (iv$^{\ast}$)
and (v$^{\ast}$) are
\begin{equation}
\{ \nu\tau\nu ^{-1}-1+\epsilon _1(1-\tau )\nu \}T
+\epsilon _1(1-\tau )(\nu ^{-1}-\nu ^{-1}\tau\nu ^{-2})U=0,
\tag{6.38}
\end{equation}
\begin{equation}
(1+\epsilon _1\tau\nu -\nu ^2)T
+\{ \epsilon _1\tau (\nu ^{-1}-\nu\tau\nu ^{-2})-(1-\nu ^2\tau\nu ^{-2})
(1+\epsilon _1\nu ^{-1}) \} U=0
\tag{6.39}
\end{equation}
respectively. For $h^{\pm} \in Z^1(\mathcal FR^{\ast}, V)^{\pm}$ as above,
we have
\begin{equation}
h^{\pm}((\widetilde\sigma\widetilde\nu )^2)
=(1+\epsilon _1\sigma\nu )S+(\nu ^{-2}+\epsilon _1\sigma\nu ^{-1})U.
\tag{6.40}
\end{equation}
Now we consider the case $\epsilon _1=-1$.
Let $B_A^-$ be the subspace of $V$ generated by
$(1-\sigma\nu )S+(\nu ^{-2}-\sigma\nu ^{-1})U$
when $S$, $T$, $U$ extend over vectors of $V$ satisfying the relations
(6.21), (6.22), (6.38), (6.39) and
\begin{equation}
(\delta -1)\{ (1-\sigma\nu )S+(\nu ^{-2}-\sigma\nu ^{-1})U\} =0.
\tag{6.41}
\end{equation}

We have $B_A^- \subset Z_A^-$. As shown in \S4.1, we have
$$
\zeta ^-(B_A^-)=\{ 0\} \quad {\rm if} \ l_1 \neq l_2, \qquad
\dim \zeta ^-(B_A^-) \leq 1\quad {\rm if} \ l_1=l_2.
$$
Using Proposition 5.2, (2), we can show that
$$
\dim Z_A^-/B_A^- \geq \dim S_{l_1+2.l_2+2}(\Gamma ) \qquad
{\rm if} \ l_2 \geq 2, \  l_1 \neq l_2 \ {\rm or \ if} \ \  l_1=l_2,  \ l_2 \geq 4.
$$
Now by numerical computations, we have verified:

\medskip

{\bf Fact 4.} Suppose $0 \leq l_2\leq l_1 \leq 20$.
Then $\dim S_{l_1+2, l_2+1}(\Gamma )=\dim Z_A^-/B_A^-$.

\medskip

The formula (6.11) can be generalized in the following way.
We put $g^-=T(\varpi )f^-$ where $g^-$ is defined by (5.4).
Let $\varphi \in H^1(R^{\ast}, V)^{\Gamma , -}$ be a corresponding
element to $f^-$. We may assume that (6.31) holds.
There exists a $1$-cochain $b \in C^1(\mathcal FR^{\ast}, V)$ such that
\begin{equation}
f^-(\pi ^{\ast}(x_1), \pi ^{\ast}(x_2))
=x_1b(x_2)+b(x_1)-b(x_1x_2), \qquad x_1, x_2 \in \mathcal FR^{\ast}.
\tag{6.42}
\end{equation}
As the intial conditions, we may assume that
$$
b(\widetilde\sigma )=0, \quad b(\widetilde\nu ^2)=0, \quad
b(\widetilde\tau )=0, \quad
b(\widetilde\nu\widetilde\sigma\widetilde\nu ^{-1})=0, \quad
b(\widetilde\nu\widetilde\tau\widetilde\nu ^{-1})=0
$$ 
for the free generators of $\mathcal FR^{\ast}$. Then the formula
(5.9) holds when $b(\widetilde\gamma _j)=0$, $1 \leq j \leq m$.
This condition holds if $\widetilde\gamma _j$ is equal to
one of the five free generators as above or their inverses.
In particular, $\psi =b\vert R^{\ast}$ is given by
\begin{equation*}
\begin{aligned}
& \psi(\widetilde\gamma _1\widetilde\gamma _2\cdots\widetilde\gamma _m) \\
= \, & c\sum _{i=1}^d \beta _i^{-1}
\varphi (\widetilde{\beta _i\gamma _1\beta _{q_1(i)}^{-1}}
\widetilde{\beta _{q_1(i)}\gamma _2\beta _{q_2(i)}^{-1}}\cdots
\widetilde{\beta _{q_{m-1}(i)}\gamma _m\beta _{q_m(i)}^{-1}}
(\widetilde{\beta _i\gamma _1\gamma _2\cdots \gamma _m\beta _{q_m(i)}^{-1}})^{-1})
\end{aligned}
\end{equation*}
provided $\widetilde\gamma _j$ is equal to
one of the five free generators of $\mathcal FR^{\ast}$ or their inverses and
$\gamma _1\gamma _2\cdots\gamma _m=1$.
The above formula is the same as (6.11) but there is one important point
about which we must be careful. This $\psi$ belongs to
$H^1(R^{\ast}, V)^{\Gamma}$ and gives a corresponding element to $g^-$
but it does not necessarily belong to $H^1(R^{\ast}, V)^{\Gamma ,-}$.
We obtain $\psi ^- \in H^1(R^{\ast}, V)^{\Gamma ,-}$
corresponding to $g^-$ by $\psi ^-=\frac {1}{2}(1-e)\psi$ (cf. \S6.8).

\medskip

{\bf Example 6.10.} We take $l_1=12$, $l_2=8$.
We have $\dim S_{14,10}(\Gamma )=2$.
Moreover we have $\zeta ^-(Z_A^-)=2$ in this case. Hence
$\zeta ^-$ gives an isomorphism of $Z_A^-/B_A^-$ into ${\mathbf C}^{l_2+1}$.
Take an eigenvector of $T(2)$ in $Z_A^-/B_A^-$ and map it by $\zeta ^-$.
Then we find
$$
R(10, \Omega )/R(8, \Omega )=50-\sqrt{106}, \qquad
$$
if $0 \neq \Omega \in S_{14,10}(\Gamma )$ satisfies
$\Omega\vert T(2)=(-2560+960\sqrt{106})\Omega$.
If $0 \neq \Omega \in S_{14,10}(\Gamma )$ satisfies
$\Omega\vert T(2)=(-2560-960\sqrt{106})\Omega$, then we have
$$
R(10, \Omega )/R(8, \Omega )=50+\sqrt{106}.
$$
Let $\Omega$ be a Hecke eigenform of $S_{14,10}(\Gamma )$.
Then $L(m, \Omega )$ is a critical value for integers in the range
$3 \leq m \leq 11$ (cf. (5.14)). We have $L(s, \Omega )=L(14-s, \Omega )$ (cf. (2.7)).
By Examples 6.6 and 6.10, we have treated all critical values
on the right of the critical line.

\medskip

{\bf Example 6.11.} We take $l_1=l_2=18$.
We have $\dim S_{20,20}(\Gamma )=7$.
Calculating the action of $T(2)$ on $Z_A^+/B_A^+$ using (6.11),
we find that the characteristic polynomial of $T(2)$ is
(we can use $Z_A^-/B_A^-$ which gives the same result)
\begin{equation*}
\begin{aligned}
& (X-97280)^2(X+840640)
(X^4-1286780X^3+19006483200X^2 \\
& +27181090390835200X-22979876427231395840000).
\end{aligned}
\end{equation*}
The irreducible factor of degree four corresponds to the base change part from
$S_{20}(\Gamma _0(5), (\frac {}{5}))$;
$X+840640$ corresponds to the base change part from $S_{20}({\rm SL}_2(\mathbf Z))$;
the factor $(X-97280)^2$ corresponds to the non base change part.
Let $\Omega \in \dim S_{20,20}(\Gamma )$ be a Hecke eigenform
in the non base change part. A calculation for the plus part yields the result
$$
R(18, \Omega )/R(10, \Omega )=39355680000, \qquad
R(16, \Omega )/R(10, \Omega )=33163650,
$$
$$
R(14, \Omega )/R(10, \Omega )=\frac{1266460}{27}, \qquad
R(12, \Omega )/R(10, \Omega )=\frac {26075}{216}.
$$
A calculation for the minus part yields the result
$$
R(17, \Omega )/R(11, \Omega )=\frac {111006792000}{803}, \qquad
R(15, \Omega )/R(11, \Omega )=\frac {54618434}{365},
$$
$$
R(13, \Omega )/R(11, \Omega )=\frac {453159}{1606}.
$$
We note that though there are two Hecke eigenforms in the non base change part,
these ratios are the same for them.
\footnote{We can show that the $L$-functions (2.47) are the same for two
Hecke eigenforms in the non base change part. In fact, let
$\Omega \neq 0$ be a Hecke eigenform in the non base change part and
let $\lambda (\mathfrak m)$ be the eigenvalue of $T(\mathfrak m)$ for
$\Omega$. For the nontrivial automorphism $\sigma$ of $F$,
there exists a Hecke eigenform $\Omega _{\sigma} \neq 0$ such that
$\Omega _{\sigma} \vert T(\mathfrak m)=\lambda (\mathfrak m^{\sigma})\Omega _{\sigma}$
(cf. [Y2], p. 1035, Remark).
Since $\Omega$ is not a base change, we have
$\lambda (\mathfrak m) \neq \lambda (\mathfrak m^{\sigma})$
for some $\mathfrak m$. Hence $\Omega _{\sigma}$ is not a constant
multiple of $\Omega$. On the other hand, $L(s, \Omega _{\sigma})$
is equal to $L(s, \Omega )$.
}

\bigskip

\centerline{\S7. Numerical examples II}

\medskip

{\bf 7.1.} In this section, we treat the case $F={\mathbf Q}(\sqrt {13})$.
We use the same notation as in the previous section.
Many results there remain valid in the present case so we will be brief.

The fundamental unit of $F$ is $\epsilon =\frac {3+\sqrt {13}}{2}$.
The elements $\sigma$, $\nu$ and $\tau$ of $\Gamma ^{\ast}$
satisfy the relations (i$^{\prime}$) $\sim$ (iv$^{\prime}$) in \S6.1 and
\begin{equation}
\nu ^2\tau\nu ^{-2}=\tau (\nu\tau\nu ^{-1})^3.
\tag{v$^{\prime}$}
\end{equation}
Though we do not know that (i$^{\prime}$) $\sim$ (v$^{\prime}$) 
are the fundamental relations, we will show that it is possible to calculate
ratios of critical values of $L$-functions rigorously.

Let $\mathcal F^{\ast}$ be the free group on three letters
$\widetilde\sigma$, $\widetilde\nu$, $\widetilde\tau$.
We define a surjective homomorphism
$\pi ^{\ast}: \mathcal F^{\ast} \longrightarrow \Gamma ^{\ast}$
by  $\pi ^{\ast}(\widetilde\sigma )=\sigma$, $\pi ^{\ast}(\widetilde\nu )=\nu$,
$\pi ^{\ast}(\widetilde\tau )=\tau$. Let $R^{\ast}$ be the kernel of $\pi ^{\ast}$.
Then $R^{\ast}$ contains the elements (i$^{\ast}$) $\sim$ (iv$^{\ast}$) in \S6.1 and
\begin{equation}
\widetilde\nu ^2\widetilde\tau\widetilde\nu ^{-2}
\{ \widetilde\tau (\widetilde\nu\widetilde\tau\widetilde\nu ^{-1})^3\} ^{-1}.
\tag{v$^{\ast}$}
\end{equation}
For every $\gamma \in \Gamma ^{\ast}$, we choose
$\widetilde\gamma \in \mathcal F^{\ast}$ so that
$\pi ^{\ast}(\widetilde\gamma )=\gamma$.
We use the same algorithm as in the previous section.

We consider $f^+$ (cf. \S5.6). We put $f^{\ast}=\widetilde T(f)$.
Then $f^{\ast} \in Z^2(\Gamma ^{\ast}, V)$ and
$f^{\ast}\vert \Gamma=f^+$ (cf. \S5.3).
Let $\varphi \in H^1(R^{\ast}, V)^{\Gamma ^{\ast}}$ be a corresponding
element to $f^{\ast}$. We may assume that (6.4) and (6.5) hold.
We may also assume that $\varphi (\widetilde\sigma ^2)=0$.
We need to analyze the process of adding $h\vert R^{\ast}$ to $\varphi$.
For $S$, $T$, $U \in V$, there exists $h \in H^1(\mathcal F^{\ast}, V)$ such that
$$
h(\widetilde\sigma )=S, \qquad h(\widetilde\tau )=T, \qquad
h(\widetilde\nu )=U.
$$
We find that the conditions for $h$ to vanish on the elements
(iv$^{\ast}$) and (v$^{\ast}$) are (6.7) and
\begin{equation*}
\begin{aligned}
& [\nu ^2-\tau \{ 1+\nu\tau\nu ^{-1}+(\nu\tau\nu ^{-1})^2\} \nu -1]T \\
+ & [(1-\nu ^2\tau\nu ^{-2})(1+\nu )-\tau \{1+\nu\tau\nu ^{-1}+(\nu\tau\nu ^{-1})^2\}
(1-\nu\tau\nu ^{-1}) ]U=0
\end{aligned}
\tag{7.1}
\end{equation*}
respectively. We put
$$
A=\varphi ((\widetilde\sigma\widetilde\nu )^2), \qquad
B=\varphi ((\widetilde\sigma\widetilde\tau )^3).
$$
Then (6.10) holds. As in the previous section, our objective is to determine $A$ explicitly.

\medskip

{\bf 7.2.} Let us consider the Hecke operators.
We put $g^{\ast}=T(\varpi )f^{\ast}$ where $g^{\ast}$
is defined by (5.4) with $\Gamma ^{\ast}$ in place of $\Gamma$.
Let $\psi \in H^1(R^{\ast}, V)^{\Gamma ^{\ast}}$ be a corresponding
element to $g^{\ast}$. We may assume that $\psi$ is given by (6.11).

We have $3=(4+\sqrt {13})(4-\sqrt {13})$ in $F$. Put
$\varpi =4-\sqrt {13}=-2\epsilon +7$, $\mathfrak p=(\varpi )$
and we consider the Hecke operator $T(\mathfrak p)=T(\varpi )$.
We may take
\begin{equation*}
\beta _1=\begin{pmatrix} 1 & 0 \\ 0 & \varpi \end{pmatrix}, \quad
\beta _2=\begin{pmatrix} 1 & 1 \\ 0 & \varpi \end{pmatrix}, \quad
\beta _3=\begin{pmatrix} 1 & \epsilon \\ 0 & \varpi \end{pmatrix}, \quad
\beta _4=\begin{pmatrix} \varpi & 0 \\ 0 & 1 \end{pmatrix}.
\end{equation*}
Using (6.11), we can compute $\psi (\widetilde\sigma ^2)$,
$\psi ((\widetilde\sigma\widetilde\nu )^2)$ and
$\psi ((\widetilde\sigma\widetilde\tau )^3)$.
Remarkably it turns out that these quantities can be expressed by
$A$ and $B$. Since this is technically the essential part of calculation,
we are going to explain the computation of $\psi ((\widetilde\sigma\widetilde\tau )^3)$
in some detail. By (6.11), we have
$$
\psi ((\widetilde\sigma\widetilde\tau )^3)=c\beta _3^{-1}Z_3,
$$
where
\begin{equation}
Z_3=\varphi ((\widetilde\sigma
\widetilde{\begin{pmatrix} \epsilon ^{-1} & 2\epsilon -7 \\ 0 & \epsilon \end{pmatrix}}
\widetilde\sigma
\widetilde{\begin{pmatrix} 1 & -2\epsilon \\ 0 & 1 \end{pmatrix}}
)^3).
\tag{7.2}
\end{equation}

For $x \in \mathcal O_F$ and $u \in E_F$ such that $x$ divides $u-1$, we put
\begin{equation*}
\begin{aligned}
& \{ x, u \} _4 \\
= &
\widetilde{\begin{pmatrix} 1 & x \\ 0 & 1 \end{pmatrix}}
\widetilde\sigma
\widetilde{\begin{pmatrix} 1 & (1-u)/x \\ 0 & 1 \end{pmatrix}}
\widetilde\sigma
\widetilde{\begin{pmatrix} 1 & -x/u \\ 0 & 1 \end{pmatrix}}
\widetilde\sigma
\widetilde{\begin{pmatrix} 1 & -u(1-u)/x \\ 0 & 1 \end{pmatrix}}
\widetilde\sigma
\widetilde{\begin{pmatrix} u^{-1} & 0 \\ 0 & u \end{pmatrix}}.
\end{aligned}
\end{equation*}
Then $\{ x, u \} _4 \in R^{\ast}$. As a quantitative version of Lemma A.6, (3)
of Appendix, we can show that
\begin{equation*}
\begin{aligned}
& \varphi (\{ x, u^e \} _4)=\varphi (\{ x, u \} _4) \\
+ \, & \sigma\begin{pmatrix} u^{-1} & u^{-e+1}(1-u^e)/x \\ 0 & u \end{pmatrix}\sigma
\begin{pmatrix} 1 & u^{e-2}x \\ 0 & 1 \end{pmatrix}
\varphi (\{ -u^{e-2}x, u^{e-1}\} _4) \\
- \, & \sigma\begin{pmatrix} u^{-1} & u^{-e+1}(1-u^e)/x \\ 0 & u \end{pmatrix}\sigma
\varphi ((\widetilde\sigma
\widetilde{\begin{pmatrix} u^{1-e} & 0 \\ 0 & u^{e-1} \end{pmatrix}})^2) \\
+ \, & \sigma\begin{pmatrix} u^{-e} & 0 \\ 0 & u^e \end{pmatrix}
\varphi ((\widetilde\sigma
\widetilde{\begin{pmatrix} u^{-e} & 0 \\ 0 & u^e \end{pmatrix}})^2) \\
- \, & \sigma\begin{pmatrix} u^{-1} & 0 \\ 0 & u \end{pmatrix}
\varphi ((\widetilde\sigma
\widetilde{\begin{pmatrix} u^{-1} & 0 \\ 0 & u \end{pmatrix}})^2)
\end{aligned}
\tag{7.3}
\end{equation*}
for $e \in \mathbf Z$. (This formula holds for any real quadratic field $F$.)
By (7.3) and using the formulas given in \S6.5, we can express
$\psi ((\widetilde\sigma\widetilde\tau )^3)$ in terms of $A$ and $B$.

\medskip

{\bf 7.3.} We assume $\varphi \in H^1(R^{\ast}, V)^{\Gamma ^{\ast},+}$
(cf. \S5.5). Then, as in \S6.6, we have
$$
(\delta\tau +1)B=0, \qquad (\delta -1)A=0.
$$

{\bf Fact 1.} Suppose $0\leq l_2\leq l_1 \leq 20$.
Then adding $h \vert R^{\ast}$, $h \in H^1(\mathcal F^{\ast}, V)$
to $\varphi$ (keeping $\varphi$ in the plus space under the action of $\delta$),
we may assume $B=0$.

\medskip

Therefore our task is to find constraints on
$A=\varphi ((\widetilde\sigma\widetilde\nu )^2)$. We put
$x=\sigma\begin{pmatrix} \epsilon ^{-1} & 2\epsilon -7 \\ 0 & \epsilon \end{pmatrix}
\sigma\begin{pmatrix} 1 & -2\epsilon \\ 0 & 1 \end{pmatrix}$ and let
\begin{equation}
Z_A^+=\{ \mathbf v \in V \mid (\sigma\nu -1)\mathbf v=0, \ (\delta -1)\mathbf v=0, \
xZ_3=Z_3 \} .
\tag{7.4}
\end{equation}
Here $Z_3$ is defined by (7.2) and the meaning of $xZ_3=Z_3$ is the same as in \S6.6.
Namely, $xZ_3=Z_3$ must hold because $x^3=1$;
since $Z_3$ can be expressed by $A$, $xZ_3=Z_3$ gives a constraint on $A$.

We consider the contribution of $H^1(\mathcal F^{\ast}, V)$ to $Z_A^+$.
Take $h \in H^1(\mathcal F^{\ast}, V)$ and put
$$
h(\widetilde\sigma )=S, \qquad h(\widetilde\nu )=U, \qquad
h(\widetilde\tau )=T.
$$
We require that $h\vert R^{\ast}$ vanishes on the elements
(i$^{\ast}$), (ii$^{\ast}$), (iv$^{\ast}$), (v$^{\ast}$).
These conditions are equivalent to (6.21), (6.22), (6.7) and
and (7.1). We have
$$
h((\widetilde\sigma\widetilde\nu )^2)=(\sigma\nu +1)(\sigma U+S).
$$
We also require that (6.23) holds.
Let $B_A^+$ be the subspace of $V$ generated by
$(\sigma\nu +1)(\sigma U+S)$ when $S$, $T$, $U$ extend over vectors of $V$
satifying the relations (6.7), (6.21), (6.22), (6.23) and (7.1). 
We have $B_A^+ \subset Z_A^+$. As shown in \S4.1, (6.24) holds.
By Proposition 5.5, we have
$$
\dim Z_A^+/B_A^+ \geq \dim S_{l_1+2.l_2+2}(\Gamma ) \qquad
{\rm if} \ l_2 \geq 4, \  l_1 \neq l_2 \ {\rm or \ if} \ \  l_1=l_2,  \ l_2 \geq 6.
$$
Now by numerical computations, we have verified:

\medskip

{\bf Fact 2.} Suppose $0 \leq l_2\leq l_1 \leq 20$.
Then $\dim S_{l_1+2, l_2+1}(\Gamma )=\dim Z_A^+/B_A^+$.

\medskip

This fact means that the constraints posed on
$A=\varphi ((\widetilde\sigma\widetilde\nu )^2)$ is enough.

\medskip

{\bf Example 7.1.} We take $l_1=l_2=6$.
We have $\dim S_{8,8}(\Gamma )=5$.
Calculating the action of $T(\mathfrak p)$ on $Z_A^+/B_A^+$ using (6.11),
we find that the characteristic polynomial of $T(\mathfrak p)$ is
$$
(X^2-40X-3957)(X^3+28X^2-2601X-71748).
$$
The quadratic factor corresponds to the non base change part;
the irreducible factor of degree three corresponds to the base change part
from $S_8(\Gamma _0(13), (\frac {}{13}))$.
Let $\Omega \in S_{8,8}(\Gamma )$ be the Hecke eigenform
such that $\Omega\vert T(\mathfrak p)=(20+\sqrt {4357})\Omega$.
Then we find
$$
R(6, \Omega )/R(4, \Omega )=70/3.
$$

\medskip

{\bf Example 7.2.} We take $l_1=l_2=8$.
We have $\dim S_{10,10}(\Gamma )=7$.
We find that the characteristic polynomial of $T(\mathfrak p)$ is
$$
(X^2-16X-42789)(X^5+X^4-66033X^3+1260423X^2+530326440X+14266185264).
$$
The quadratic factor corresponds to the non base change part.
Let $\Omega \in S_{10,10}(\Gamma )$ be the Hecke eigenform
such that $\Omega\vert T(\mathfrak p)=(8+\sqrt {42853})\Omega$.
Then we find
$$
R(7, \Omega )/R(5, \Omega )=50.
$$

\medskip

{\bf Example 7.3.} We take $l_1=l_2=10$.
We have $\dim S_{12,12}(\Gamma )=11$.
We find that the characteristic polynomial of $T(\mathfrak p)$ is
\begin{equation*}
\begin{aligned}
& (X-252)(X^4+252X^3-496198X^2-116604684X+25202349477) \\
& (X^6+244X^5-665334X^4-129598956X^3+109163403621X^2 \\
& +14522233287672X-255121008509808).
\end{aligned}
\end{equation*}
The irreducible factor of degree four corresponds to the non base change part;
$X-252$ corresponds to the base change part from $S_{12}({\rm SL}_2(\mathbf Z))$
and the irreducible factor of degree six corresponds to the base change part from
$S_{12}(\Gamma _0(13), (\frac {}{13}))$.
Put
$$
f(X)=X^4+252X^3-496198X^2-116604684X+25202349477.
$$
Let $\theta$ be a root of $f(X)$ and put $K=\mathbf Q(\theta )$.
We find that $K$ contains a quadratic subfield
$F=\mathbf Q(\sqrt {7\cdot 5167})$. Put $d=7\cdot 5167$.
Then a root of $f(X)$ is given by
$$
\psi =-(63+\sqrt {d})+\sqrt {223837-360\sqrt {d}}.
$$
We have
$$
N(223837-360\sqrt {d})=13\cdot 563\cdot 6205151.
$$
This number and the quadratic fields in Examples 7.1 and 7.2
are consistent with the table given in Doi-Hida-Ishii [DHI].

For the Hecke eigenform $\Omega \in S_{12,12}(\Gamma )$ such that
$\Omega \vert T(\varpi )=\psi\Omega$, we find
$$
R(10, \Omega )/R(6, \Omega )=\frac {3732099+18663\sqrt {d}}{5},
$$
$$
R(8, \Omega )/R(6, \Omega )=\frac {24367+121\sqrt {d}}{20}.
$$

\bigskip

\centerline {\S8. A comparison of two methods}

\medskip

In [Sh3], Shimura gave a method to calculate critical values of
$D(s, f, g)$ for two elliptic modular forms $f$ and $g$.
Here $D(s, f, g)$ is the Rankin-Selberg convolution of $f$ and $g$.
Shortly later he gave a generalization to the case of Hilbert modular forms ([Sh4]).
Taking one argument in the convoluted $L$-function as a suitable Eisenstein series,
this method enables us to calculate the ratios of critical values of $L(s, \Omega )$
for a Hilbert modular form $\Omega$. We call this technique method A.
We call the cohomological technique method B, which was initiated in [Sh1]
and studied in this paper when $[F:\mathbf Q]=2$.
It is interesting to compare A and B.

(0) Method A is more general and conceptually simpler.
It has the advantage to give the relation of the product of
the plus and minus periods to the Petersson norm. It is applicable
also to modular forms of half integral weights.

(1) If $n=[F:\mathbf Q]>2$, the method B has to calculate
$H^n(\Gamma , V)$, which is beyond the reach at present.
Therefore when $[F:\mathbf Q]>2$, A is definitively superior than B.

(2) Suppose that $[F:\mathbf Q]=2$. The method B is still incomplete.
But in the cases well worked out, $F={\mathbf Q}(\sqrt {5})$ for example,
B has the advantage that we can write a program which calculates everything
by machine. It can also be used to calculate the characteristic polynomials
of Hecke operators. (In this respect, it is desirable to solve the problem
mentioned at the end of subsection 6.5.) We employed essentially a single program
to obtain examples in section 7. Therefore in some cases at least,
B will have the advantage over A. But in general the method A is conceptually simpler.

In Doi-Goto [DG] and Doi-Ishii [DI], the authors gave interesting examples
of critical values of $D(s, f, g)$ for Hilbert modular forms $f$ and $g$.
Their interests was the relation of this value to the congruences between
Hilbert modular forms. However they did not give examples of
critical values of $L(s, \Omega )$. Recently Dr. K. Okada calculated
the ratios of critical values of $L(s, \Omega )$ and confirmed
the numerical value of Example 7.1 by method A. He obtained
one more example for $F={\mathbf Q}(\sqrt {17})$.

(3) Suppose that $F=\mathbf Q$. The method B is developed into
the theory of modular symbols which is presently used to calculate
characteristic polynomials of Hecke operators. For the $L$-values,
the author doesn't know which is faster. But the calculation of [Sh1]
reviewed in the introduction suggests that B would not be more complex
than A.

\bigskip

\centerline{\bf Appendix. Generators and relations}

\medskip

Let $F$ be a real quadratic field and $\epsilon$ be the fundamental unit of $F$.
Let $\{ 1, \omega \}$ be an integral basis of $\mathcal O_F$, i.e.,
$\mathcal O_F=\mathbf Z \oplus \mathbf Z\omega$. We write
\begin{equation}
\epsilon ^2=A+B\omega, \qquad
\epsilon ^2\omega =C+D\omega .
\tag{A.1}
\end{equation}
We put $\Gamma ={\rm PSL}(2, \mathcal O_F)$,
$\widetilde\Gamma ={\rm SL}(2, \mathcal O_F)$,
$$
\widetilde P=\left\{ \begin{pmatrix} a & b \\ 0 & a^{-1} \end{pmatrix} \bigg| \
a \in E_F, b \in \mathcal O_F \right\} , \qquad
P=\widetilde P/\{ \pm 1_2 \} .
$$
We define elements of $\widetilde\Gamma$ by
$$
\sigma =\begin{pmatrix} 0 & 1 \\ -1 & 0 \end{pmatrix}, \quad
\mu =\begin{pmatrix} \epsilon & 0 \\ 0 & \epsilon ^{-1} \end{pmatrix}, \quad
\tau =\begin{pmatrix} 1 & 1 \\ 0 & 1 \end{pmatrix}, \quad
\eta =\begin{pmatrix} 1 & \omega \\ 0 & 1 \end{pmatrix}.
$$
Then it is known that $\sigma$, $\mu$, $\tau$ and $\eta$ generate $\widetilde\Gamma$
(cf. Vaser$\rm {\check{s}}$tein [V]). This fact can be proved in elementary way if $\mathcal O_F$ is
a Euclidean ring, $F=\mathbf Q(\sqrt {5})$ for example.
We use same letters $\sigma$, $\mu$, $\tau$ and $\eta$ for their classes in $\Gamma$,
since this will cause no confusion. Now we have relations among them:
\begin{equation}
\sigma ^2=1.
\tag{i}
\end{equation}
\begin{equation}
(\sigma\tau ) ^3=1.
\tag{ii}
\end{equation}
\begin{equation}
(\sigma\mu ) ^2=1.
\tag{iii}
\end{equation}
\begin{equation}
\tau\eta =\eta\tau .
\tag{iv}
\end{equation}
\begin{equation}
\mu\tau\mu ^{-1}=\tau ^A\eta ^B.
\tag{v}
\end{equation}
\begin{equation}
\mu\eta\mu ^{-1}=\tau ^C\eta ^D.
\tag{vi}
\end{equation}
If we can take $\omega =\epsilon$ and
$-\epsilon ^{-1}=A^{\prime}+B^{\prime}\epsilon$, then we have
\begin{equation}
\sigma\eta\sigma=\tau ^{A^{\prime}}\eta ^{B^{\prime}}\sigma \eta ^{-1}\mu .
\tag{vii}
\end{equation}
The relations (ii) and (vii) follow from
\begin{equation}
\sigma\begin{pmatrix} 1 & t \\ 0 & 1 \end{pmatrix}\sigma
=\begin{pmatrix} 1 & -t^{-1} \\ 0 & 1 \end{pmatrix}\sigma
\begin{pmatrix} -t & 1\\ 0 & -t^{-1} \end{pmatrix},
\quad t \in E_F.
\tag{A.2}
\end{equation}
It is easy to see that $\mu$, $\tau$ and $\eta$ generate $P$ and (iv) $\sim$ (vi) are
their fundamental relations.

The purpose of this appendix is to prove the following theorem.

\proclaim Theorem A.1. Let $F=\mathbf Q(\sqrt {5})$ and
$\Gamma ={\rm PSL}(2, \mathcal O_F)$. We take $\omega =\epsilon$.
The fundamental relations satisfied by the generators
$\sigma$, $\mu$, $\tau$ and $\eta$ are (i) $\sim$ (vii).
\par
We note that if $F=\mathbf Q(\sqrt {5})$ then
$A=1$, $B=1$, $C=1$, $D=2$, $A^{\prime}=1$, $B^{\prime}=-1$.
The relations (i) to (vi) and (A.2) hold for any real quadratic field.
Our theorem states that the minimal relations are enough when
$F={\mathbf Q}(\sqrt {5})$. This minimality will be satisfied
by some more real quadratic fields with small discriminants but
will not hold in general.

We begin by preliminary considerations on generators and relations of $\Gamma$.
~\footnote{For this part, we do not assume $F=\mathbf Q(\sqrt{5})$.}
Since $\Gamma$ is generated by $P$ and $\sigma$, every relation among elements of $P$ and $\sigma$
takes the form
$$
p_1\sigma p_2\sigma \cdots p_m \sigma =1, \qquad p_i \in P, \ 1 \leq i \leq m.
$$
Using (i) and (iii) $\sim$ (vi), this relation can be written as
$$
\begin{pmatrix} 1 & x_1 \\ 0 & 1 \end{pmatrix}\sigma
\begin{pmatrix} 1 & x_2 \\ 0 & 1 \end{pmatrix}\sigma \cdots
\begin{pmatrix} 1 & x_m \\ 0 & 1 \end{pmatrix}\sigma
=\begin{pmatrix} u & 0 \\ 0 & u^{-1} \end{pmatrix} , \quad
x_i \in \mathcal O_F, \ u \in E_F.
$$
We call a relation of this type an {\it $m$ terms relation} counting
the number of $\sigma$ involved.

\proclaim Lemma A.2. Using relations (i) and (iii) $\sim$ (vi), every three terms relation
can be reduced to (A.2).
\par
{\bf Proof.} If we have a two terms relation
$$
\begin{pmatrix} 1 & x_1 \\ 0 & 1 \end{pmatrix}\sigma
\begin{pmatrix} 1 & x_2 \\ 0 & 1 \end{pmatrix}\sigma
=\begin{pmatrix} u & 0 \\ 0 & u^{-1} \end{pmatrix},
$$
we have $x_1=x_2=0$, $u=\pm 1$. Hence the two terms relation reduces to (i).
Let
$$
\begin{pmatrix} 1 & x_1 \\ 0 & 1 \end{pmatrix}\sigma
\begin{pmatrix} 1 & x_2 \\ 0 & 1 \end{pmatrix}\sigma
\begin{pmatrix} 1 & x_3 \\ 0 & 1 \end{pmatrix}\sigma
=\begin{pmatrix} u & 0 \\ 0 & u^{-1} \end{pmatrix}
$$
be a three terms relation. Then we see that $x_2=\pm u \in E_F$. Using (A.2), we have
$\sigma\begin{pmatrix} 1 & x_2 \\ 0 & 1 \end{pmatrix}\sigma =p_1\sigma p_2$
with some $p_1$, $p_2 \in P$ and the three terms relation in question reduces
to a two terms relation. This completes the proof.

\proclaim Lemma A.3. Assume that we can take $\omega =\epsilon$.
The relation (A.2) can be reduced to the relations (i) $\sim$ (vii).
In other words, the relation (A.2) for $t \in E_F$ can be reduced
to the relations (A.2) for $t=1$, $\epsilon$ using relations (i) and (iii) $\sim$ (vi).
\par
{\bf Proof.} We write the relation (A.2) as $\{ t \}$.
Using (i), the relation (iii) implies the relation
$\begin{pmatrix} u & 0 \\ 0 & u^{-1} \end{pmatrix}\sigma
=\sigma \begin{pmatrix} u^{-1} & 0 \\ 0 & u \end{pmatrix}$ for $u \in E_F$.
Then we obtain the relation $\{ -t \}$ taking the inverse of the both sides of (A.2),
using (i), (iv) $\sim$ (vi). Taking the conjugate by $\mu$ of both sides of (A.2),
we obtain the relation $\{ \epsilon ^{-2}t\}$ using (i), (iii) $\sim$ (vi).
Since $E_F$ is generated by $\epsilon$ and $\pm 1$, this completes the proof.

\medskip

Next we consider the four terms relation.
\begin{equation}
\begin{pmatrix} 1 & x_1 \\ 0 & 1 \end{pmatrix}\sigma
\begin{pmatrix} 1 & x_2 \\ 0 & 1 \end{pmatrix}\sigma
\begin{pmatrix} 1 & x_3 \\ 0 & 1 \end{pmatrix}\sigma
\begin{pmatrix} 1 & x_4 \\ 0 & 1 \end{pmatrix}\sigma
=\begin{pmatrix} u & 0 \\ 0 & u^{-1} \end{pmatrix}
\tag{A.3}
\end{equation}
We write the relation (A.3) as $\{ x_1, x_2, x_3, x_4; u \}$.

\proclaim Lemma A.4. The four terms relation (A.3) reduces to
(i) $\sim$ (vi) and (A.2) if $x_i \in E_F$ for some $i$, $1 \leq i \leq 4$.
\par
{\bf Proof.} Suppose that $x_2 \in E_F$. By (A.2), we have
$\sigma\begin{pmatrix} 1 & x_2 \\ 0 & 1 \end{pmatrix}\sigma =p_1\sigma p_2$
with some $p_1$, $p_2 \in P$. Using this expression, we find that (A.3) reduces to
a three terms relation. We write (A.3) as
$$
\begin{pmatrix} 1 & x_2 \\ 0 & 1 \end{pmatrix}\sigma
\begin{pmatrix} 1 & x_3 \\ 0 & 1 \end{pmatrix}\sigma
\begin{pmatrix} 1 & x_4 \\ 0 & 1 \end{pmatrix}\sigma
=\sigma\begin{pmatrix} 1 & -x_1 \\ 0 & 1 \end{pmatrix}
\begin{pmatrix} u & 0 \\ 0 & u^{-1} \end{pmatrix} .
$$
Using (i) $\sim$ (vi), the right-hand side can be written as
$\begin{pmatrix} u^{-1} & 0 \\ 0 & u \end{pmatrix}\sigma
\begin{pmatrix} 1 & -u^{-2}x_1 \\ 0 & 1 \end{pmatrix}$.
Hence $\{ x_1, x_2, x_3, x_4; u \}$ is equivalent to
$\{ x_2, x_3, x_4, u^{-2}x_1; u^{-1} \}$ under (i) $\sim$ (vi).
By this cyclic rotation, any $x_i$ can be brought to the second position
at the cost of multiplying by a unit. Hence the assertion follows.

\medskip

For $u \in E_F$, $x \in \mathcal O_F$, we have the relation
\begin{equation*}
\begin{aligned}
& \begin{pmatrix} 1 & x \\ 0 & 1 \end{pmatrix}\sigma
\begin{pmatrix} 1 & (1-u)/x \\ 0 & 1 \end{pmatrix}\sigma
\begin{pmatrix} 1 & -x/u \\ 0 & 1 \end{pmatrix}\sigma
\begin{pmatrix} 1 & -u(1-u)/x \\ 0 & 1 \end{pmatrix}\sigma \\
= & \begin{pmatrix} u & 0 \\ 0 & u^{-1} \end{pmatrix}
\end{aligned}
\tag{A.4}
\end{equation*}
if $x$ divides $u-1$.

\proclaim Lemma A.5. Under (i) $\sim$ (vi) and (A.2), the four terms relation
(A.3) can be reduced to (A.4) with some $x$ and $u$.
\par
{\bf Proof.} We see easily that the four terms relation (A.3) is equivalent to a relation of the form
\begin{equation}
\sigma\begin{pmatrix} 1 & x \\ 0 & 1 \end{pmatrix}\sigma
=\begin{pmatrix} 1 & y_1 \\ 0 & 1 \end{pmatrix}\sigma
\begin{pmatrix} 1 & y_2 \\ 0 & 1 \end{pmatrix}\sigma
\begin{pmatrix} 1 & y_3 \\ 0 & 1 \end{pmatrix}
\begin{pmatrix} h & 0 \\ 0 & h^{-1} \end{pmatrix} .
\tag {A.3$^\prime$}
\end{equation}
Here $x$, $y_i \in \mathcal O_F$, $1 \leq i \leq 3$ and $h \in E_F$.
By a direct computation, we get
$$
h(y_1y_2-1)=-\omega , \qquad hy_2=\omega x, \qquad h^{-1}(y_2y_3-1)=-\omega ,
$$
where $\omega =\pm 1$. Putting $u=\omega h^{-1}$, we have
$$
y_2=ux, \qquad y_1=\frac {1-u}{ux}, \qquad y_3=\frac {1-u^{-1}}{ux}.
$$
Hence we see that $x$ divides $u-1$ and that (A.3$^\prime )$ is equivalent to
\begin{equation*}
\begin{aligned}
& \sigma\begin{pmatrix} 1 & x \\ 0 & 1 \end{pmatrix}\sigma \\
= & \begin{pmatrix} 1 & (1-u)/ux \\ 0 & 1 \end{pmatrix}\sigma
\begin{pmatrix} 1 & ux \\ 0 & 1 \end{pmatrix}\sigma
\begin{pmatrix} 1 & (1-u^{-1})/ux \\ 0 & 1 \end{pmatrix}
\begin{pmatrix} u^{-1} & 0 \\ 0 & u \end{pmatrix} .
\end{aligned}
\tag {A.3$^{\prime\prime}$}
\end{equation*}
On the other hand, under (i) $\sim$ (vi), (A.4) is equivalent to
\begin{equation*}
\begin{aligned}
& \sigma\begin{pmatrix} 1 & -x \\ 0 & 1 \end{pmatrix}\sigma \\
= & \begin{pmatrix} 1 & (1-u)/x \\ 0 & 1 \end{pmatrix}\sigma
\begin{pmatrix} 1 & -x/u \\ 0 & 1 \end{pmatrix}\sigma
\begin{pmatrix} 1 & -u(1-u)/x \\ 0 & 1 \end{pmatrix}
\begin{pmatrix} u & 0 \\ 0 & u^{-1} \end{pmatrix} .
\end{aligned}
\tag {A.4$^\prime$}
\end{equation*}
We obtain (A.3$^{\prime\prime})$ from (A.4$^\prime$) by substituting
$x$ by $-x$ and $u$ by $u^{-1}$. This completes the proof.

\medskip

We denote the four terms relation (A.4) by $\{ x, u \}$.
We have $\{ x, u \}=\{x, (1-u)/x, -x/u, -u(1-u)/x; u \}$.
Under (i) $\sim$ (vi), the relation of the form (A.3$^\prime$) is
equivalent to $\{ x, u \}$ and the relation $\{ x_1, x_2, x_3, x_4; u \}$ is
equivalent to $\{ x_2, x_3, x_4, u^{-2}x_1; u^{-1} \}$
(cf. the proofs of Lemmas A.4 and A.5). Therefore $\{ x, u \}$ is equivalent to
$\{ (1-u)/x, u^{-1} \}$ under (i) $\sim$ (vi).
By Lemma A.4, $\{ x, u \}$ is reducible to (i) $\sim$ (vi) and (A.2) if
$x \in E_F$ or $(1-u)/x \in E_F$. 

\proclaim Lemma A.6. Assuming (i) $\sim$ (vi) and (A.2), the following assertions hold.
\begin{enumerate}
\item [(1)] $\{ x, u \}$ is equivalent to $\{ -x, u^{-1} \}$.
\item [(2)] $\{ x, u \}$ is equivalent to $\{ t^2x, u \}$ for every $t \in E_F$.
\item [(3)] We assume the four terms relation $\{ x, u \}$.
Then $\{ x, u^e \}$ is equivalent to $\{ u^ex, u^{1-e} \}$ for $e \in \mathbf Z$.
\item [(4)] $\{ x, u \}$ is equivalent to $\{ (1-u)/x, u^{-1} \}$.
\item [(5)] Suppose that $(x)=(2)$. Then $\{ x, u \}$ is equivalent to $\{ x, -u \}$.
\end{enumerate}
\par
{\bf Proof.} We write $\{ -x, u^{-1} \}$ in the form of (A.3$^{\prime\prime}$).
Taking the inverses of both sides, we obtain (1).
We obtain (2) taking the conjugates of both sides by
$\begin{pmatrix} t^{-1} & 0 \\ 0 & t \end{pmatrix}$.
To prove (3), we set the right-hand side of (A.3$^{\prime\prime}$)
is equal for $\{ x, u \}$ and for $\{ x, u^{e} \}$. By a simple computation,
we find that the resulting equality is
\begin{equation*}
\begin{aligned}
& \sigma\begin{pmatrix} 1 & u^ex \\ 0 & 1 \end{pmatrix}\sigma \\
= & \begin{pmatrix} 1 & (u^{e-1}-1)/u^ex \\ 0 & 1 \end{pmatrix}\sigma
\begin{pmatrix} 1 & ux \\ 0 & 1 \end{pmatrix}\sigma
\begin{pmatrix} 1 & (u^{-1}-u^{e-2})/x \\ 0 & 1 \end{pmatrix}
\begin{pmatrix} u^{e-1} & 0 \\ 0 & u^{1-e} \end{pmatrix} ,
\end{aligned}
\end{equation*}
which is $\{ u^ex, u^{1-e} \}$. Hence we obtain (3).
We noted (4) already in the discussion before Lemma A.6.
To prove (5), we set the right-hand side of (A.3$^{\prime\prime}$)
is equal for $\{ x, u \}$ and for $\{ x, -u \}$. The resulting equality is
\begin{equation*}
\begin{aligned}
& \sigma\begin{pmatrix} 1 & ux \\ 0 & 1 \end{pmatrix}\sigma \\
= & \begin{pmatrix} 1 & -2/ux \\ 0 & 1 \end{pmatrix}\sigma
\begin{pmatrix} 1 & -ux \\ 0 & 1 \end{pmatrix}\sigma
\begin{pmatrix} 1 & -2/ux \\ 0 & 1 \end{pmatrix}
\begin{pmatrix} -1 & 0 \\ 0 & -1 \end{pmatrix} ,
\end{aligned}
\end{equation*}
Since $-2/ux \in E_F$, this relation reduces to a three terms relation
by Lemma A.4. In view of Lemma A.2, this completes the proof.

\medskip

{\bf Remark A.7.} Suppose that $(1-u)/x \in E_F$. Then, by Lemma A.4,
$\{ tx, u \}$ can be reduced to (i) $\sim$ (vi) and (A.2) for every $t \in E_F$.
By (1) and (3) of Lemma A.6, we see that $\{ x, u^e \}$ can be reduced to
(i) $\sim$ (vi) and (A.2) for all $e \in \mathbf Z$.

\medskip

The following Lemma is of some interest though it will not be used in this paper.

\proclaim Lemma A.8. Suppose that there exist sequences of integers
$x_0$, $x_1$, $\ldots$, $x_k \in \mathcal O_F$ and units
$u_0$, $u_1$, $\ldots$, $u_k \in E_F$ such that
$$
x_{i-1}x_i=1-u_i, \quad 1 \leq i \leq k.
$$
We assume that $u_i=u_{i-1}^{m_i}$, $1 \leq i \leq k$ with a nonzero
interger $m_i$. If $(1-u_0)/x_0 \in E_F$, then the four terms relation
$\{ x_k, u_k \}$ reduces to (i) $\sim$ (vi) and (A.2).
\par
{\bf Proof.} Using Lemma A.6, the reducibility of $\{ tx_i, u_i^e \}$,
$t \in E_F$, $e \in \mathbf Z$ can be shown easily by induction on $i$.

\medskip

Let $G$ be a group with generators $\sigma _1, \ldots , \sigma _m$.
Let $\mathcal F$ be a free group on the free generators
$\widetilde\sigma _1, \ldots , \widetilde\sigma _m$.
Then we can define a surjective homomorphism $\pi : \mathcal F \longrightarrow G$ by
$\pi (\widetilde\sigma _i)=\sigma _i$, $1 \leq i \leq m$.
Let $R$ be the kernel of $\pi$.
Next let $S$ be a finite subset of $G$ which generates $G$.
For $\gamma \in S$, we prepare a symbol $[\gamma ]$ and let $\mathcal F^{\prime}$
be the free group on the free generators $[\gamma ]$, $\gamma \in S$.
We can define a surjective homomorphism
$\pi ^{\prime}: \mathcal F^{\prime} \longrightarrow G$ by
$\pi ^{\prime}([\gamma ])=\gamma$, $\gamma \in S$. Let $R^{\prime}$ be
the kernel of $\pi ^{\prime}$. Clearly
$([\gamma _1][\gamma _2])^{-1}[\gamma _1\gamma _2] \in R^{\prime}$
if $\gamma _1$, $\gamma _2$, $\gamma _1\gamma _2 \in S$.
We assume that $R^{\prime}$ is generated by the elements of this form and
their conjugates.

Now for every $\gamma \in S$, we take and fix an expression
$$
\gamma =\sigma _{i_1}^{\epsilon _1}\cdots \sigma _{i_k}^{\epsilon _k},
\qquad i_j \in [1, m], \quad \epsilon _j=\pm 1
$$
and put
$\widetilde\gamma =\widetilde\sigma _{i_1}^{\epsilon _1}\cdots
\widetilde\sigma _{i_k}^{\epsilon _k}$.
(If $\gamma =\sigma _i \in S$, we put $\widetilde\gamma =\widetilde\sigma _i$.)
By the universality of the free group, there exists a homomorphism
$\varphi : \mathcal F^{\prime} \longrightarrow \mathcal F$
which satisfies $\varphi ([\gamma ])=\widetilde\gamma$, $\gamma \in S$.
Then we have $\pi ^{\prime}=\pi \circ \varphi$.
Let $R_0$ be the normal subgroup of $\mathcal F$ generated by
$(\widetilde\gamma _1\widetilde\gamma _2)^{-1}
\widetilde{\gamma _1\gamma _2}$, $\gamma _1$, $\gamma _2$,
$\gamma _1\gamma _2 \in S$ and their conjugates. We have
$R_0 \subset R$. Since $\varphi (R^{\prime}) \subset R_0$ by the assumption,
$\varphi$ induces the homomomorphism
$\bar\varphi : \mathcal F^{\prime}/R^{\prime} \longrightarrow \mathcal F/R_0$
which satisfies $\bar\varphi (g \mod R^{\prime})=\varphi (g) \mod R_0$,
$g \in \mathcal F^{\prime}$.

\proclaim Lemma A.9. Let the notation be the same as above.
If $\sigma _i \in S$, $1 \leq i \leq m$, then we have $R_0=R$.
\par
{\bf Proof.} Define a homomorphism
$\pi _0: \mathcal F/R_0 \longrightarrow G$ by
$\pi _0(h \mod R_0)=\pi (h)$, $h \in \mathcal F$.
Since $(\pi _0\circ\bar\varphi )(g \mod R^{\prime})
=(\pi\circ\varphi )(g)=\pi ^{\prime}(g)$, $g \in \mathcal F^{\prime}$,
$\pi _0\circ\bar\varphi$ is injective. Hence
$\pi _0\vert \bar\varphi (\mathcal F^{\prime}/R^{\prime})$
is injective. We can write
$\bar\varphi (\mathcal F^{\prime}/R^{\prime})=H/R_0$ with
a subgroup $H$ of $\mathcal F$. Now the assumption of the Lemma implies
$H=\mathcal F$. Therefore $\pi _0$ is injective and we obtain $R_0=R$.

\medskip

For the proof of Theorem A.1, we use the following theorem of Macbeath
(cf. Theorem 1 of [Mac] and also Theorem 1.1 of [Sw]).

\proclaim Theorem M. Let $X$ be a path connected Hausdorff topological space
and $\Gamma$ be a group which acts on $X$ as homeomorphisms.
We assume that the fundamental group $\pi _1(X)$ of $X$ is trivial.
Let $V$ be a path connected open subset of $X$ such that $X=\Gamma V$.
Define a subset $S$ of $\Gamma$ by
$$
S=\{ \gamma \in \Gamma \mid V \cap \gamma V \neq \emptyset \} .
$$
Then $S$ generates $\Gamma$.
\footnote {This fact is an old result of Siegel, cf. [Si1].}
Let $\mathcal F$ be the free group
which has the symbols $[\sigma ]$, $\sigma \in S$ as free generators.
Define a homomorphism $\pi : \mathcal F \longrightarrow \Gamma$ by
$\pi ([\sigma ])=\sigma$. Let $R$ be the kernel of $\pi$.
Then $R$ is generated by $([\sigma ][\tau ])^{-1}[\sigma\tau ]$
and their conjugates, where $\sigma$ and $\tau$ are elements of $S$ which satisfy
\begin{equation}
V \cap \sigma V \cap \sigma\tau V \neq \emptyset .
\tag{$\ast$}
\end{equation}
In other words, $\Gamma$ has a presentation $\Gamma =\mathcal F/R$.
\par

Swan ([Sw]) generalized this theorem to the case where $\pi _1(X) \neq 1$
and obtained generators and relations for ${\rm SL}(2, \mathcal O_K)$,
for several imaginary quadratic fields $K$ with small discriminants.

Let the notation be the same as in Theorem M.
For a subset T of X, we put
$$
S(T)=\{ \gamma \in \Gamma \mid T \cap \gamma T \neq \emptyset \} .
$$
Let $D$ be a closed subset of $X$ such that $\Gamma D=X$.

\proclaim Lemma A.10. Suppose in addition that the topological space $X$ is normal.
Then we have
$$
\cap _{U \supset D, \ U \ {\rm is \ open}} \ S(U)=S(D).
$$
\par
{\bf Proof.} Clearly the left-hand side contains the right-hand side.
Pick an element $\gamma$ of the left-hand side. Assume that
$D\cap \gamma D=\emptyset$. Since $X$ is normal, we can find
open subsets $U$ and $U^{\prime}$ of $X$ so that
$$
U \supset D, \qquad U^{\prime} \supset \gamma D, \qquad
U \cap U^{\prime}=\emptyset .
$$
Put $U^{\prime\prime}=U \cap \gamma ^{-1}U^{\prime}$. Then we have
$U^{\prime\prime} \supset D$,
$U^{\prime\prime} \cap \gamma U^{\prime\prime}\subset U \cap U^{\prime}=\emptyset$.
This is a contradiction and we complete the proof.

\medskip

Next we assume that $S(D)$ is finite and that $S(U)$ is finite for an open set $U$
which contains $D$. We put
$$
S(D)=\{ \gamma _1, \ldots , \gamma _m \} , \qquad
S(U)=\{ \gamma _1, \ldots , \gamma _m, \gamma _{m+1}, \ldots , \gamma _n \}
$$
assuming $S(U) \supsetneqq S(D)$. By Lemma A.10, for every $\gamma _i$,
$i>m$, there exists an open set $U_i \supset D$ such that
$\gamma _i \notin S(U_i)$. Put
$V=U \cap (\cap _{i=m+1}^n U_i)$.
Then we have $\gamma _i \notin S(V)$.
Therefore we conclude that $S(D)=S(V)$ for an open set $V$
which contains $D$. This means that we may replace $S$ to $S(D)$ in Theorem M
if such a $V$ is path connected.
(Note that in Theorem M, $([\sigma ][\tau ])^{-1}[\sigma\tau ] \in R$
for $\sigma$, $\tau \in S$ such that $\sigma\tau \in S$.
Thus the condition ($\ast$) may be dropped.
However ($\ast$) reduces the number of relations and can be essential
for the practical purpose.)

\medskip

Now let $F$ be a totally real field of degree $n$. Let us review the fundamental domain
of $\Gamma ={\rm PSL}(2, \mathcal O_F)$ acting on $\mathfrak H^n$ (cf. [Si2]).
Let $\sigma _1, \ldots , \sigma _n$ be all the isomorphisms of $F$ into
$\mathbf R$. For $a \in F$, we put $a^{(i)}=a^{\sigma _i}$.
Take an integral basis of ${\mathcal O_F}$ so that
$$
\mathcal O_F=\mathbf Z\omega _2+\mathbf Z\omega _2 +\cdots +\mathbf Z\omega _n
$$
and let $\epsilon _1, \ldots , \epsilon _{n-1}$ be generators
of a free part of $E_F$.
For $x=(x_1, \ldots , x_n) \in \mathbf C^n$, we put $N(x)=x_1\cdots x_n$.
For simplicity, we assume that the class number of $F$ is one.
Take $z=(z_1, \ldots , z_n) \in \mathfrak H^n$. Put
$z_j=x_j+iy_j$, $x_j$, $y_j \in \mathbf R$.
We define the local coordinates of $z$ relative to the cusp
$\infty$ by the formulas (cf. [Si2], p. 249)
\begin{equation}
Y_1\log \vert \epsilon _1^{(k)}\vert +
\cdots +Y_{n-1}\log \vert \epsilon _{n-1}^{(k)}\vert
=\frac {1}{2}\log \frac {y_k}{\root n \of {N(y)}}, \qquad 1 \leq k \leq n-1.
\tag{A.5}
\end{equation}
\begin{equation}
X_1\omega _1^{(l)}+\cdots +X_n\omega _n^{(l)}=x_l, \qquad 1 \leq l \leq n,
\tag{A.6}
\end{equation}
Here $y=(y_1, \ldots , y_n)$. We put
$$
D_{\infty}= \{ z \in \mathfrak H^n \mid
-\frac {1}{2}\leq Y_i<\frac {1}{2}, \ 1 \le i \leq n-1, \quad
-\frac {1}{2}\leq X_j<\frac {1}{2}, \ 1 \le j \leq n \} .
$$
Then $D_{\infty}$ is a fundamental domain of $P$.
($P$ is the subgroup of $\Gamma$ consisting of all elements
which are represented by upper triangular matrices.)
We define
\begin{equation*}
\begin{aligned}
D=\{ z \in \overline {D_{\infty}} \mid & N(\vert cz+d\vert )\geq 1 \ \
\text {whenever} \\
& \text {$c$ and $d$ are relatively prime integers of $\mathcal O_F$} \} .
\end{aligned}
\tag{A.7}
\end{equation*}
Here $\overline {D_{\infty}}$ denote the closure of $D_{\infty}$ and
$\vert cz+d\vert =(\vert c^{(1)}z_1+d^{(1)}\vert , \ldots , \vert c^{(n)}z_n+d^{(n)}\vert )$.
Then $D$ satisfies that (cf. [Si2], p. 266--268):
\begin{enumerate}
\item  $D$ is a closed subset of $\mathfrak H^n$ such that $\Gamma D=\mathfrak H^n$.
\item Two distinct interior points of $D$ cannot be transformed
each other by an element of $\Gamma$.
\item There are only finitely many $\gamma \in \Gamma$ such that
$D \cap \gamma D \neq \emptyset$. Furthermore $D$ and $\gamma D$,
$\gamma \neq 1$ can intersect only on the boundary of $D$.
\end{enumerate}

Now we assume that $[F:\mathbf Q]=2$. We may assume that
$\omega _1=1$, $\omega _2=\omega$, $\epsilon ^{(1)}=\epsilon$.
Then we have
\begin{equation*}
\begin{aligned}
& D=\bigg\{ z \in \mathfrak H^2 \mid
\epsilon ^{-2}\leq\frac {y_2}{y_1} \leq \epsilon ^2, \\
& -\frac {1}{2} \leq \frac {1}{\omega -\omega ^{\prime}}
(\omega ^{\prime}x_1-\omega x_2) \leq \frac {1}{2}, \quad
-\frac {1}{2} \leq \frac {1}{\omega -\omega ^{\prime}}
(x_1-x_2) \leq \frac {1}{2}, \\
& N(\vert cz+d\vert )\geq 1 \
\text {whenever $c$ and $d$ are relatively prime integers of $\mathcal O_F$} \bigg\} .
\end{aligned}
\tag{A.8}
\end{equation*}
Here $\omega ^{\prime}$ denotes the conjugate of $\omega$.

Hereafter in this section, we assume that $F=\mathbf Q(\sqrt {5})$.
We take $\omega =\epsilon$. The next lemma is the essential ingredient
of the proof of Theorem A.1.

\proclaim Lemma A.11. Let $F=\mathbf Q(\sqrt {5})$ and take $\omega =\epsilon$.
Put $S=\{ \gamma \in \Gamma \mid D \cap \gamma D \neq \emptyset \}$.
Then $S$ is a finite set and we have $S\subset S_0\sqcup S_1\sqcup S_2$, where
$$
S_0=P, \qquad
S_1=\{ \gamma =\begin{pmatrix} a & b \\ c & d \end{pmatrix} , \ c \in E_F \},
$$
$$
S_2=\left\{ \gamma =
\begin{pmatrix} \pm\epsilon ^3 & b \\ 2\epsilon & \pm\epsilon ^3 \end{pmatrix} , \qquad
\begin{pmatrix} \pm 1 & b \\ 2\epsilon ^{-2}&\pm 1 \end{pmatrix} \right\} .
$$
Here $\pm$ can be taken arbirarily and
$b \in \mathcal O_F$ is chosen so that $\det \gamma =1$.
($S_2$ consists of eight elements.)
\par

We give a proof of Theorem A.1 assuming Lemma A.11.

\medskip

{\bf Proof of Theorem A.1.} We consider
$\mathfrak H^2 \subset \mathbf C^2$ and let $d$ denote the Euclidean metric
induced by this embedding. For $\delta >0$, we put
\begin{equation*}
D_{\delta}=\{ z \in \mathfrak H^2 \mid d(z, D)<\delta \} .
\end{equation*}
We see easily that $D$ is path connected. Let $z \in D_{\delta}$.
Then there exists $z_1 \in D$ such that $d(z, z_1)<\delta$.
Hence $z$ is connected by a path to $z_1$. Therefore $D_{\delta}$ is
path connected. By using the argument of Lemma A.10, we see that
$\cap _{\delta >0} S(D_{\delta})=S$.
Moreover we can show without difficulty that $S(D_{\delta})$ is finite
when $\delta$ is sufficiently small. 
Therefore $S(D_{\delta})=S$ when $\delta$ is sufficiently small
and Theorem M can be applied with $S$ given in Lemma A.11.

For $\gamma \in S$, we prepare a symbol $[\gamma]$ and consider
the free group $\mathcal F^{\prime}$ on the free generators $[\gamma ]$.
By Theorem M, it is sufficient to show that
$[\gamma _2]^{-1}[\gamma _1]^{-1}[\gamma _1\gamma _2]$, $\gamma _1$, $\gamma _2$,
$\gamma _1\gamma _2 \in S$ can be reduced to a three term relation.
We put $S_i^{\prime}=S \cap S_i$, $0 \leq i \leq 2$.
We can check easily that $\sigma$, $\mu$, $\tau$, $\eta \in S$.
Hence Lemma A.9 is applicable. Let $\mathcal F$ be the free group on the free generators
$\widetilde\sigma$, $\widetilde\mu$, $\widetilde\tau$ and $\widetilde\eta$.
We define a homomorphism $\pi : \mathcal F \longrightarrow \Gamma$ by
$\pi (\widetilde\sigma )=\sigma$, $\pi (\widetilde\mu )=\mu$,
$\pi (\widetilde\tau )=\tau$, $\pi (\widetilde\eta )=\eta$.
For $\gamma \in S$, we define $\widetilde\gamma \in \mathcal F$ such that
$\pi (\widetilde\gamma )=\gamma$ as follows.

If $\gamma \in P$, we write $\gamma =\mu ^a\tau ^b\eta ^c$. Then we define
$\widetilde\gamma =\widetilde\mu ^a\widetilde\tau ^b\widetilde\eta ^c$.
In particular, this rule applies to an element $\gamma \in S_0^{\prime}$.
We have
\begin{equation}
\begin{pmatrix} a & b \\ c & d \end{pmatrix}
=\begin{pmatrix} 1 & c^{-1}a \\ 0 & 1 \end{pmatrix}
\begin{pmatrix} 0 & 1 \\ -1 & 0 \end{pmatrix}
\begin{pmatrix} -c & -d \\ 0 & -c^{-1} \end{pmatrix} ,
\qquad c \in E_F.
\tag{A.9}
\end{equation}
Hence $\gamma \in S_1^{\prime}$ can be written as
$\gamma =p_1\sigma p_2$, $p_1$, $p_2 \in P$. We fix such an expression and define
$\widetilde\gamma =\widetilde p_1\widetilde\sigma \widetilde p_2$.
Suppose $\gamma \in S_2^{\prime}$. We write $\gamma$ in the form
$\gamma =\begin{pmatrix} u & \beta \\ 2\epsilon ^m & u^{\ast} \end{pmatrix}$,
$u$, $u^{\ast} \in E_F$, $\beta \in \mathcal O_F$, $m \in \mathbf Z$.
We have
\begin{equation}
\begin{pmatrix} u & \beta \\ 2\epsilon ^m & u^{\ast} \end{pmatrix}
=\sigma\begin{pmatrix} 1 & -2u^{-1}\epsilon ^m \\ 0 & 1 \end{pmatrix}
\sigma\begin{pmatrix} -u & -\beta \\ 0 & -u^{-1} \end{pmatrix} .
\tag{A.10}
\end{equation}
We fix this expression $\gamma =\sigma p_1\sigma p_2$, $p_1$, $p_2 \in P$
and define
$\widetilde\gamma =\widetilde\sigma\widetilde p_1\widetilde\sigma \widetilde p_2$.

By Lemma A.9, it is sufficient to show that
$\widetilde\gamma _2^{-1}\widetilde\gamma _1^{-1}\widetilde{\gamma _1\gamma _2}$
reduces to a three terms relation (under (i) $\sim$ (vi) and (A.2)) when
$\gamma _1$, $\gamma _2$, $\gamma _1\gamma _2 \in S$.
We see that there cannot arise the case where
all of $\gamma _1$, $\gamma _2$, $\gamma _1\gamma _2$ belong to $S_2^{\prime}$,
by inspecting the $(2,1)$-component of $\gamma _1\gamma _2$.
This implies that if two of $\gamma _1$, $\gamma _2$, $\gamma _1\gamma _2$
belong to $S_2^{\prime}$, then the other one must belong to $S_0^{\prime}$.
Therefore $\widetilde\gamma _2^{-1}\widetilde\gamma _1^{-1}\widetilde{\gamma _1\gamma _2}$
defines at most a four terms relation. We may assume that
$\widetilde\gamma _2^{-1}\widetilde\gamma _1^{-1}\widetilde{\gamma _1\gamma _2}$
defines a four terms relation. Then one of $\gamma _1$, $\gamma _2$, $\gamma _1\gamma _2$
belongs to $S_2^{\prime}$. By (A.10), this relation takes the form
(A.3$^{\prime}$) with $x \in \mathcal O_F$ such that $(x)=(2)$.
As shown in the proof of Lemma A.5, it suffices to consider the four terms relation
$\{ x, u \}$ for $u \in E_F$ such that $x$ divides $u-1$ . Now the group
$E_{(2)}=\{ u \in E_F \mid u \equiv 1 \mod 2 \}$ is generated by
$-1$ and $\epsilon ^3$. By $\epsilon ^3-1=2\epsilon$ and Remark A.7,
we see that $\{ x, \epsilon ^{3e} \}$ is reducible to (i) $\sim$ (vi) and (A.2)
for $e \in \mathbf Z$. By Lemma A.6, (5), $\{ x, -\epsilon ^{3e} \}$ is reducible to
(i) $\sim$ (vi) and (A.2). This completes the proof.

\medskip

Now we are going to prove Lemma A.11.
We consider an element $\gamma \in \Gamma$ such that
for a point $z \in D$, $\gamma z \in D$ holds, i.e.,
$D \cap \gamma ^{-1}D \neq \emptyset$.
\footnote{Since $S_i$, $i=0$, $1$, $2$ is stable under $\gamma \mapsto\gamma ^{-1}$,
it suffices to determine $\gamma$ which satisfies $D \cap \gamma ^{-1}D \neq \emptyset$.}.
We put $\gamma =\begin{pmatrix} a & b \\ c & d \end{pmatrix}$,
$z^{\prime}=\gamma z$, $z^{\prime}=(z_1^{\prime}, z_2^{\prime})$,
$z_j^{\prime}=x_j^{\prime}+iy_j^{\prime}$, $j=1$, $2$,
$y^{\prime}=(y_1^{\prime}, y_2^{\prime})$. We have
$$
N(y^{\prime})=\frac {N(y)}{N(\vert cz+d\vert )^2}.
$$
Hence $N(y^{\prime})\leq N(y)$. Changing the roles of $z$ and $z^{\prime}$,
we have $N(y)\leq N(y^{\prime})$. Hence we see that
$N(y^{\prime})=N(y)$ and
\begin{equation}
N(\vert cz+d\vert )=1.
\tag{A.11}
\end{equation}
Since we are assuming that $F=\mathbf Q(\sqrt {5})$, $\omega =\epsilon$, we have
$$
x_1=X_1+\frac {1+\sqrt {5}}{2}X_2, \quad
x_2=X_1+\frac {1-\sqrt {5}}{2}X_2, \quad
-\frac {1}{2} \leq X_1 \leq \frac {1}{2}, \ -\frac {1}{2} \leq X_2 \leq \frac {1}{2}.
$$
Then $x_1x_2=X_1^2-X_2^2+X_1X_2$ and we see that
\begin{equation}
\vert x_1x_2 \vert \leq \frac {5}{16}, \qquad
\vert x_1\vert \leq \frac {3+\sqrt {5}}{4}, \qquad \vert x_2\vert \leq \frac {1+\sqrt {5}}{4}.
\tag{A.12}
\end{equation}
Since $z \in D$, we have
\begin{equation}
N(\vert z\vert )^2=(x_1^2+y_1^2)(x_2^2+y_2^2) \geq 1.
\tag{A.13}
\end{equation}
Put $k=y_1y_2$. Since $\epsilon ^{-2}\leq y_1/y_2\leq \epsilon ^2$, we have
$\epsilon ^{-1}\sqrt {k} \leq y_1$, $y_2 \leq \epsilon\sqrt {k}$. Then by (A.13), we have
$$
k^2+(x_1^2+x_2^2)\epsilon ^2k +x_1^2x_2^2-1 \geq 0.
$$
We consider the equation with respect to $t$:
\begin{equation}
t^2+(x_1^2+x_2^2)\epsilon ^2t +x_1^2x_2^2-1=0.
\tag{A.14}
\end{equation}
Let $\xi$ be the positive root of (A.14) and let
$\kappa ^{\ast}=\min \xi$. Here the minimum is taken with respect to
$X_1$ and $X_2$, regarding $x_1$ and $x_2$ as the functions of $X_1$ and $X_2$;
$X_1$ and $X_2$ extend over the domain $-1/2 \leq X_1, X_2 \leq 1/2$.
Let $\kappa$ be the positive root of the equation
$$
t^2+\frac {7(3+\sqrt{5})}{8}t-\frac {15}{16}=0.
$$
This is the positive root of (A.14) when $X_1=X_2=1/2$,
$x_1=(3+\sqrt{5})/4$, $x_2=(3-\sqrt {5})/4$.
We have $\kappa =0.19622\cdots$. By elementary but somewhat tedious calculation,
which we omit the details, we can show that $\kappa ^{\ast}= \kappa$.
Hence we have
\begin{equation}
y_1y_2 \geq \kappa=0.19622\cdots .
\tag{A.15}
\end{equation}

If $c=0$, then $\gamma \in S_0$. It suffices to show that
$\gamma \in S_1 \sqcup S_2$ assuming $c \neq 0$.
By (A.11), we have
\begin{equation}
\vert N(c)\vert y_1y_2 \leq 1.
\tag{A.16}
\end{equation}
By (A.15), we have $\vert N(c)\vert \leq 1/\kappa$. Therefore
$\vert N(c)\vert =1$ or $4$ or $5$. If $\vert N(c)\vert =1$, then
$c \in E_F$ and $\gamma \in S_1$. Hereafter we assume $\vert N(c)\vert =4$ or $5$.
By (A.15) and (A.16), noting $\epsilon ^{-2}\leq y_1/y_2\leq \epsilon ^2$, , we obtain
\begin{equation}
\epsilon ^{-1}\sqrt {\kappa} \leq y_1, y_2 \leq \frac {\epsilon}{\sqrt {\vert N(c)\vert}} .
\tag{A.17}
\end{equation}
Since $N(\vert z\vert )\geq 1$, we have $(x_1^2+y_1^2)(x_2^2+y_2^2) \geq 1$.
Using $y_1y_2 \leq 1/\vert N(c)\vert$, we have
\begin{equation}
x_1^2y_2^4-(1-x_1^2x_2^2-\frac {1}{N(c)^2})y_2^2+\frac {x_2^2}{N(c)^2} \geq 0.
\tag{A.18}
\end{equation}
If $x_1=0$, we obtain
$$
y_1^2x_2^2 \geq 1-\frac {1}{N(c)^2} \geq 1-\frac {1}{16}
$$
from $N(\vert z\vert )\geq 1$ and (A.16). By (A.12), we have
$$
y_1 \geq \sqrt {1-\frac {1}{16}}\cdot \frac {2}{\epsilon}=1.19681\cdots .
$$
This contradics (A.17). Hence we have $x_1 \neq 0$.

First we exclude the case $\vert N(c)\vert =5$. To this end,
we assume $\vert N(c)\vert =5$ and consider the equation (cf. (A.18))
\begin{equation}
x_1^2t^2-(1-x_1^2x_2^2-\frac {1}{25})t+\frac {x_2^2}{25}=0.
\tag{A.19}
\end{equation}
Let $f(t)$ be the polynomial of $t$ on the left-hand side. For
$t_0=\epsilon ^{-2}\kappa$, we have
$$
f(t_0)\leq (\frac {\epsilon +1}{2})^2t_0^2-(1-\frac {25}{256}-\frac {1}{25})t_0
+\frac {1}{25}(\frac {\epsilon}{2})^2=-0.02882\cdots <0
$$
using (A.12). Let $\eta _1>\epsilon ^{-2}\kappa >\eta _2$ be the roots of the equation (A.19).
By (A.17) and (A.18), we must have $y_2 \geq \sqrt {\eta _1}$.
We note that (cf. (A.17))
\begin{equation}
y_1, y_2 \leq \frac{\epsilon}{\sqrt{5}}=0.72360\cdots .
\tag{A.20}
\end{equation}
We consider $\eta _1$ as a function of $X_1$ and $X_2$ defined in the domain
$-1/2 \leq X_1, X_2 \leq 1/2$.
First we consider $\eta _1$ on the subdomain defined by the condition $x_1>0$.
It is not difficult to check that $\eta _1$ is monotone decreasing
with respect to the both arguments $X_1$ and $X_2$.
For $X_1=1/2$, $X_2=0.4985$, we have $\sqrt {\eta _1}=0.72377\cdots$.
For $X_1=0.4985$, $X_2=1/2$, we have $\sqrt {\eta _1}=0.72389\cdots$.
In view of (A.20), we must have $X_1$, $X_2>0.4985$.
Similarly, in the subdomain $x_1<0$, we must have $X_1$, $X_2<-0.4985$.

First we consider the case $X_1$, $X_2>0.4985$. For relatively prime integers
$\alpha$, $\beta \in \mathcal O_F$, we have (cf. (A.8))
$N(\vert \alpha z+\beta \vert ) \geq 1$. Take $\alpha =2$, $\beta =-\epsilon ^2$.
We have
$$
\vert 2x_1-\epsilon ^2\vert \leq 0.03(1+\epsilon ), \qquad
\vert 2x_2-\epsilon ^{-2}\vert \leq 0.03(1+\vert \epsilon ^{\prime}\vert ).
$$
Here $\epsilon ^{\prime}=(1-\sqrt {5})/2$ is the conjugate of $\epsilon$.
Then we find
\begin{equation*}
\begin{aligned}
& N(\vert 2z-\epsilon ^2\vert )^2
= \{ (2x_1-\epsilon ^2)^2+4y_1^2 \}\{ (2x_2-\epsilon ^{-2})^2+4y_2^2 \} \\
= & \, 16y_1^2y_2^2+4y_1^2(2x_2-\epsilon ^{-2})^2
+4y_2^2(2x_1-\epsilon ^2)^2+(2x_1-\epsilon ^2)^2(2x_2-\epsilon ^{-2})^2 \\
\leq & \, \frac {16}{25} +4y_1^2\{ 0.03(1+\epsilon )\} ^2
+4y_2^2 \{ 0.03(1+\vert \epsilon ^{\prime}\vert )\} ^2 \\
& \hskip 5em +\{ 0.03(1+\epsilon )\} ^2\{ 0.03(1+\vert \epsilon ^{\prime}\vert )\} ^2.
\end{aligned}
\end{equation*}
Since $y_1$, $y_2 \leq 0.72360\cdots$, this contradicts
$N(\vert 2z-\epsilon ^2\vert ) \geq 1$.
When $X_1$, $X_2 <-0.4985$, we obtain a contradiction similarly by taking
$\alpha =2$, $\beta =\epsilon ^2$. Thus we have shown that the case
$\vert N(c)\vert =5$ cannot occur.

It remains to show that $\gamma \in S_2$ assuming $\vert N(c)\vert =4$.
We can write $c=\pm 2\epsilon ^m$ with $m \in \mathbf Z$.
Changing $\gamma$ to $-\gamma$ if necessary, we may assume that $c=2\epsilon ^m$.
We put $z^{\prime}=(z_1^{\prime}, z_2^{\prime})=\gamma z$,
$z_j^{\prime}=x_j^{\prime}+iy_j^{\prime}$, $j=1$, $2$. Since $z=\gamma ^{-1}z^{\prime}$,
$\gamma ^{-1}=\begin{pmatrix} d & -b \\ -c & a \end{pmatrix}$, the estimate (A.17) holds
also for $y_1^{\prime}$ and $y_2^{\prime}$. We have
\begin{equation}
\epsilon ^{-1}\sqrt{\kappa} =0.27376\cdots \leq y_1, y_2, y_1^{\prime}, y_2^{\prime}
\leq \frac {\epsilon}{\sqrt{\vert N(c)\vert}}=0.80901\cdots .
\tag{A.21}
\end{equation}
We have
$$
\vert c^{(j)}z_j+d^{(j)}\vert ^2=\frac {y_j}{y_j^{\prime}}, \qquad j=1, 2.
$$
Hence we obtain
\begin{equation}
\epsilon ^{-2}\sqrt {\kappa}\sqrt {\vert N(c)\vert} \leq \vert c^{(j)}z_j+d^{(j)}\vert ^2
\leq \frac {\epsilon ^2}{\sqrt {\kappa}\sqrt {\vert N(c)\vert}}, \qquad j=1, 2.
\tag{A.22}
\end{equation}
In particular, we have
$$
(c^{(j)})^2y_j^2 \leq \frac {\epsilon ^2}{\sqrt {\kappa}\sqrt {\vert N(c)\vert}}, \qquad j=1, 2.
$$
Using (A.21), we obtain
\begin{equation}
\vert c^{(j)}\vert \leq \epsilon ^2\kappa ^{-3/4}\vert N(c)\vert ^{-1/4}=6.27915\cdots ,
\qquad j=1, 2.
\tag{A.23}
\end{equation}
From (A.23), we obtain $m=0$, $\pm 1$, $\pm 2$.

Next we are going to restrict possibilities of $d$. A preliminary table of listing
all possible $d$ can be obtained by (A.22) and (A.23).
By (A.11) and (A.21), we have
\begin{equation}
\{ (2\epsilon ^mx_1+d^{(1)})^2+4\epsilon ^{2m}\cdot\epsilon ^{-2}\kappa \}
\{ (2(\epsilon ^{\prime})^m x_2+d^{(2)})^2+4\epsilon ^{-2m}\cdot \epsilon ^{-2}\kappa \} \leq 1.
\tag{A.24}
\end{equation}
We consider the equation (cf. (A.18))
\begin{equation}
x_1^2t^2-(1-x_1^2x_2^2-\frac {1}{16})t+\frac {x_2^2}{16}=0.
\tag{A.25}
\end{equation}
Let $g(t)$ be the polynomial of $t$ on the left-hand side. For
$t_0=\epsilon ^{-2}\kappa$, we can check $g(t_0)<0$.
Let $\eta _1>t_0>\eta _2$ be the roots of $g(t)$. By (A.18) and (A.21), we have
$y_2 \geq \sqrt {\eta _1}$. As in the case where $\vert N(c)\vert =5$,
we consider $\eta _1$ as a function of $X_1$ and $X_2$ defined in the domain
$-1/2 \leq X_1, X_2 \leq 1/2$.
On the subdomain defined by the condition $x_1>0$,
we check that $\eta _1$ is monotone decreasing
with respect to the both arguments $X_1$ and $X_2$.
For $X_1=1/2$, $X_2=0.39$, we have $\sqrt {\eta _1}=0.81291\cdots$.
For $X_1=0.38$, $X_2=1/2$, we have $\sqrt {\eta _1}=0.81101\cdots$.
In view of (A.21), we must have $X_1>0.38$, $X_2>0.39$.
Similarly, in the subdomain $x_1<0$, we must have $X_1<-0.38$, $X_2<-0.39$.
Let $V$ be the closed domain
$$
V=\{ (X_1, X_2) \mid 0.38 \leq \vert X_1\vert \leq 1/2, \
0.39 \leq \vert X_2\vert \leq 1/2 \}
$$
and consider the function
$$
f(X_1, X_2)=
\{ (2\epsilon ^mx_1+d^{(1)})^2+4\epsilon ^{2m-2}\kappa \}
\{ (2(\epsilon ^{\prime})^m x_2+d^{(2)})^2+4\epsilon ^{-2m-2}\kappa \}
$$
on $V$. By (A.24), we see that:
\begin{equation}
\text {The minimum of $f(X_1, X_2)$ on $V$ does not exceed $1$.}
\tag{C1}
\end{equation}
Next let $\xi$ be the positive root of (A.14). Since $y_1y_2 \geq \xi$, we have
$y_1$, $y_2 \geq \epsilon ^{-1}\sqrt {\xi}$. By (A.11), we obtain another inequality:
\begin{equation*}
\begin{aligned}
& (2\epsilon ^mx_1+d^{(1)})^2(2(\epsilon ^{\prime})^m x_2+d^{(2)})^2
+4\epsilon ^{-2m-2}\xi (2\epsilon ^mx_1+d^{(1)})^2 \\
+ & \, 4\epsilon ^{2m-2}\xi (2(\epsilon ^{\prime})^mx_2+d^{(2)})^2
+16\xi ^2 \leq 1.
\end{aligned}
\tag{A.26}
\end{equation*}
We regard $x_1$, $x_2$ and $\xi$ as the functions of $X_1$ and $X_2$ and
let $g(X_1, X_2)$ be the function on the left-hand side of (A.26).
Then (A.26) implies:
\begin{equation}
\text {The minimum of $g(X_1, X_2)$ on $V$ does not exceed $1$.}
\tag{C2}
\end{equation}
By numerical computations using a computer, we find the following:

For $m=0$, (C1) leaves possibilities $d=\pm 1$, $\pm \epsilon$, $\pm \epsilon ^2$,
$\pm \epsilon ^{-1}$. If combined with (C2), the only possibility is $d=\pm \epsilon ^2$.
For $m=1$, (C1) leaves possibilities $d=\pm 1$, $\pm \epsilon$, $\pm \epsilon ^2$,
$\pm \epsilon ^3$. If combined with (C2), the only possibility is $d=\pm \epsilon ^3$.
For $m=2$, (C1) leaves possibilities $d=\pm \epsilon$, $\pm \epsilon ^2$, $\pm \epsilon ^3$, 
$\pm \epsilon ^4$. If combined with (C2), the only possibility is $d=\pm \epsilon ^4$.
For $m=-1$, (C1) leaves possibilities $d=\pm 1$, $\pm \epsilon$, $\pm \epsilon ^{-1}$,
$\pm \epsilon ^{-2}$. If combined with (C2), the only possibility is $d=\pm \epsilon$.
For $m=-2$, (C1) leaves possibilities $d=\pm 1$, $\pm \epsilon ^{-1}$, $\pm \epsilon ^{-2}$,
$\pm \epsilon ^{-3}$. If combined with (C2), the only possibility is $d=\pm 1$.

Thus, in every case where $c=2\epsilon ^m$, we have $d=\pm \epsilon ^n$
with $n$ depending only on $m$.
Changing the roles of $z$ and $z^{\prime}$ and noting that
$-\gamma ^{-1}=\begin{pmatrix} -d & b \\ c & -a \end{pmatrix}$,
we see that $a$ must have the same form $a=\pm \epsilon ^n$.
(Here the $\pm$ sign is arbitrary but $n$ is the same for $d$ and $a$.)
By $\det \gamma =1$, we have $ad \equiv 1 \mod 2$,
which implies $n \equiv 0 \mod 3$. Therefore only the cases
$m=1$, $-2$ can survive and we see that $\gamma \in S_2$.
This completes the proof of Lemma A.11.

\noindent Department of Mathematics, Kyoto University, Kyoto 606-8502, Japan

\noindent E-mail: yoshida@math.kyoto-u.ac.jp


\bigskip

\begin{thebibliography}{A}

\bibitem [B]{B} D. Blasius, Hilbert modular forms and the Ramanujan conjecture,
Noncommutative geometry and number theory, 35--56, Aspects Math., 37, Vieweg, 2006.

\bibitem [BW]{BW} A. Borel and N. Wallach, Continuous cohomology,
discrete subgroups, and representations of reductive groups, Ann. Math. Studies 94,
Princeton University Press, 1980.

\bibitem [CE]{CE} H. Cartan and S. Eilenberg, Homological algebra,
Princeton University Press, 1956.

\bibitem [DG]{DG} K. Doi and K. Goto, On the $L$-series associated with
modular forms, Memoirs of Institute of Science and Engineering,
Ritsumeikan Univ., 52 (1993), 1--19 (in Japanese).

\bibitem [DHI]{DHI} K. Doi, H. Hida and H. Ishii, Discriminant of Hecke fields
and twisted adjoint $L$-values for $GL(2)$, Inv. Math. 134 (1998), 547--577.

\bibitem [DI]{DI} K. Doi and H. Ishii, Hilbert modular $L$-values and discriminant of
Hecke's fields, Memoirs of Institute of Science and Engineering,
Ritsumeikan Univ., 53 (1994), 1--12.

\bibitem [E]{E} B. Eckman, Cohomology of groups and tranfer,
Ann. of Math. 58 (1953), 481--493.

\bibitem [Ha]{Ha} G. Harder, Eisenstein cohomology of arithmetic groups.
The case ${\rm GL}_2$, Inv. Math. 89 (1987), 37--118.

\bibitem [Hi1]{Hi1} H. Hida, On abelian varieties with complex multiplication
as factors of the abelian variety attached to Hilbert modular forms,
Japanese J. Math. 5 (1979), 157--208.

\bibitem [Hi2]{Hi2} H. Hida, $p$-ordinary cohomology groups for
${\rm SL}(2)$ over number fields, Duke Math. J.,
69 (1993), 259--314.

\bibitem [Hi3]{Hi3} H. Hida, On the critical values of $L$-functions of
${\rm GL}(2)$ and ${\rm GL}(2)\times{\rm GL}(2)$, Duke Math. J.,
74 (1994), 431--529.

\bibitem [HW]{HW} G. H. Hardy and E. M. Wright, An introduction to
the theory of numbers, Oxford University Press, fifth edition, 1979.

\bibitem[JL]{JL} H. Jacquet and R. P. Langlands, Automorphic forms on $GL(2)$,
Lecture notes in mathematics 114, Springer-Verlag, 1970.

\bibitem[K]{K} A. G. Kurosh, The theory of groups, English edition, two volumes, Chelsea,
1955, 1956.

\bibitem[KS]{KS} M. Kuga and G. Shimura, On vector differential forms attached to
automorphic forms, J. Math. Soc. Japan, 12 (1960), 258--270
 (= Collected Papers of Goro Shimura I, [60a]).

\bibitem [Mac]{Mac} A. M. Macbeath, Groups of homeomorphisms of a simply connected space,
Ann. of Math. 79 (1964), 473--488.

\bibitem[Man]{Man} Y. I. Manin, Periods of parabolic forms and $p$-adic Hecke series,
Math. USSR Sbornik 21 (1973), 371--393.

\bibitem[MM]{MM} Y. Matsushima and S. Murakami, On vector valued harmonic forms
and automorphic forms on symmetic Riemannian manifolds,
Ann. of Math. 78 (1963), 365--416.

\bibitem[MS]{MS} Y. Matsushima and G. Shimura, On the cohomology groups
attached to certain vector valued differential forms on the product of the upper half plane,
Ann. of Math. 78 (1963), 417--449 (= Collected Papers of Goro Shimura I, [63c]).

\bibitem[PARI2]{PARI2} PARI/GP, version {\tt 2.3.4}, Bordeaux, 2008,
{http://pari.math.u-bordeaux.fr/}.

\bibitem [Sc]{Sc} O. Schreier, Die Untergruppen der freie Gruppen, Abh. Math.
Sem. Univ. Hamburg 5 (1927), 161--183.

\bibitem [Se1]{Se1} J-P. Serre, Corps locaux, deuxi$\grave {\rm e}$me \'edition, Hermann, 1968.

\bibitem [Se2]{Se1} J-P. Serre, Cohomologie des groupes discrets,
Ann. of Math. Studies 70 (1971), 77--169 (=\OE{}uvre II, 88).

\bibitem [Shi]{Shi} H. Shimizu, On discontinuous groups operating on the product of
the upper half planes, Ann. of Math. 77 (1963), 33--71.

\bibitem [Sh1]{Sh1} G. Shimura, Sur les int\'egrales attach\'ees aux formes automorphes,
J. Math. Soc. Japan 11 (1959), 291--311 (= Collected Papers I, [59c]).

\bibitem [Sh2]{Sh2} G. Shimura, Introduction to the Arithmetic Theory of
Automorphic Functions, Iwanami Shoten and Princeton University Press, 1971.

\bibitem [Sh3]{Sh3} G. Shimura, The special values of the zeta functions associated with
cusp forms, Comm. pure and applied Math. 29 (1976), 783--804 (=Collected Papers II, [76b]).

\bibitem [Sh4]{Sh4} G. Shimura, The special values of the zeta functions associated with
Hilbert modular forms, Duke Math. J. 45 (1978), 637--679 (=Collected Papers III, [78c]).

\bibitem [Sh5]{Sh5} G. Shimura, The critical values of certain Dirichlet series attached to
Hilbert modular forms, Duke Math. J. 63 (1991), 557--613 (=Collected Papers IV, [91]).

\bibitem [Sh6]{Sh6} G. Shimura, Eisenstein series and zeta functions on
symplectic groups, Inv. Math. 119 (1995), 539--584 (=Collected Papers IV, [95a]).

\bibitem [Sh7]{Sh7} G. Shimura, Arithmeticity in the theory of automorphic forms,
Math. Surveys and Monogr. vol. 82, American Mathematical Society, 2000.

\bibitem [Si1]{Si1} C. L. Siegel, Discontinuous groups, Ann. of Math. 44(1943), 674--689
(= Gesammelte Abhandlungen III, No. 43).

\bibitem [Si2]{Si2} C. L. Siegel, Lectures on advanced analytic number theory,
Tata Institute, 1961.

\bibitem [Su]{Su} M. Suzuki, Group Theory I, Grundlehren der mathematischen
Wissenshaften 247, Springer Verlag, 1982.

\bibitem [Sw]{Sw} R. W. Swan, Generators and relations for certain special linear groups,
Advances in Math. 6 (1971), 1--77.

\bibitem [V]{V} L. N. Vaser$\rm {\check{s}}$tein,
On the group $SL_2$ over Dedekind rings of arithmetic type,
Math. USSR Sbornik 18 (1972), 321--332.

\bibitem [Y1]{Y1} H. Yoshida, On the zeta functions of Shimura varieties and periods
of Hilbert modular forms, Duke Math. J. 75 (1994), 121--191.

\bibitem [Y2]{Y2} H. Yoshida, On a conjecture of Shimura concerning periods of Hilbert
modular forms, Amer. J. Math. 117 (1995), 1019--1038.

\bibitem [Y3]{Y3} H. Yoshida, Absolute CM-periods, Math. Surveys and Monogr. vol. 106,
American Mathematical Society, 2003.

\end{thebibliography}
\end{document}